Fine Resolution of $k$-transversal Cones

M. Stiefenhofer

ABSTRACT. Tougeron's implicit function theorem and Hensel's lemma are well known representatives concerning $2k$-approximation/$k$-nondegeneracy implying existence of solutions with identity of order $k$. This note aims to extend this principle to equations $G[z] = 0$ in Banach spaces, using $k$-transversality concepts, which may geometrically be interpreted as generalized cones spanned by submanifolds, each characterized by a certain expansion rate.

The number of manifolds in the cone, as well as their expansion rates, is recursively increased until an appropriate desingularization of the cone is build up with linearization expressed by first $k + 1$ derivatives of the singular operator at the base point.

Along these lines, a well-defined submersion is constructed in the cone with uniformly bounded inverse when approaching the singularity. The techniques are restricted to curves, possibly touching by high order the singular locus of $G$, but ultimately traversing it, in this way defining an isolated singularity of the operator family given by the linearization along the curve.

The fine resolution of the cone by the manifolds represents an improvement compared to measuring the variation of the nonlinear operator exclusively by the overall behaviour of the determinant.

In case of finite dimensions, each half-cone is characterized by a constant topological degree that can be used to investigate a solution curve in general position with respect to secondary bifurcation.

The core of all considerations is given by some characteristic patterns, valid in the system of undetermined coefficients that allow for detailed analysis of the power series resulting from plugging the power series of the ansatz into the power series of the nonlinear operator.

*Contents*

## 1. Introduction

Given an equation $G[x,y]=0,\ G\in C^q(\mathbb{K}^n\times\mathbb{K}^m,\mathbb{K}^m),\ \mathbb{K}=\mathbb{R},\mathbb{C}$, with base solution $G[0,0]=0$, which is locally embedded within a smooth approximation $y_0(x)$ satisfying

$$\| G[x,y_0(x)] \| = O\big( \|x\|^{2k+1}\big) \qquad (A_0)$$

$$det\big\{ G_y[x,y_0(x)] \big\} \neq O\big( \|x\|^{k+1}\big), \qquad (N_0)$$

i.e. along $y_0(x)$ the map value $G[x,y_0(x)]$ varies slowly away from zero by order of $2k+1$, compared to the fast change of the $y$-derivative $G_y[x,y_0(x)]$ that varies by order of $k$ (and not slower). Under these assumptions, it is well known that further zeros $y(x)$ of $G[x,y]$ exist in the vicinity of the approximation $y_0(x)$, i.e. we obtain $G[x,y(x)]=0$ and identity of order $k$ between exact solution curve $y(x)$ and approximation $y_0(x)$ by

$$\| y(x)-y_0(x) \| = O\big( \|x\|^{k+1}\big). \qquad (I_0)$$

Note that the approximation condition $(A_0)$ defines a degeneracy condition of order $2k$, i.e. we require the first $2k$ derivatives of $G[x,y_0(x)]$ to vanish at $x=0$, whereas $(N_0)$ defines a nondegeneracy condition of order $k$, finally allowing a contraction mapping argument to ascertain the existence of $y(x)$ as well as the identity condition $(I_0)$ of order $k$. The proof is essentially based on the ansatz

$$y = y_0(x) + det\big\{ G_y[x,\ y_0(x)] \big\}\cdot \bar{y}\,, \qquad (1.1)$$

as well as factoring out of the adjoint matrix of the Jacobian $G_y[\cdot]$ yielding a remainder equation appropriate for application of the implicit function theorem with respect to the new variable $\bar{y}$.

Formulas $(A_0),(N_0),(I_0)$ may be summarized as follows: $2k$-approximation combined with $k$-nondegeneracy implies $k$-identity between existing solutions and given approximation. A constellation of this kind occurs in several places of algebra and analysis, applied to mappings of quite different properties with respect to smoothness and type of participating spaces. Compare Newton's Lemma [G], Tougeron's implicit function theorem [BH], [H], [R] or Hensel's Lemma [F]. Nevertheless, the techniques of proof are quite similar.

Now in this note, the constellation is applied to $G[z]=0,\ G\in C^q(B,\bar{B})$ with $G[0]=0$ and $B,\bar{B}$ real or complex Banach spaces. The linearization $G'[0]\in L[B,\bar{B}]$ is assumed to be bounded, but not necessarily a Fredholm operator, i.e. in general it is not supposed that reduction to finite dimensions by Lyapunov-Schmidt is possible.

When moving from the finite dimensional setting of $(A_0),(N_0),(I_0)$ to general Banach spaces, then the determinant and the adjoint matrix of the Jacobian are no longer present in a straightforward way, implying the nondegeneracy condition $(N_0)$ and the proof by itself to be reformulated appropriately. In some more detail, the proof will be characterized by replacing (1.1) by

$$z = z_0(\varepsilon) + A_\varepsilon\cdot n^c \qquad (1.2)$$

with $z_0(\varepsilon)$ an approximation of order $2k$, now parametrized by the external parameter $\varepsilon\in\mathbb{R}$ or $\varepsilon\in\mathbb{C}$ and $A_\varepsilon\in L[B,B]$ an $\varepsilon$-dependent family of linear operators acting between appropriate subspaces of $B$. The action of $A_\varepsilon$ may be interpreted as a blow-up coordinate transformation of $B$ with singularity given at $\varepsilon=0$, due to $A_0=0\in L[B,B]$.

In some more detail, the subspaces are successively built up by a resolution of the two Banach spaces $B$ and $\bar{B}$ according to the following diagram.

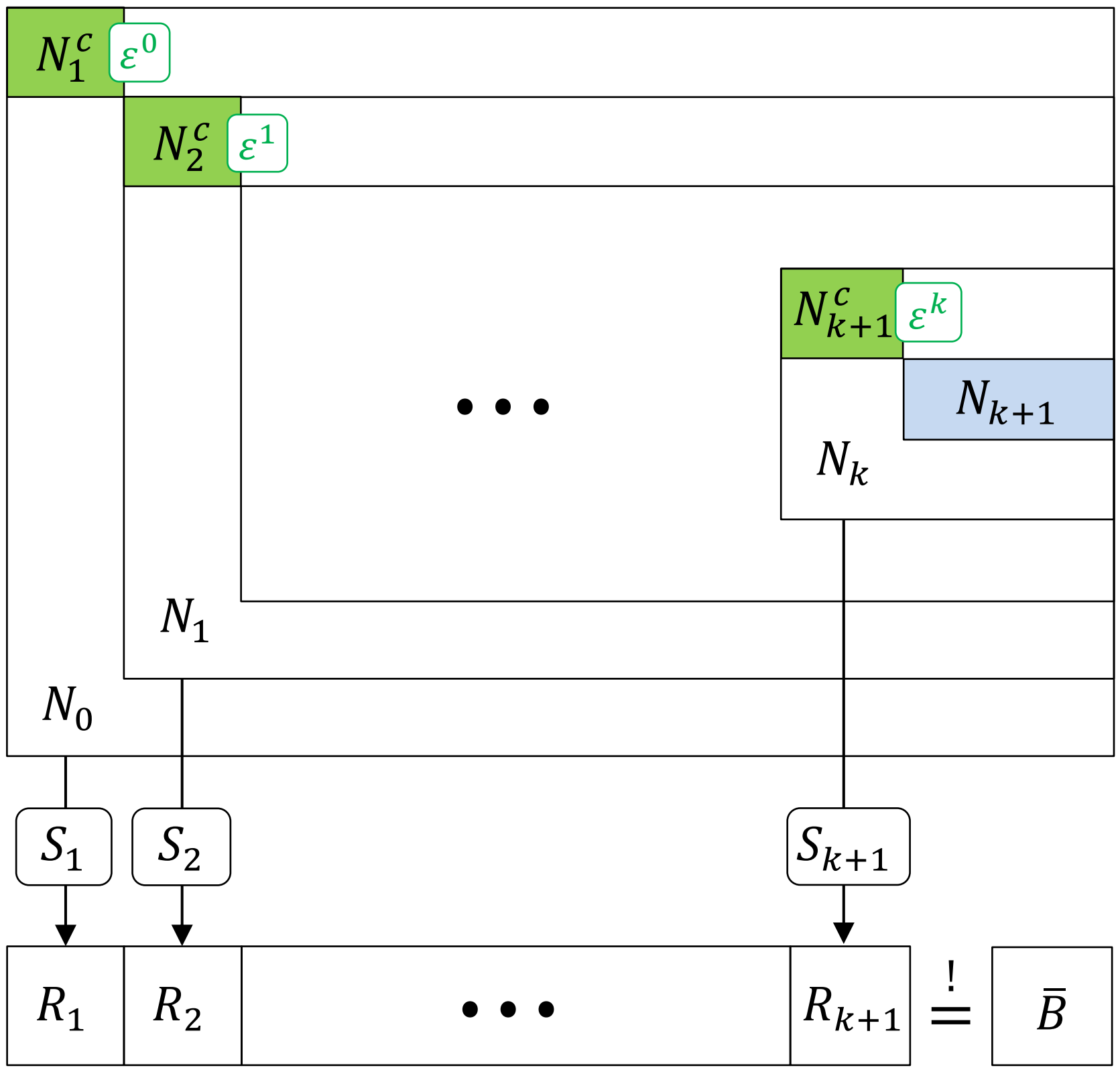

Figure 1 : Resolution of the Banach spaces $B$ and $\bar{B}$.

Thus, within $B = N_0$ a filtration $B = N_0 \supset N_1 \supset \cdots \supset N_{k+1}$ is recursively set up by kernels $N_i$ of appropriate linear operators $S_1, S_2, \dots, S_{k+1}$, where the ranges $R_1, R_2, \dots, R_{k+1}$ of these operators are constructed in a way to ensure, step by step, an increase of the direct sum $R_1 \oplus R_2 \oplus \cdots \oplus R_{k+1}$ within $\bar{B}$. Then, the nondegeneracy condition of the implicit function theorem is satisfied, as soon as a value of $k$ is reached, such that the direct sum agrees with $\bar{B}$ according to

$$R_1 \oplus R_2 \oplus \cdots \oplus R_{k+1} = \bar{B}, \tag{1.3}$$

in this way replacing the *det*-nondegeneracy condition ($N_0$) from above. The implicit function theorem is applied to the variable $n^c$ in (1.2), where $n^c$ is build up by the complements $N_1^c, \dots, N_{k+1}^c$ of the kernels, i.e. $N_0 = N_1^c \oplus N_1, \dots, N_k = N_{k+1}^c \oplus N_{k+1}$.

If the direct sum in (1.3) does not come to an end, then the attempt was unsuccessful. For example, this situation occurs, if we try to approximate a solution curve that is completely contained in the singular locus of $G[z] = 0$.

Now, from a broader perspective, we try to construct a resolution of the above type with respect to an arbitrary curve $z_0(\varepsilon)$ in $B$ with $z_0(0) = 0$, which is not necessarily a high order approximation of $G[z] = 0$. Geometrically, this process is characterized by giving each of the subspaces of figure 1 a precise interpretation, as depicted in figure 2.

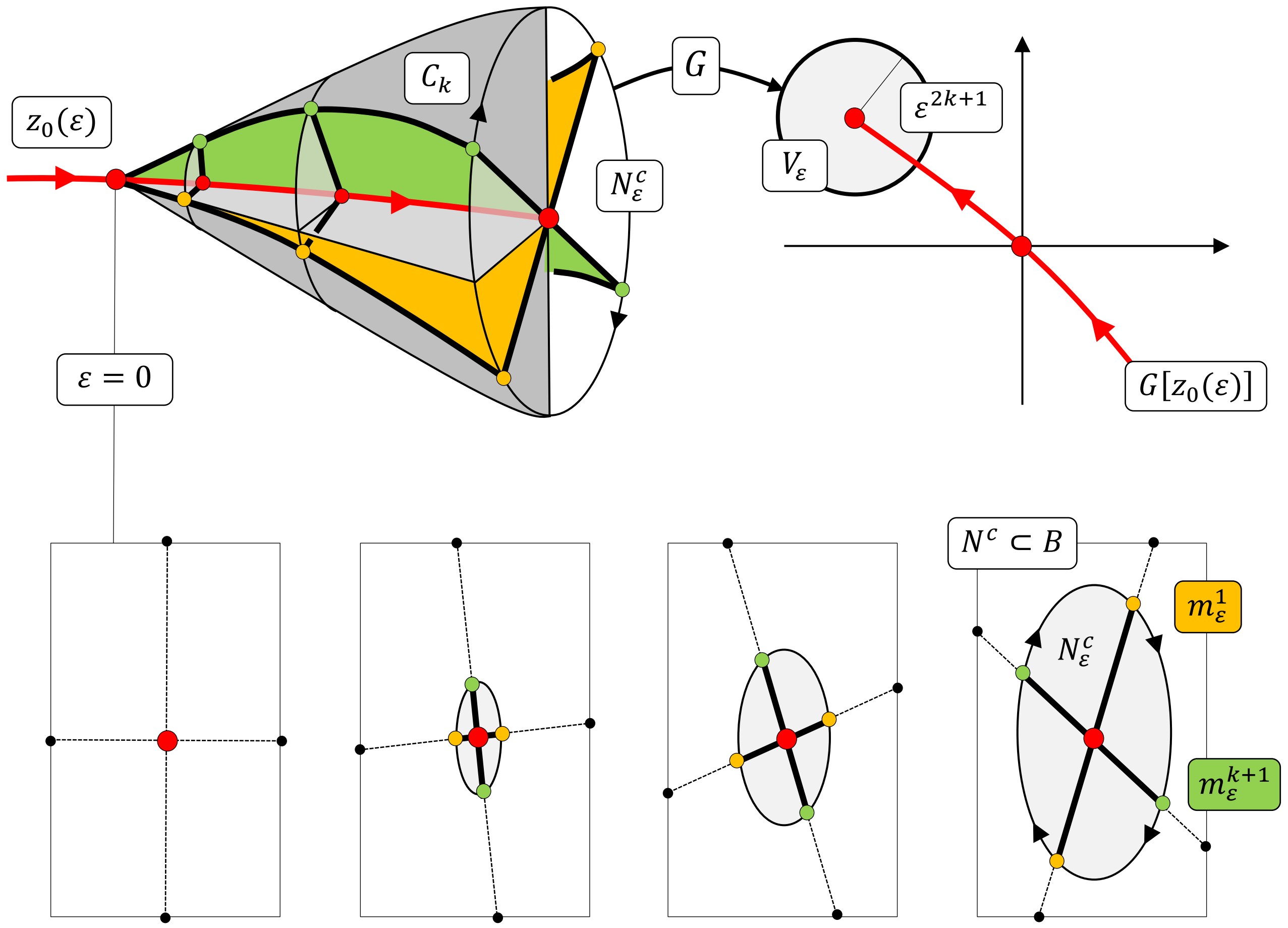

Figure 2 : Geometrical interpretation of a fine resolution of a *k-transversal cone*.

The aim is to construct in $B$ some sort of generalized cone $C_k$ around the given curve, characterized by a simple fibration along the base $z_0(\varepsilon)$ with local fibers $N_\varepsilon^c$, as depicted top-left in the figure. For simplicity, only a half-cone of $C_k$ is drawn.

The construction is complete, as soon as for fixed $\varepsilon \neq 0$, the fibre $N_\varepsilon^c$ is mapped by the leading $\varepsilon$-term of the linearization onto an open neighborhood $V_\varepsilon$ of $G[z_0(\varepsilon)]$ in the image space $\bar{B}$, as indicated top-right. This constellation occurs, as soon as the direct sum condition (1.3) is satisfied. Then we call $C_k$ *a k-transversal or submersive cone.*

For $\varepsilon \neq 0$ fixed, the local fiber $N_\varepsilon^c$ is defined using a subspace $N^c \subset B$ that is restricted to an $\varepsilon$-dependent neighborhood of the origin yielding $N_\varepsilon^c$, as depicted bottom-right. In addition, $N^c$ is further splitted into $k+1$ $\varepsilon$-dependent subspaces $m_\varepsilon^1, \dots, m_\varepsilon^{k+1}$ satisfying $m_\varepsilon^1 \oplus \cdots \oplus m_\varepsilon^{k+1} = N^c$, where the subspaces are derived from the resolution of figure 1.

When moving $\varepsilon$ towards zero, as indicated bottom from right to left, then the subspaces $m_\varepsilon^i$ are slightly rotating around the center line, thereby defining smooth Banach manifolds transversally intersecting each other along the center line $z_0(\varepsilon)$. The restriction of the $m_\varepsilon^i$ to $N_\varepsilon^c$ yields the localized manifolds as depicted top-left by orange and green surfaces. In this sense, the generalized cone $C_k$ is spanned by $k+1$ ruled surfaces.

The importance of the smooth manifolds defined by $m_\varepsilon^i$ is given by the fact that they provide us with detailed information about the speed of variation of $G$ when moving away from $z_0(\varepsilon)$ with-in one of the manifolds. More precisely, the norms of partial derivatives in the directions of the smooth manifolds are exactly given by

$$\| G_{m_\varepsilon^1}[\, z_0(\varepsilon)\,] \| = O(\, |\varepsilon|^0) \;, \quad \cdots \quad , \| G_{m_\varepsilon^{k+1}}[\, z_0(\varepsilon)\,] \| = O\big(\, |\varepsilon|^k\big) \,, \tag{1.4}$$

i.e. when moving through the curve $z_0(\varepsilon)$, thereby monitoring successively the partial derivatives in the directions of $m_\varepsilon^1, \dots, m_\varepsilon^{k+1}$, then we obtain highest values in the direction of $m_\varepsilon^1$ by order of exactly $|\varepsilon|^0$ and lowest values in the direction of $m_\varepsilon^{k+1}$ by order of exactly $|\varepsilon|^k$. In this sense a stepwise fine resolution concerning fast and slow variation of $G$ is constructed in the cone. The partial derivatives are explicitly given using the linear operators $S_1, S_2, \dots, S_{k+1}$ from figure 1.

Further, along the manifolds $m_\varepsilon^1, \dots, m_\varepsilon^{k+1}$, the cone is shrinking to zero with different orders of magnitude with respect to $\varepsilon$, implying an ellipsoid $N_\varepsilon^c$ to occur for small values of $\varepsilon$. More precisely, along $m_\varepsilon^1, \dots, m_\varepsilon^{k+1}$, the cone is shrinking with radius of order $|\varepsilon|^{2k+1}, \dots, |\varepsilon|^{k+1}$, implying the product of partial derivatives from (1.4) times radius of cone to be identical of order $O\big(|\varepsilon|^{2k+1}\big)$ in every direction of the cone, i.e. the larger the partial derivative in a certain direction, the smaller the radius of the cone in this direction. This balancing implies that the image of a fiber under $G$ results in the open ball $V_\varepsilon$ comprising $G[z_0(\varepsilon)]$ in $\bar{B}$ with homogenous radius of order $|\varepsilon|^{2k+1}$.

As already mentioned, a submersive cone, satisfying the rate conditions (1.4), can typically be derived with respect to an arbitrary curve $z_0(\varepsilon)$ in $B$ with $z_0(0) = 0$. However, in the next step, let us add some approximation properties. First, if the distance of $G[z_0(\varepsilon)]$ to the origin is greater than the radius of the ball $V_\varepsilon$, as shown in figure 2 top-right, then obviously no element of $N_\varepsilon^c$ is mapped by $G$ to the orgin of $\bar{B}$ and no solution of $G[z] = 0$ exists within the fiber $N_\varepsilon^c$. On the other hand, if the distance of $G[z_0(\varepsilon)]$ to the origin is smaller or equal to the radius of the ball $V_\varepsilon$ according to $\|G[z_0(\varepsilon)]\| = O\big(|\varepsilon|^{2k+1}\big)$, then at least one point within $N_\varepsilon^c$ is mapped to the origin in $\bar{B}$, i.e. for every $\varepsilon \neq 0$ at least one solution of $G[z] = 0$ occurs in $N_\varepsilon^c$.

Now, if $C_k$ is a $k$-transversal cone, then by the implicit function theorem, we will see that all solutions in the cone are given by a smooth manifold $z(\varepsilon, p),\ \varepsilon \neq 0$ with parameter $p$ taken from a certain subspace $P$ of $B$ with dimension of $P$ between zero and infinity. In case of $\dim(P) = 0$, exactly one solution curve exists in the cone.

On the other hand, if $\|G[z_0(\varepsilon)]\| \neq O\big(|\varepsilon|^{2k+1}\big)$, we can be sure that the cone represents for $\varepsilon \neq 0$ a solution free domain in $B$.

The construction of the cone $C_k$ depends exactly on first $k$ coefficients of the Taylor expansion $z_0(\varepsilon) = \varepsilon \cdot \bar{z}_1 + \cdots + \frac{1}{k!}\varepsilon^k \cdot \bar{z}_k + O\big(|\varepsilon|^{k+1}\big)$ and on $k+1$ first derivatives $G_0^1, \dots, G_0^{k+1}$ of $G$ at $z = 0$. In particular, a $k$-transversal cone cannot be destroyed by perturbations of $G[z]$ of order $O\big(\|z\|^{k+2}\big)$. In this sense, the nondegeneracy condition remains valid under perturbations of order $O\big(\|z\|^{k+2}\big)$. Further, if solutions $z(\varepsilon, p)$ exist within the cone, then the solutions cannot be destroyed by perturbations of order $O\big(\|z\|^{2k+1}\big)$. This corresponds to the fact that the approximation condition is of order $2k$. Hence, perturbations of $G[z]$ between $O\big(\|z\|^{k+2}\big)$ and $O\big(\|z\|^{2k}\big)$ may destroy solutions, but not the transversality of the cone.

We also note, that by a linear diffeomorphic transformation of $N_\varepsilon^c$ and a near identity transformation of $V_\varepsilon$, the ruled surfaces from figure 2 defined by $m_\varepsilon^1, \dots, m_\varepsilon^{k+1}$ may further be straightened out and $V_\varepsilon$ may be transformed to an exact ball, i.e. by setting $G_\varepsilon[z] := G[z_0(\varepsilon) + z]$, an $RL$-transformation yields the linearization

$$\mathcal{J}_\varepsilon \circ G_\varepsilon \circ A_\varepsilon = \varepsilon^{2k+1} \cdot L \tag{1.5}$$

with an $\varepsilon$-independent linear mapping $L \in GL[N^c, \bar{B}]$, the linear bijection $A_\varepsilon \in GL[N^c, N^c]$, $\varepsilon \neq 0$ from(1.2), as well as a near identity transformation $\mathcal{J}_\varepsilon = I_{\bar{B}} + O(|\varepsilon|)$ of the image space $\bar{B}$. Formula (1.5) represents some sort of normal form of a $k$-transversal cone. Possibly, (1.5) may be compared to general linearization techniques in Banach spaces derived in [HM], [BH].

Up to this point, all spaces are allowed to be Banach spaces of infinite dimensions. Now, when restricting to finite dimensions, then the information (1.4) about the behaviour of the linearization along $z_0(\varepsilon)$ can directly be used to perform an $\varepsilon$-dependent decoupling of the system on the linear level and to calculate the overall behaviour of the determinant of the linearization with respect to the complement $N^c$ by

$$det\{\, G_{N^c}[\, z_0(\varepsilon)\,]\,\} = \varepsilon^{\chi} \cdot \underbrace{r(\varepsilon)}_{\neq 0} \qquad with \quad \chi := 1 \cdot \underbrace{dim(m_\varepsilon^2)}_{\geq 0} + \cdots + k \cdot \underbrace{dim\big(m_\varepsilon^{k+1}\big)}_{\geq 1}\,. \tag{1.6}$$

Obviously, the characteristic number $\chi$ satisfies $\chi \geq k$, thus offering the possibility to ascertain solution curves in case of $\|G[z_0(\varepsilon)]\| = O\big(|\varepsilon|^{2k+1}\big)$ and $det\{G_{N^c}[z_0(\varepsilon)]\} \neq O(|\varepsilon|^{\chi+1})$, which means that the determinant is allowed to vary much slower than required by standard $(N_0)$ according to $det\{G_{N^c}[z_0(\varepsilon)]\} \neq O\big(|\varepsilon|^{k+1}\big)$. Note that the only constellation with $\chi = k$ occurs in case of $m_\varepsilon^2 = \cdots = m_\varepsilon^k = \{0\}$ and $dim\big(m_\varepsilon^{k+1}\big) = 1$.

This is well known, that measuring the partial derivatives in complementary directions to the approximate solution $z_0(\varepsilon)$ in terms of the determinant may not be optimal [F], [BK], [S4]. In particular, in [F] a generalization of Hensel's Lemma is obtained by resolving space appropriately under consideration of different rates of variation. Fisher also remarks that a more geometric interpretation and application of his concept would be desirable. We hope to perform a small step in this direction.

The construction of an appropriate filtration is also performed in [BK], considering equations between modules or even free abelian groups. We are quite confident that our approach may also work in a more general setting. The techniques used, completely rely on rather basic constructions centered around chain rule and/or substitution of power series into power series with corresponding system of undetermined coefficients. In particular, the techniques apply in infinite dimensions, such as Banach spaces or possibly in cord spaces [BH] as well.

From a more functional analytical point of view, the direct sum condition (1.3) may represent some sort of generalization of the $k$-transversality concept concerning families of Fredholm operators $L(\lambda)$, as introduced in [E], [EL], [L1], [LM] with respect to a known trivial solution curve $y = 0$ of $G[\lambda, y] = L(\lambda)y + O(\|y\|^2)$. In that context the leading term of the determinant of the Jacobian with respect to $y$ is calculated by $det\big\{G_y[\,\lambda, 0]\big\} = \lambda^{\chi} \cdot r(\lambda)$, where the exponent $\chi$ is derived by comparable methods as in the paper at hand. In [EL], the exponent $\chi$ is introduced as *generalized algebraic multiplicity of the operator family* $L(\lambda)$. Further, the leading term of the determinant is used to calculate the behaviour of the topological degree that ascertains secondary bifurcation, as soon as sign change of the degree occurs along the trivial solution $y = 0$.

Obviously by (1.6), these lines of reasoning can be transferred to a $k$-transversal cone, supposed it contains a unique solution curve $z(\varepsilon)$. The corresponding constellation is indicated in figure 3.

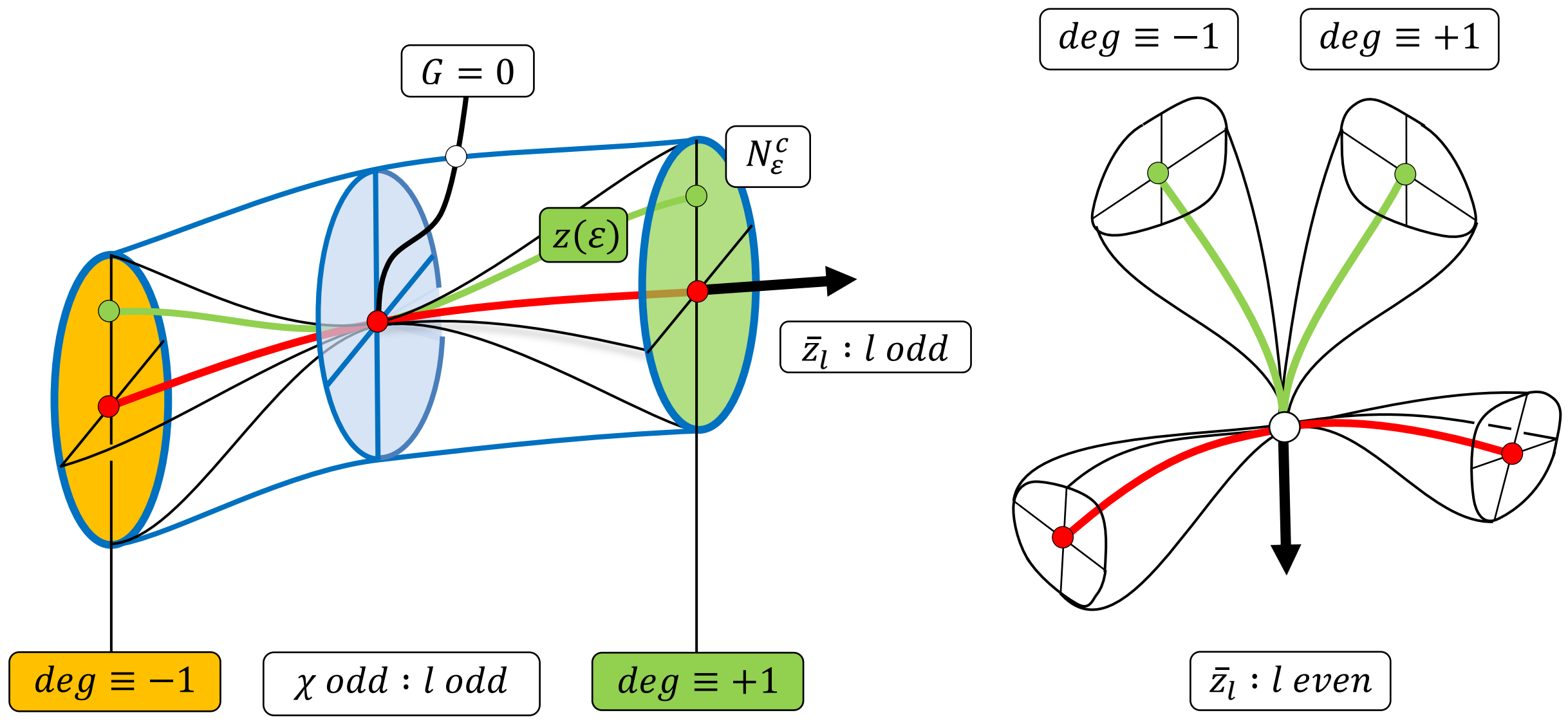

Figure 3 : Topological degree and secondary bifurcation.

First, the topological degree of $G$ with respect to the section $N_\varepsilon^c$ is given by the sign of the determinant at the solution $z(\varepsilon)$. Now, it is interesting to note that the characteristic number $\chi$ from (1.6) is a property of the complete cone, implying $det\{G_{N^c}[z(\varepsilon)]\} = \varepsilon^\chi \cdot r(\varepsilon),\ r(\varepsilon) \neq 0$ along the solution curve $z(\varepsilon)$ as well. Hence, supposed $\mathbb{K} = \mathbb{R}$, we obtain constant degree $deg \equiv -1$ or $deg \equiv +1$ in each half-cone.

Now, if $\chi$ is odd, then different signs occur in the two half-cones and, if additionally the leading coefficient $\bar{z}_l \neq 0$ of $z_0(\varepsilon) = \frac{1}{l!}\varepsilon^l \cdot \bar{z}_l + \cdots + \frac{1}{k!}\varepsilon^k \cdot \bar{z}_k + \cdots$ satisfies $l\ odd$, then a cylindrical construction can be set up, as indicated in the left diagram of figure 3 by blue color. The cylindrical construction is designed to connect the two half-cones (with different degrees), thus implying by homotopy invariance of the degree at least one point on the cylinder to exist with $G = 0$. Finally, by shrinking of the cylinder, a continuum of points with $G = 0$ is shown to emanate from $z = 0$. In this sense a $k$-transversal cone is perfectly adapted to deal locally with secondary bifurcation based on topological arguments. Far reaching local and global bifurcation results can be found in [LM] and [L2] in the context of $G[\lambda, y] = L(\lambda)y + O(\|y\|^2)$.

If $l \geq 2$ is *even*, as shown by the cusp curve on the right side of figure 3 (green curve), then the half cones may have different topological degrees, but the cylindrical construction fails and no results concerning secondary bifurcation are possible, at least based on standard degree theory.

Summarizing, within the Banach space $B$, we aim to construct different cones out of the singularity at $z = 0$ that are characterized by some kind of homogeneity. Above all, homogeneity with respect to $G$ expansion rates as in (1.4), yielding the possibility to establish complete linearization of $G$ within the cone as in (1.5), as well as offering the possibility to use topological techniques concerning secondary bifurcation.

From a mathematical point of view, it is quite clear that each cone is simply blown-up to a cylinder, when approaching the singularity from the direction of the cone as $\varepsilon \to 0$, i.e. a certain desingularization of the singularity with respect to the cone is performed. During this process, the operator norm of the inverse has to explode. However, this explosion is uniformly bounded within the complete cone by $O(|\varepsilon|^{-k})$, i.e. in the notion of [L1] a necessary condition for a $k$-*transversal cone* to occur is given by a *generalized algebraic eigenvalue of order* $k$.

Reversely, we will see that under appropriate assumptions, explosion of the operator norm of the inverse by $O(|\varepsilon|^{-k})$ is sufficient too, concerning the existence of a *$k$-transversal cone*. Hence, building up a *$k$-transversal cone $C_k$* may also be interpreted as a constructive method for calculation of the algebraic order $k$ belonging to the operator norm of the inverse of an isolated singularity.

All this shows that our techniques are restricted to curves $z_0(\varepsilon)$, possibly touching by high order the singular locus of $G$, but ultimately traversing it, in this way defining an isolated singularity of an operator family. Hence, the constructed cone merely represents a detailed, quantitative description of a domain, comprising the curve $z_0(\varepsilon)$, where submersion holds.

Except the considerations concerning the topological degree, we do not suppose finite dimensions or other Fredholm properties, but usually we assume all subspaces of figure 1 to be closed with continuous projections.

One may hope to understand the singularity in a better way, if the conditions in some or all of the cones comprising solution curves are known. A result of this kind appears in [DR], where it is shown that all solutions of $G[z] = 0$ are found, supposed that only a finite number of solution curves exist with each curve satisfying the nondegeneracy condition (1.3) with $k = 1$, i.e. every cone is spanned by exactly $k + 1 = 2$ ruled surfaces. However, this result only applies in case of $G'[0]$ to be a Fredholm operator. In [ J ], the result is partly extended to $k = 2$.

From a technical point of view, the paper intends to combine the well-known concept from algebra of *$2k$-approximation/$k$-nondegeneracy* with the basic concept of *$k$-transversality* from functional analysis introduced in [EL].

In section 2, some basic patterns, valid within the system of undetermined coefficients are introduced and a first existence result concerning solutions of $G[z] = 0$ is proved and summarized in Theorem 1.

Section 3 contains the main results of the paper concerning uniqueness, perturbation, expansion rates and norm of the inverse, comprised in Theorem 2. In Corollary 1, we give conditions for a *generalized algebraic eigenvalue of order $k$* [L1] to deliver a *$k$-transversal cone*.

In section 4, we restrict to finite dimensions, possibly after performing a Lyapunov-Schmidt reduction, for dealing with secondary bifurcation from a given solution curve in general position, based on Brouwer's degree and summarized within Corollary 2. Section 4 also contains an example, which aims to convince the reader that formulas concerning general position might be useful. Other examples can be found in [S3].

Sections 5, 6 and 7 contain complete proofs.

The main results and proofs are not new, see [S1]. The interpretation of some of the results may be new. In addition, the patterns concerning the system of undetermined coefficients might be of some interest.

## *2. Patterns within the System of undetermined Coefficients*

In this section, we basically focus on the description of some characteristic patterns within the system of undetermined coefficients. As an application, we present a first theorem concerning the existence of local solution curves through $z = 0$, parametrized by the external parameter $\varepsilon \in \mathbb{K}$ . Depending on the dimension of the smallest subspace $N_{k+1}$ within the filtration of $B$ in figure 1 (blue box), the solution curves will also depend on further parameters implying the existence of a smooth solution manifold of higher (or infinite) dimension.

In some more detail, an equation $G[z] = 0, G \in C^q(B,\bar{B})$ with $G[0] = 0$ and $B, \bar{B}$ real or complex Banach spaces is considered, with the aim of finding solution curves $G[z(\varepsilon)] = 0$ through the origin with $\varepsilon \in \mathbb{K} = \mathbb{R},\ \mathbb{C},\ |\varepsilon| \ll 1$.

The basic technique to derive $z(\varepsilon)$ is very simple. For $k \geq 1$, insert the ansatz $z = \varepsilon \cdot z_1 + \cdots + \frac{1}{(2k+1)!}\,\varepsilon^{2k+1} \cdot z_{2k+1}$ into $G$ and perform a Taylor expansion at $\varepsilon = 0$ up to order $2k+1$ by

$$G[\,\varepsilon \cdot z_1 + \cdots + \frac{\varepsilon^{2k+1}}{(2k+1)!} \cdot z_{2k+1}\,] = \sum_{i=1}^{2k+1} \frac{1}{i!}\,\varepsilon^i \cdot T^i[\,z_i, \ldots, z_1] \;+\; \varepsilon^{2k+2} \cdot r[\,\varepsilon, z_{2k+1}, \ldots, z_1\,] \tag{2.1}$$

with smooth remainder map $r[\cdot]$ , $q \geq 2k+2$ and coefficients $T^i[z_i, \ldots, z_1],\ i \geq 1$, building up the so called system of undetermined coefficients according to

$$\begin{aligned} T^1[\,z_1] \quad &= \; G_0^1 \cdot z_1 = 0 \\ T^2[\,z_2, z_1] \quad &= \; G_0^1 \cdot z_2 \;+\; G_0^2 \cdot z_1^2 = 0 \\ T^3[\,z_3, z_2, z_1] \; &= \; G_0^1 \cdot z_3 \;+\; 3G_0^2 \cdot z_1 z_2 + \; G_0^3 \cdot z_1^3 = 0 \\ \vdots \qquad\quad & \qquad\qquad \vdots \end{aligned} \tag{2.2}$$

Here $G_0^\beta$ denotes the $\beta$-th derivative of $G$ at $z = 0$. In general, by higher order chain rule [AMR], we obtain for $i \geq 1$ the expression

$$T^i[\,z_i, \ldots, z_1] = \sum_{\beta=1}^{i} G_0^\beta \sum_{\substack{n_1+\cdots+n_i=\beta \\ 1\cdot n_1+\cdots+i\cdot n_i=i}} \frac{i!}{n_1! \cdots n_i!} \prod_{\tau=1}^{i} \left(\frac{1}{\tau!}\,z_\tau\right)^{n_\tau} \in \bar{B} \tag{2.3}$$

with $T^i$ depending explicitly from $z_1, \ldots, z_i$ and $G_0^1, \ldots, G_0^i$. In this sense, the $i$-th coefficient $T^i$ within the $\varepsilon$-expansion (2.1) is given by a sum composed of first $i$ coefficients of the ansatz $z = \varepsilon \cdot z_1 + \cdots$ and first $i$ coefficients of the Taylor expansion of $G$ at $z = 0$. Moving $k$ towards infinity, this process may also be interpreted as plugging the $\varepsilon$-power series of the ansatz into the $z$-power series of $G$. This viewpoint is taken in [S4], where power series solutions of differential algebraic equations are investigated.

In the next step, we summarize some results concerning the system of undetermined coefficients. First, by direct inspection of (2.3), the following linearity structure with respect to $[T^{2k}, \ldots, T^{k+1}]$ and higher order coefficients $[z_{2k}, \ldots, z_{k+1}]$ is valid for $k \geq 1$

$$\begin{pmatrix} T^{2k}[z_{2k}, \ldots, z_1] \\ \vdots \\ T^{k+1}[z_{k+1}, \ldots, z_1] \end{pmatrix} = \Delta^k(z_{k-1}, \ldots, z_1) \cdot \begin{pmatrix} z_{2k} \\ \vdots \\ z_{k+1} \end{pmatrix} + \mathrm{I}^k(z_k, \ldots, z_1) \tag{2.4}$$

with $\Delta^1(z_1) = G_0^1$ in case of $k = 1$. Here $\Delta^k(z_{k-1}, \ldots, z_1) \in L[B^k, \bar{B}^k]$ denotes an upper triangular matrix operator with elements $\Delta_{i,j}^k(z_{k-1}, \ldots, z_1) \in L[B, \bar{B}]$ and $\mathrm{I}^k(z_k, \ldots, z_1)$ an element in $\bar{B}^k$, both explicitly defined in section 5.

In addition and for later use, we state the following relation for $k \geq 1$

$$\begin{pmatrix} T^k[z_k,\dots,z_1] \\ \vdots \\ T^1[z_1] \end{pmatrix} = \underbrace{\left(\Gamma^k\right)^{-1}}_{=:\,\Gamma^{-k}} \cdot \Delta^k(z_{k-1},\dots,z_1)\cdot \Gamma^k \cdot \begin{pmatrix} z_k \\ \vdots \\ z_1 \end{pmatrix} \tag{2.5}$$

$$with \quad \Gamma^k = diag\left[\,\Gamma_k^k,\dots,\Gamma_1^k\,\right] \quad and \quad \Gamma_i^k = \binom{k+i}{i-1},\ \ i=1,\dots,k\,,$$

i.e. the first $k$ coefficients $[T^1,\dots,T^k]$ of the power series (2.1) can completely be expressed by the linear part $\Delta^k$ of the next $k$ coefficients $[T^{k+1},\dots,T^{2k}]$ from (2.4). In this sense, $\Delta^k$ delivers some sort of recurrency structure that is playing a key role concerning the investigation of the system of undetermined coefficients.

Finally, coefficients above $T^{2k}$ can be formulated for $k\ge 0, l\ge 0$ by lower ones according to

$$T^{2k+1+l}(\,z_{2k+1+l},\dots,z_{k+1+l},\ \ z_{k+l},\dots,z_{k+1},\ \ z_k,\dots,z_1\,)$$

$$= \left[\,T_{z_0}^0(z_0)\ \ T_{z_1}^2(z_1)\ \cdots\ T_{z_k}^{2k}(z_k,\dots,z_1)\,\right]\cdot \boxed{C^{2k+l}}\cdot \begin{pmatrix} z_{2k+1+l} \\ \vdots \\ z_{k+1+l} \end{pmatrix} + R_{2k+1+l}(z_{k+l},\dots,z_1) \tag{2.6}$$

with

$$\boxed{C^{2k+l}} = Diag\left[\,\gamma_0^{2k+l}\ \cdots\ \gamma_k^{2k+l}\,\right]\ \ and\ \ \gamma_t^{2k+l} := \binom{2k+1+l}{t}\cdot\binom{2t}{t}^{-1},\ \ t=0,\dots,k \tag{2.7}$$

Here, $T_{z_0}^0(z_0) := G_0^1$ and $T_{z_i}^{2i}(z_i,\dots,z_1)$, $i=1,\dots,k$ denote partial derivatives of $k+1$ low order coefficients $T^2$, $T^4,\dots,T^{2k}$ with derivatives taken with respect to low $z$-coefficients $z_1,\dots,z_k$ respectively, implying a second recursion principle. Formula (2.6) may be interpreted as an extended version of a formula by Hurwitz [Hu] expressing high order derivatives of an expansion by lower ones. For comparable results see [DL], Lemma 2.2 and [VFZ], Theorems 3.5, 3.6.

The proofs of (2.4), (2.5) follow by simple and elementary calculation from (2.1), omitted in the paper at hand. The proof of (2.6) can be found in section 7 formulas (7.1)-(7.6).

Next, assume the existence of $[\bar z_{2k},\dots,\bar z_1]$ with $T^{2k}[\bar z_{2k},\dots,\bar z_1] = \cdots = T^1[\bar z_1] = 0$, implying

$$\|\,G[\,\varepsilon\cdot\bar z_1+\cdots+\tfrac{1}{(2k)!}\,\varepsilon^{2k}\cdot\bar z_{2k}\,]\,\| = O(\,|\varepsilon|^{2k+1}\,) \tag{2.8}$$

by (2.1) with $z_{2k+1}=0$. Thus, the curve $z_0(\varepsilon)=\varepsilon\cdot\bar z_1+\cdots+\frac{1}{(2k)!}\,\varepsilon^{2k}\cdot\bar z_{2k}$ defines an approximate solution curve of $G[z]=0$ of order $2k$, which obviously corresponds to the approximation condition $(A_0)$ of the introduction.

Now, for given $[\bar z_k,\dots,\bar z_1]\in B^k$, there exists a solution $[\bar z_{2k},\dots,\bar z_{k+1}]\in B^k$ of (2.4) if and only if the range $R[\Delta^k(\bar z_{k-1},\dots,\bar z_1)]$ contains the inhomogeneity $\mathrm{I}^k(\bar z_k,\dots,\bar z_1)$. But then, the solution $[\bar z_{2k},\dots,\bar z_1]$ of $T^{2k}=\cdots=T^1=0$ is embedded within an affine subspace of solutions given by $[z_k,\dots,z_1]=[\bar z_k,\dots,\bar z_1]$ and

$$\begin{pmatrix} z_{2k} \\ \vdots \\ z_{k+1} \end{pmatrix} = \begin{pmatrix} \bar z_{2k} \\ \vdots \\ \bar z_{k+1} \end{pmatrix} + \begin{pmatrix} b_{2k} \\ \vdots \\ b_{k+1} \end{pmatrix} \quad with \quad \begin{pmatrix} b_{2k} \\ \vdots \\ b_{k+1} \end{pmatrix} \in N[\,\Delta^k(\bar z_{k-1},\dots,\bar z_1)\,] =: N[\Delta^k]\,. \tag{2.9}$$

Here $N[\Delta^k]$ denotes the nullspace of the operator in brackets. Note also, if the leading coefficients $[\bar z_k,\dots,\bar z_1]$ are fixed, then all solutions of $T^{2k}=\cdots=T^1=0$ are given by (2.9).

Further, every element $[b_{2k},\dots,b_{k+1}]$ from $N[\Delta^k]$ implies an approximate solution curve of order $2k$ within $B$ according to

$$z=\varepsilon\cdot\bar{z}_1+\dots+\frac{1}{k!}\varepsilon^k\cdot\bar{z}_k+\frac{1}{(k+1)!}\varepsilon^{k+1}\cdot(\bar{z}_{k+1}+b_{k+1})+\dots+\frac{1}{(2k)!}\varepsilon^{2k}\cdot(\bar{z}_{2k}+b_{2k})$$

$$+\frac{1}{(2k+1)!}\varepsilon^{2k+1}\cdot b_{2k+1}\,, \tag{2.10}$$

where the last summand of order $2k+1$ and $b_{2k+1}\in B$ does not disturb this property and is merely added for further reasoning. Note that $b_{2k+1}$ occurs the first time within $T^{2k+1}[\cdot]$ .

All the approximate solution curves in (2.10), obtained by variation of $[b_{2k+1},\dots,b_{k+1}]\in B\times N[\Delta^k]$, share the same leading coefficients up to order $k$, defining within the Banach space $B$ some sort of cone of approximations of order $2k$, as qualitatively indicated in figure 4.

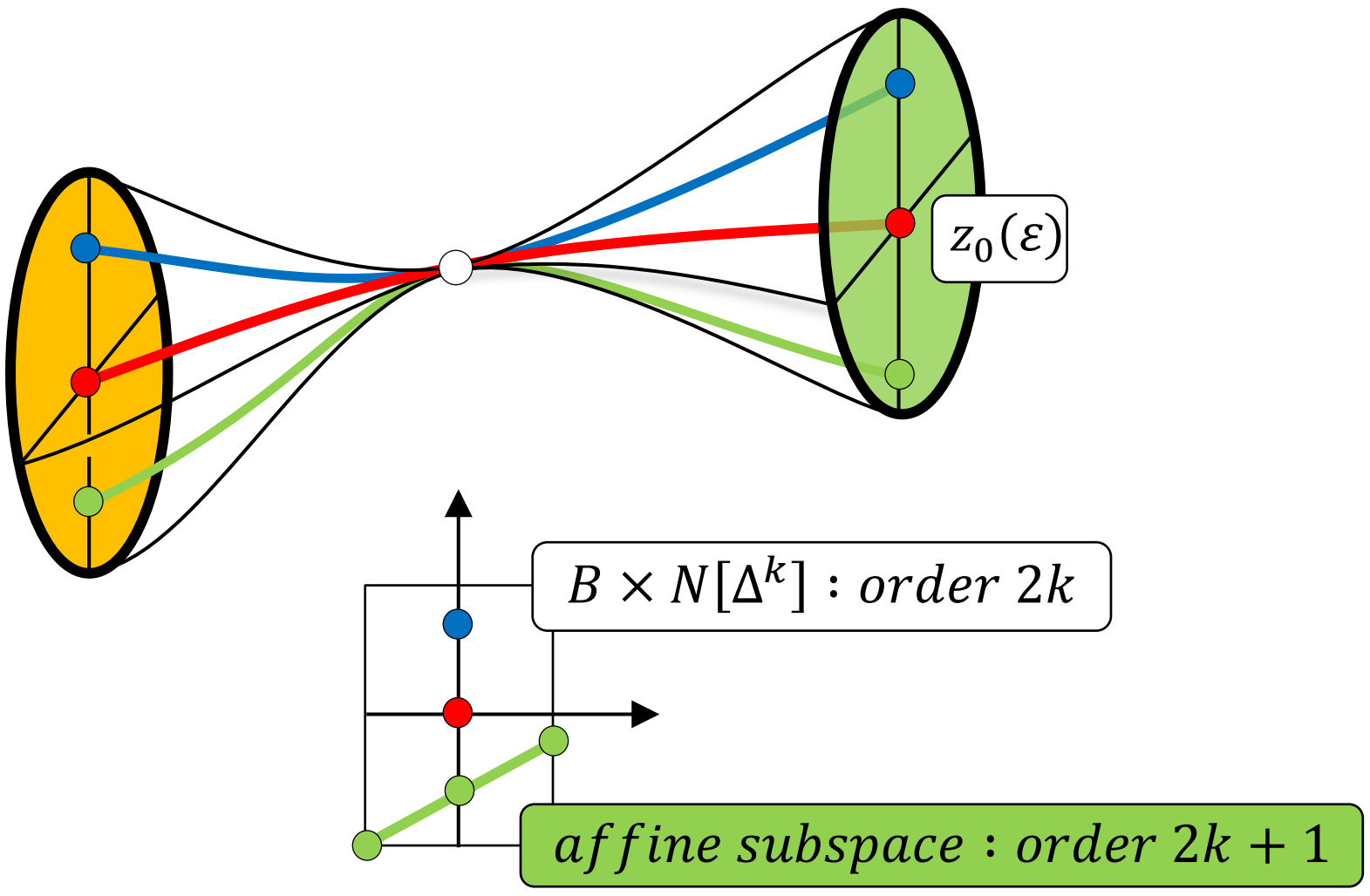

Figure 4 : Cone of approximations of order $2k$ (red, blue) and $2k+1$ (green).

The red center line $z_0(\varepsilon)$ is obtained with $[b_{2k+1},\dots,b_{k+1}]=0$, whereas the blue line represents an approximation with $[b_{2k+1},\dots,b_{k+1}]\neq 0$.

It remains to show that, within this cone, even true solutions exist, if some further structure is given to the cone that finally allows for application of the implicit function theorem with respect to $[b_{2k+1},\dots,b_{k+1}]$. In this sense, the initial approximation $z_0(\varepsilon)$ of order $2k$ is only needed to construct the affine subspace (2.9), in this way offering variability to the coefficients $z_{k+1}=\bar{z}_{k+1}+b_{k+1},\dots,z_{2k+1}=b_{2k+1}$ within (2.10). Note also, that this construction precisely explains the drop from initial approximation of order $2k$ down to true solutions agreeing merely up to order of $k$ with the initial approximation $z_0(\varepsilon)=\varepsilon\cdot\bar{z}_1+\dots+\frac{1}{(2k)!}\varepsilon^{2k}\cdot\bar{z}_{2k}$.

Now, plugging the affine subspace (2.9), (2.10) into (2.1) and cancelling out $\varepsilon^{2k+1}$, we end up with the smooth remainder equation

$$T^{2k+1}[\underbrace{z_{2k+1}}_{=\,b_{2k+1}},\begin{pmatrix}\bar{z}_{2k}\\ \vdots\\ \bar{z}_{k+1}\end{pmatrix}+\begin{pmatrix}b_{2k}\\ \vdots\\ b_{k+1}\end{pmatrix},\bar{z}_k,\dots,\bar{z}_1]+\varepsilon\cdot r[\varepsilon,\,b_{2k+1},\,b_{2k},\dots,b_{k+1}] \tag{2.11}$$

$$=:H(\varepsilon,\,b_{2k+1},\,b_{2k},\dots,b_{k+1})=0$$

and $H \in C^{q-2k-1}\big(U \times B \times N[\Delta^k],\ \bar{B}\big),\ U = \{\varepsilon \in \mathbb{K} : |\varepsilon| \ll 1\}$. Further, direct inspection of (2.3) shows the linearity property

$$T^{2k+1}[\, z_{2k+1}, \dots, z_1 \,] = W^{2k+1}(z_k, \dots, z_1) \cdot \begin{pmatrix} z_{2k+1} \\ \vdots \\ z_{k+1} \end{pmatrix} + R^{2k+1}(z_k, \dots, z_1) \in \bar{B} \qquad (2.12)$$

with $W^{2k+1}(z_k, \dots, z_1) \in L[B^{k+1}, \bar{B}]$ and $R^{2k+1}(z_k, \dots, z_1) \in \bar{B}$, which are again defined explicitly in section 5. Therefore, the remainder equation (2.11) reads for $\varepsilon = 0$

$$H(0,\, b_{2k+1},\, b_{2k}, \dots, b_{k+1}) = T^{2k+1}[\, b_{2k+1},\ \begin{pmatrix} \bar{z}_{2k} \\ \vdots \\ \bar{z}_{k+1} \end{pmatrix} + \begin{pmatrix} b_{2k} \\ \vdots \\ b_{k+1} \end{pmatrix},\ \bar{z}_k, \dots, \bar{z}_1 \,]$$

$$= W^{2k+1}(\bar{z}_k, \dots, \bar{z}_1) \cdot \left\{ \begin{pmatrix} 0 \\ \bar{z}_{2k} \\ \vdots \\ \bar{z}_{k+1} \end{pmatrix} + \begin{pmatrix} b_{2k+1} \\ b_{2k} \\ \vdots \\ b_{k+1} \end{pmatrix} \right\} + R^{2k+1}(\bar{z}_k, \dots, \bar{z}_1)$$

$$= W^{2k+1}(\bar{z}_k, \dots, \bar{z}_1) \cdot \begin{pmatrix} b_{2k+1} \\ \vdots \\ b_{k+1} \end{pmatrix} + \bar{R}^{2k+1}(\bar{z}_{2k}, \dots, \bar{z}_{k+1},\, \bar{z}_k, \dots \bar{z}_1) = 0\,. \qquad (2.13)$$

Now, if the linear mapping

$$W^{2k+1}(\bar{z}_k, \dots, \bar{z}_1) \in L[\, B \times N[\Delta^k],\ \bar{B}\,]$$

is surjective, then the mapping can also be restricted to a direct complement $N_c$ of its closed kernel to obtain a bijection $W^{2k+1}(\bar{z}_k, \dots, \bar{z}_1) \in L[N_c, \bar{B}]$ by decomposition of the Banach space $B \times N[\Delta^k]$ according to

$$B \times N[\Delta^k] = {}^{B \times N[\Delta^k]}\!\Big/_{N[\, W^{2k+1}(\bar{z}_k, \dots, \bar{z}_1)\,]} \oplus N[\, W^{2k+1}(\bar{z}_k, \dots, \bar{z}_1)\,]$$

$$=: N_c \oplus N[\, W^{2k+1}(\bar{z}_k, \dots, \bar{z}_1)\,]\,. \qquad (2.14)$$

Hence, for $\varepsilon = 0$, the remainder equation $H(0, b_{2k+1}, b_{2k}, \dots, b_{k+1}) = T^{2k+1}[\cdot] = 0$ is given by the linear equation (2.13) with corresponding affine solution subspace

$$\begin{pmatrix} b_{2k+1} \\ \vdots \\ b_{k+1} \end{pmatrix} = \begin{pmatrix} \bar{b}_{2k+1} \\ \vdots \\ \bar{b}_{k+1} \end{pmatrix} + \begin{pmatrix} c_{2k+1} \\ \vdots \\ c_{k+1} \end{pmatrix} \quad with \quad \begin{pmatrix} c_{2k+1} \\ \vdots \\ c_{k+1} \end{pmatrix} \in N[\, W^{2k+1}(\bar{z}_k, \dots, \bar{z}_1)\,] \qquad (2.15)$$

and $\bar{b} := \big(\bar{b}_{2k+1}, \dots, \bar{b}_{k+1}\big)^T = -W^{2k+1}(\bar{z}_k, \dots, \bar{z}_1)^{-1} \cdot \bar{R}^{2k+1}(\bar{z}_{2k}, \dots, \bar{z}_1) \in N_c$. Combining (2.9) and (2.15), we obtain by (2.10) a first improvement within the cone of approximations of order $2k$, i.e. a family of solution curves with approximation of order $2k+1$ parametrized by $\varepsilon$ and $N[W^{2k+1}(\bar{z}_k, \dots, \bar{z}_1)]$. In figure 4, the affine solution subspace (2.15) is indicated by the green line within $B \times N[\Delta^k]$, implying within the cone a family of approximations of order $2k+1$. Only one of these approximations of order $2k+1$ is depicted in the cone by a green line.

Now, at this point we have to decide, whether only to look for power series solutions of $G[z] = 0$ or to look for $C^{q-2k-1}$ smooth solution manifolds by use of the implicit function theorem.

If only power series solutions are looked for, then we are already done, because it can be shown that from now on, every equation $T^{2k+2}, T^{2k+3}, \dots$ assumes a structure similar to (2.13), with leading operator $W^{2k+2}, W^{2k+3}, \dots$ adopting the surjectivity of $W^{2k+1}(\bar{z}_k, \dots, \bar{z}_1)$, hence ensuring solvability of the system of undetermined coefficients up to infinity. In this manner, the family of approximations of order $2k+1$ is further lifted to approximations of order $2k+2,\ 2k+3, \dots$ .

In some more detail, this aspect is treated in [S3], Corollary 2, where approximate solution curves of a $k$-transversal cone are lifted to arc space $X_\infty$. We do not want to focus on power series solutions in this paper, but note that the key result needed, is obviously given by formula (2.6), supplying knowledge of leading operators up to infinity. Compare also remark 2) in section 3.

If $C^{q-2k-1}$ smooth solution manifolds are looked for, then suppose the complement $N_c$ to be closed, implying by bounded inverse theorem $W^{2k+1}(\bar{z}_k, \dots, \bar{z}_1) \in GL[N_c, \bar{B}]$ and the solution $\bar{b} \in N_c$ of $H(0,\cdot) = 0$ can uniquely be continued within $N_c$ to $\varepsilon \neq 0$ by use of the implicit function theorem.

In some more detail, for every $n := (c_{2k+1}, \dots, c_{k+1})^T \in N[W^{2k+1}(\bar{z}_k, \dots, \bar{z}_1)] =: N$ from (2.15), there exists a locally defined function $(b_{2k+1}, \dots, b_{k+1})\ (\varepsilon, n) \in N_c, |\varepsilon| \ll 1$, of class $C^{q-2k-1}$ with

$$G\left[\ \varepsilon \cdot \bar{z}_1 + \dots + \tfrac{1}{k!}\varepsilon^k \cdot \bar{z}_k \ + \sum_{i=k+1}^{2k+1} \tfrac{1}{i!}\varepsilon^i \cdot (\ \bar{z}_i + b_i(\varepsilon, n)\ )\ \right] = 0 \qquad (2.16)$$

and $\bar{z}_{2k+1} = 0$. In addition, if $n \in N$ is restricted to a bounded domain $\mathcal{B}_\gamma(0) := \{n \in N : \|n\| < \gamma\}$, $\gamma > 0$, then we obtain smoothness according to $(b_{2k+1}, \dots, b_{k+1})\ (\varepsilon, n) \in C^{q-2k-1}(\mathcal{B}_{\delta_1}(0) \times \mathcal{B}_\gamma(0),\ \mathcal{B}_{\delta_2}(\bar{b}))$ with $\mathcal{B}_{\delta_1}(0) := \{\varepsilon \in \mathbb{K} : |\varepsilon| < \delta_1\}$ and $\mathcal{B}_{\delta_2}(\bar{b}) := \{n_c \in N_c : \|n^c - \bar{b}\| < \delta_2\}$ denoting sufficiently small balls in $\mathbb{K}$ and $N_c$ around $\varepsilon = 0$ and $n_c = \bar{b}$ respectively. Note also that $\gamma$ is allowed to be chosen arbitrary large, but finite.

Moreover, within the open and bounded domain $\mathcal{B}_{\delta_1}(0) \times [\mathcal{B}_{\delta_2}(\bar{b}) \oplus \mathcal{B}_\gamma(0)]$, the remainder equation $H(\varepsilon, b_{2k+1},\ b_{2k}, \dots, b_{k+1}) = 0$ from (2.11), (2.15) has no other solutions and summarizing the following existence and uniqueness theorem is shown.

**Theorem 1 :** Assume $G \in C^q(B, \bar{B}),\ q \geq 2k+2$ with $G[0] = 0$ and the existence of $[\bar{z}_{2k}, \dots, \bar{z}_{k+1}, \bar{z}_k, \dots \bar{z}_1]$ satisfying approximation of order $2k$ according to

$$T^{2k}[\ \bar{z}_{2k}, \dots, \bar{z}_{k+1}, \bar{z}_k, \dots \bar{z}_1\ ] = \dots = T^1[\ \bar{z}_1\ ] = 0\,, \qquad (A_1)$$

as well as nondegeneracy of order $k$ by

$$R[\ W^{2k+1}(\bar{z}_k, \dots, \bar{z}_1)\ ] = \bar{B}\,. \qquad (N_1)$$

Then, for $\varepsilon \in \mathcal{B}_{\delta_1}(0)$ and $(b_{2k+1}, \dots, b_{k+1}) \in \mathcal{B}_{\delta_2}(\bar{b}) \oplus \mathcal{B}_\gamma(0)$, all solutions of $G[z] = 0$ within the corresponding cone of approximations of order $2k$ from (2.10) are given by

$$z(\varepsilon, n) = \varepsilon \cdot \bar{z}_1 + \dots + \tfrac{1}{k!}\varepsilon^k \cdot \bar{z}_k \ + \sum_{i=k+1}^{2k+1} \tfrac{1}{i!}\varepsilon^i \cdot (\ \bar{z}_i + b_i(\varepsilon, n)\ ) \qquad (I_1)$$

with $b_i(\varepsilon, n) \in C^{q-2k-1}(\ \mathcal{B}_{\delta_1}(0) \times \mathcal{B}_\gamma(0),\ \mathcal{B}_{\delta_2}(\bar{b})\ )$ and $(b_{2k+1}, \dots, b_{k+1})\ (0, n) = \bar{b} \in N_c$.

From $(I_1)$ we see that the solution curves agree with the initial approximation $z_0(\varepsilon) = \varepsilon \cdot \bar{z}_1 + \cdots + \frac{1}{(2k)!}\varepsilon^{2k} \cdot \bar{z}_{2k}$ (center line of cone) only with respect to the $\varepsilon$-derivatives $\bar{z}_1, \ldots, \bar{z}_k$ implying identity of order $k$ with respect to $\varepsilon$, analogously to condition $(I_0)$ from the introduction.

In advance, we also note that the operator $W^{2k+1}(\bar{z}_k, \ldots, \bar{z}_1)$ in $(N_1)$ merely depends on first $k+1$ coefficients $G_0^1, \ldots, G_0^{k+1}$ of the Taylor expansion of $G$ at $z = 0$. Again, this corresponds exactly to condition $(N_0)$ from the introduction, where first $k+1$ derivatives of $G$ have to be calculated in the origin for deciding about the validity of the $det$-condition $(N_0)$.

Finally, both of the approximation conditions $(A_1)$ and $(A_0)$ depend on first $2k$ derivatives of $G$ in the origin, in this way yielding conditions of order $2k$.

**Remark** : It should be noted that instead of $T^{2k+1}[\,\cdot\,] + O(|\varepsilon|) = 0$ in (2.11), we can also work with the remainder equation

$$T^{2k}[\, z_{2k}, \ldots, z_{k+1},\, z_k,\, \bar{z}_{k-1}, \ldots, \bar{z}_1 \,] + O(\,|\varepsilon|\,) = 0\,,$$

implying solutions of the form

$$z(\varepsilon, n) = \varepsilon \cdot \bar{z}_1 + \cdots + \frac{1}{(k-1)!}\varepsilon^{k-1} \cdot \bar{z}_{k-1} + \frac{1}{k!}\varepsilon^k \cdot (\,\bar{z}_k + b_k(\varepsilon, n)\,) + \sum_{i=k+1}^{2k} \frac{1}{i!}\varepsilon^i \cdot (\,\bar{z}_i + b_i(\varepsilon, n)\,)\,,$$

which represent a certain generalization of Theorem 1, because the coefficient $z_k = \bar{z}_k + b_k(0, n)$ is now allowed to vary in the vicinity of $\bar{z}_k$. But also when performing this extension, first a basic solution of $T^{2k+1} = 0$ is constructed, which is further continued to $\varepsilon \neq 0$ using surjectivity of the same operator $W^{2k+1}(\bar{z}_k, \ldots, \bar{z}_1)$.

In spite of Theorem 1, a lot of questions remain open. First, the cone of approximations within the Banach space $B$ is only poorly described by the image of the smooth operator in (2.10) mapping $(\varepsilon, b_{2k+1}, \ldots, b_{k+1}) \in \mathcal{B}_{\delta_1}(0) \times \left[\mathcal{B}_{\delta_2}(\bar{b}) \oplus \mathcal{B}_\gamma(0)\right] \subset \mathbb{K} \times B^{k+1}$ into $B$. In fact, it is even not clear, whether each half-cone defines an open domain in $B$ around the center line $z_0(\varepsilon)$.

Secondly, the family $z(\varepsilon, n)$ of solution curves obviously shows a strong redundancy due to possible reparametrizations of a solution curve, e.g. by $\varepsilon \to \varepsilon + \varepsilon^{k+1}$, running through the same solution orbit in $B$, but with different parametrization within the family $z(\varepsilon, n)$. It would be desirable to have some sort of minimal parametrization that covers all solutions in the cone.

Thirdly, Theorem 1 does not at all show the iterative aspect of successively increasing $k$ until the direct sum $R_1 \oplus R_2 \oplus \cdots \oplus R_{k+1} \subset \bar{B}$ satisfies $R_1 \oplus R_2 \oplus \cdots \oplus R_{k+1} = \bar{B}$ from (1.3), i.e. the link between the surjectivity condition $R[W^{2k+1}(\bar{z}_k, \ldots, \bar{z}_1)] = \bar{B}$ and the fine resolution of a $k$-transversal cone is not yet performed.

And finally, it has not become clear that the cone with corresponding fibration can be constructed without supposing any approximation properties. In the next section, our aim is to fill these gaps. Nevertheless, we think it is worth noting Theorem 1 separately, due to its remarkable simplicity of proof.

## *3. Main Results*

Now, to answer these questions, take $k$ leading coefficients $[\bar{z}_1, \ldots \bar{z}_k] = [0, \ldots 0, \bar{z}_l, \ldots \bar{z}_k] \in B^k$ of an arbitrary curve $z_0(\varepsilon) = \varepsilon \cdot \bar{z}_1 + \cdots + \frac{1}{k!}\varepsilon^k \cdot \bar{z}_k + O\left(|\varepsilon|^{k+1}\right)$ satisfying $\bar{z}_l \neq 0,\ 1 \leq l \leq k$ and define a sequence of linear operators $S_1, S_2, \ldots, S_{k+1}$ with corresponding direct sums along the following lines of reasoning.

As preliminary, set $N_0 := B$, $R_0^c := \bar{B}$ and $S_1 := G_0^1 \in L[N_0, R_0^c]$ with closed kernel $N_1 := N[S_1] \subset N_0$, range $R_1 := R[S_1] \subset R_0^c$ and associated direct sum decompositions of $B$ and $\bar{B}$ according to

$$\begin{array}{ccccc} B = N_0 = & N_1^c & \oplus & N_1 \\ & \uparrow & & \\ & \boxed{S_1} & & \\ & \downarrow & & \\ \bar{B} = R_0^c = & R_1 & \oplus & R_1^c \end{array} \tag{3.1}$$

All subspaces are assumed to be closed with continuous projections. The bijection $S_1 = G_0^1$ between $N_1^c$ and $R_1$ is indicated in (3.1) by arrows. In addition, choose the complement $N_1^c$ such that $\bar{z}_l \notin N_1^c$. Next define

$$S_2 := 2\, P_{R_1^c}\, G_0^2 \bar{z}_1|_{N_1} \in L[\, N_1, R_1^c\,] \tag{3.2}$$

with continuous projection $P_{R_1^c}$ evaluated with respect to decomposition (3.1). Note that $S_2$ is mapping from the kernel $N_1$ of $S_1$ into the complement $R_1^c$ of the image $R_1$ of $S_1$, i.e. the second operator $S_2$ is just creating values in the subspace $R_1^c$ that is not reached by the first mapping.

As for $S_1$, we assume for $S_2$ a decomposition of $N_1$ and $R_1^c$ by closed subspaces and continuous projections using kernel $N_2 := N[S_2] \subset N_1$ and range $R_2 := R[S_2] \subset R_1^c$ according to $N_1 = N_2^c \oplus N_2$ and $R_1^c = R_2 \oplus R_2^c$. And again choose $N_2^c$ such that $\bar{z}_l \notin N_2^c$.

We note in advance that if the given curve $z_0(\varepsilon)$ is an approximation of order $2k$, then the construction simplifies somewhat due to $\bar{z}_l \in N_1, N_2, \ldots$, i.e. $\bar{z}_l \neq 0$ is contained in every kernel $N_i$ and no care is necessary with respect to chosen complements $N_i^c$.

Continuing this way up to $S_{k+1} \in L[N_k, R_k^c]$, a sequence of successively defined closed kernels and ranges is constructed satisfying $B = N_0 \supset N_1 \supset \cdots \supset N_k \supset N_{k+1}$, as well as defining direct sum decompositions of $B$ and $\bar{B}$ according to

$$\begin{array}{l} B = \overbrace{N_1^c \oplus N_2^c \oplus \cdots \oplus N_{k+1}^c}^{=:\,N^c} \oplus N_{k+1} \quad \textit{with} \quad N_{k+1} = P_{k+1} \oplus \{\bar{z}_{l,k+1}\} \\ \quad\;\; \uparrow \qquad \uparrow \qquad\qquad \uparrow \\ \quad\; \boxed{S_1} \;\; \boxed{S_2} \qquad \boxed{S_{k+1}} \\ \quad\;\; \downarrow \qquad \downarrow \qquad\qquad \downarrow \\ \bar{B} = R_1 \oplus R_2 \oplus \cdots \oplus R_{k+1} \oplus R_{k+1}^c \end{array} \tag{3.3}$$

Note that construction (3.3) is also possible without assuming closedness of all subspaces and continuity of all projections. But then, only power series solutions may be looked for and the application of the implicit function theorem is not possible. Compare also section 2.

The complete construction merely depends on first $k$ coefficients $\bar{z}_1, \ldots, \bar{z}_k$ of $z_0(\varepsilon)$ and on first $k+1$ derivatives $G_0^1, \ldots, G_0^{k+1}$ of $G$ at $z = 0$. It may well happen to obtain $S_i = 0 \in L[N_{i-1}, R_{i-1}^c]$ for some $1 \leq i \leq k+1$ with corresponding decompositions simplifying according to $N_{i-1} = \{0\} \oplus N_i$ and $R_{i-1}^c = \{0\} \oplus R_i^c$. In case of $G_0^1 = \cdots = G_0^{k+1} = 0$, even the case $\bar{B} = \{0\} \oplus \cdots \oplus \{0\} \oplus R_{k+1}^c$ occurs, which simply means that up to given $k \geq 1$ no progress is made with respect to building up the direct sum $\bar{B} = R_1 \oplus R_2 \oplus \cdots \oplus R_{k+1}$.

Next, considering $\bar{z}_l \notin N^c$, choose the projection of $\bar{z}_l$ onto $N_{k+1}$, i.e. $\bar{z}_{l,k+1} \neq 0$, allowing us to refine the direct sum of $B$ in (3.3) according to $N_{k+1} = P_{k+1} \oplus \{\bar{z}_{l,k+1}\}$ with closed complement $P_{k+1} \subset N_{k+1}$, as already indicated in (3.3). Then, the cone $C_k$ is defined by the map

$$\underbrace{Z_k\,(\varepsilon, n_1^c, \dots,\ n_{k+1}^c, p_{k+1})}_{\in \mathbb{K}\times N_1^c \times\cdots\times N_{k+1}^c \times P_{k+1}} = z_0(\varepsilon) + \underbrace{\left[\frac{\varepsilon^{2k+1}}{(2k+1)!} \cdots \frac{\varepsilon^{k+1}}{(k+1)!}\right] \cdot \overbrace{\begin{pmatrix} I_B & * & * \\ & \ddots & * \\ & & I_B \end{pmatrix}}^{=:\ \boxed{\widehat{M}_{k+1}}}}_{=:\ A_\varepsilon = [\,A_\varepsilon^1 \ \cdots\ A_\varepsilon^{k+1}\,]} \cdot \begin{pmatrix} n_1^c \\ \vdots \\ n_{k+1}^c + p_{k+1} \end{pmatrix}. \tag{3.4}$$

The curve $z_0(\varepsilon)$ defines the center line of the cone with leading term given by $\bar{z}_l$. The second summand is defined by a linear mapping $A_\varepsilon$ satisfying for $\varepsilon \neq 0$

$$A_\varepsilon \in GL[\,N_1^c \times \cdots \times N_k^c \times (N_{k+1}^c \oplus P_{k+1})\,,\ N_1^c \oplus \cdots \oplus N_k^c \oplus N_{k+1}^c \oplus P_{k+1}\,]\,. \tag{3.5}$$

Hence, the range of $A_\varepsilon$ defines by (3.3) for every $\varepsilon \neq 0$ a direct complement to the center line $z_0(\varepsilon)$ implying a simple fibration with base $z_0(\varepsilon)$ in $B$, i.e. at every point of $z_0(\varepsilon)$ a direct complement to $z_0(\varepsilon)$ is attached.

The elements of the matrix operator $\widehat{M}_{k+1} \in GL[B^{k+1}, B^{k+1}]$ are given by the composition of multilinear mappings, mainly derived from (3.3), where the basic structure of the bijection $\widehat{M}_{k+1}$ is given by an upper tridiagonal matrix with diagonal exclusively composed of the identity map, as already indicated in (3.4). The precise definitions are given in section 5.

Then, by restriction of $n_i^c \in N_i^c$ and $p_{k+1} \in P_{k+1}$ to arbitrary open, but bounded domains $U_i^c \subset N_i^c$ and $U_{k+1}^p \subset P_{k+1}$, each comprising zero, the image $U_\varepsilon \coloneqq R[A_\varepsilon(U_1^c \times \cdots \times U_k^c \times (U_{k+1}^c \oplus U_{k+1}^p))]$ is shrinking to zero as $\varepsilon \to 0$, thereby creating by open mapping theorem the open and cone-like neighborhood $C_k$ of $z_0(\varepsilon)$ in $B$ for $\varepsilon \neq 0$ according to

$$C_k = \left\{ z \in B\ :\ z = Z_k(\varepsilon, n^c, p),\ \ n^c \in U_1^c \times \cdots \times U_{k+1}^c,\ \ p_{k+1} \in U_{k+1}^p,\ \ 0 < |\varepsilon| \ll 1 \right\}. \tag{3.6}$$

In some more detail, the operator $Z_k$ delivers by upper triangularity of $\widehat{M}_{k+1}$ in (3.4) a fine structure within the cone caused by the manifolds that occur by individually applying the operator $Z_k$ to the domains $U_1^c \subset N_1^c, \dots, U_{k+1}^c \subset N_{k+1}^c$ and $U_{k+1}^p \subset P_{k+1}$ according to

$$\begin{aligned}
Z_k(\varepsilon, n_1^c, 0, \dots, 0) &= z_0(\varepsilon) + \frac{1}{(2k+1)!}\varepsilon^{2k+1} \cdot n_1^c &&=: Z_k^1(\varepsilon, n_1^c) \\
Z_k(\varepsilon, 0, n_2^c, \dots, 0) &= z_0(\varepsilon) + \frac{1}{(2k)!}\varepsilon^{2k} \cdot n_2^c &&=: Z_k^2(\varepsilon, n_2^c) \\
Z_k(\varepsilon, 0, 0, n_3^c, \dots, 0) &= z_0(\varepsilon) + \frac{1}{(2k+1)!}\varepsilon^{2k-1} \cdot n_3^c + O\big(|\varepsilon|^{2k}\big) &&=: Z_k^3(\varepsilon, n_3^c) \\
\vdots \qquad & \ \ \vdots && \vdots \qquad \vdots \\
Z_k(\varepsilon, 0, \dots, 0, n_{k+1}^c, 0) &= z_0(\varepsilon) + \frac{1}{(k+1)!}\varepsilon^{k+1} \cdot n_{k+1}^c + O\big(|\varepsilon|^{k+2}\big) &&=: Z_k^{k+1}(\varepsilon, n_{k+1}^c) \\
Z_k(\varepsilon, 0, \dots, 0, 0, p_{k+1}) &= z_0(\varepsilon) + \frac{1}{(k+1)!}\varepsilon^{k+1} \cdot p_{k+1} + O\big(|\varepsilon|^{k+2}\big) &&=: Z_k^0(\varepsilon, p_{k+1})
\end{aligned} \tag{3.7}$$

By construction, for $\varepsilon \neq 0$ sufficiently small, the different manifolds meet transversally along the center line $z_0(\varepsilon)$. Further, within $Z_k^{k+1}(\varepsilon, n_{k+1}^c)$ and $Z_k^0(\varepsilon, p_{k+1})$, the variables $n_{k+1}^c$ and $p_{k+1}$ are by leading order multiplied with $\varepsilon^{k+1}$, which denotes shrinking of these manifolds by order of $O\big(|\varepsilon|^{k+1}\big)$ as $\varepsilon \to 0$. The remaining manifolds show faster shrinking between $O\big(|\varepsilon|^{k+2}\big)$ and $O\big(|\varepsilon|^{2k+1}\big)$. The higher order terms within (3.7) cause the rotation of the manifolds around $z_0(\varepsilon)$, as indicated in figure 2.

Note also that the behaviour of the first and the second manifold $Z_k^1(\varepsilon, n_1^c)$ and $Z_k^2(\varepsilon, n_2^c)$ is rather simple, showing no rotation at all. Now, if $k = 1$, then only these manifolds occur in the system and the complexity of the general scheme is not completely seen. The case $k = 1$ is treated in many papers, maybe most prominent [CR].

Now, when mapping $C_k$ from $B$ to $\bar{B}$ by $G$, we will see that the map values $G[z_0(\varepsilon)]$ of the center line are perturbed in the following way

$$G[\,Z_k(\varepsilon, n^c, p_{k+1})\,] = G[\,z_0(\varepsilon)\,] \;+\; \varepsilon^{2k+1} \cdot \hat{L}_{k+1} \cdot \begin{pmatrix} n_1^c \\ \vdots \\ n_{k+1}^c + p_{k+1} \end{pmatrix} + O\big(\,|\varepsilon|^{2k+2}\,\big) \tag{3.8}$$

with linear mapping $\hat{L}_{k+1} \in L[B^{k+1}, \bar{B}]$ satisfying $R\big[\hat{L}_{k+1}|_{N^c}\big] = R_1 \oplus R_2 \oplus \cdots \oplus R_{k+1} \subset \bar{B}$.

Hence, an open neighborhood of $G[z_0(\varepsilon)]$ in $\bar{B}$, created by the leading $\varepsilon^{2k+1}$ term, can only be expected to occur, if $R_{k+1}^c = \{0\}$ is supposed in (3.3), i.e. the surjectivity of $\hat{L}_{k+1}$ with respect to the complement $N^c$ is needed to obtain by leading $\varepsilon$-term a constellation as depicted in figure 2 top-right.

The surjectivity of $\hat{L}_{k+1}$ shows to be equivalent with the surjectivity of $W^{2k+1}(\bar{z}_k, \dots, \bar{z}_1)$ from the last section. Also note that a necessary condition for zeros of $G$ to exist in the cone $C_k$ obviously reads $\|G[z_0(\varepsilon)]\| = O\big(|\varepsilon|^{2k+1}\big)$.

We call the cone $C_k$ from (3.6) a *k-transversal cone with respect to the map G,* if one of the following equivalent conditions is satisfied

$$R\big[\,\hat{L}_{k+1}\,|_{N^c}\big] = \bar{B} \quad\Leftrightarrow\quad R_{k+1}^c = \{0\} \quad\Leftrightarrow\quad \bar{B} = R_1 \oplus R_2 \oplus \cdots \oplus R_{k+1}\,. \tag{$N_2$}$$

The condition $(N_2)$ is now representing nondegeneracy of order $k$. Note that $k$-transversality is a property defined by the image of the cone $C_k$ under the map $G$.

**Theorem 2 :** Assume $G \in C^q(B, \bar{B}), q \geq 2k+2$ with $G[0] = 0$ and a $C^q$-curve

$$z_0(\varepsilon) \;=\; \varepsilon \cdot \bar{z}_1 + \cdots + \tfrac{1}{k!}\,\varepsilon^k \cdot \bar{z}_k \;+\; O\big(\,|\varepsilon|^{k+1}\big) \tag{3.9}$$

with $\bar{z}_l \neq 0,\; 1 \leq l \leq k$, to be the first coefficient different from 0. Then, the following results are valid for $\varepsilon \in \mathbb{K}$ sufficiently small and $k \geq 1$.

(i) $\quad G[\,z\,] = 0,\; z \in C_k \quad\Longrightarrow\quad \|\,G[\,z_0(\varepsilon)\,]\,\| = O\big(\,|\varepsilon|^{2k+1}\,\big)$

Now assume $C_k$ to be a *k-transversal cone* and set $p := p_{k+1}$.

(ii) Suppose approximation of order $2k$ by $\|G[\,z_0(\varepsilon)\,]\| = O\big(\,|\varepsilon|^{2k+1}\,\big)$.

Then all solutions of $G[z] = 0$ in $C_k$ are given by

$$z(\varepsilon, p) := Z_k[\,\varepsilon, n^c(\varepsilon, p), p\,] \;=\; z_0(\varepsilon) \;+\; [\,A_\varepsilon^1 \;\cdots\; A_\varepsilon^{k+1}\,] \cdot \begin{pmatrix} n_1^c(\varepsilon, p) \\ \vdots \\ n_{k+1}^c(\varepsilon, p) + p \end{pmatrix}$$

with $n_i^c(\varepsilon, p) \in C^{q-2k-1}\big(\,\mathcal{B}_{\delta_1}(0) \times U_{k+1}^p, N_i^c\,\big)\,,\; i = 1, \dots k+1$.

The map $z(\varepsilon, p)$ defines a regular Banach manifold in $B$ for $\varepsilon \neq 0$.

(iii) Perturbations of $G[z]$ of order $O\big(\|z\|^{2k+1}\big)$ cannot destroy the family $z(\varepsilon, p)$ of solution curves from (ii), only varying $\varepsilon$-derivatives of $z(\varepsilon, p)$ above the order of $k$.

Perturbations of $G[z]$ of order $O\big(\|z\|^{k+2}\big)$ cannot destroy the cone $C_k$ .

(iv) There exists an $(\varepsilon, p)$-dependent near identity transformation $\mathcal{I}_{\varepsilon,p} = I_{\bar{B}} + O(|\varepsilon|)$ in $\bar{B}$, linearizing the map (3.8) from the cone $C_k \subset B$ to $\bar{B}$ according to

$$\mathcal{I}_{\varepsilon,p} \circ G_{\varepsilon,p} \circ A_\varepsilon = \varepsilon^{2k+1} \cdot \hat{L}_{k+1}$$

with linear bijection $\hat{L}_{k+1} \in GL[N^c, \bar{B}]$, linear bijection $A_\varepsilon \in GL[N^c, N^c]$ given by (3.4) and $G_{\varepsilon,p}$ defined by $G_{\varepsilon,p}[\, z\,] := G[\, z_0(\varepsilon) + z + A_\varepsilon^{k+1} \cdot p\,]$.

(v) The operator norms of partial derivatives of $G$ along $z_0(\varepsilon)$ with respect to the subspaces $R\big[\, A_\varepsilon^i(N_i^c)\,\big],\ i = 1, \dots k+1$ satisfy

$$\|\, G'[\, z_0(\varepsilon)\,]_{\,|\, R[\, A_\varepsilon^1\,]}\, \| = O(\, |\varepsilon|^0\,)\ ,\ \ \cdots\ \ , \|\, G'[\, z_0(\varepsilon)\,]_{\,|\, R[\, A_\varepsilon^{k+1}\,]}\, \| = O\big(\, |\varepsilon|^k\big)\, .$$

The linearization with respect to the complement $N^c$ satisfies for $\varepsilon \neq 0$

$$G_{N^c}[\, z_0(\varepsilon)\,] := G'[\, z_0(\varepsilon)\,]_{\,|\, N^c} \in GL[\, N^c, \bar{B}\,]$$

with operator norm of the inverse limited by

$$\|\, G_{N^c}[\, z_0(\varepsilon)\,]^{-1}\, \| = O\big(\, |\varepsilon|^{-k}\,\big) \quad as \quad \varepsilon \to 0\, .$$

(vi) Suppose $dim(N^c) < \infty$. Then, the characteristic number $\chi \geq 0$ is well defined by

$$\chi := 1 \cdot dim(N_2^c) + \cdots + k \cdot dim(N_{k+1}^c)$$

and at every point $Z(\varepsilon, n^c, p)$ of the cone $C_k$ , the determinant with respect to $N^c$ reads

$$det\{\, G_{N^c}[\, Z(\varepsilon, n_c, p)\,]\,\} = \varepsilon^\chi \cdot r(\varepsilon, n_c, p)$$

with smooth map $r(\varepsilon, n_c, p) \neq 0$.

**Remarks : 1)** Note that only part *(ii) Existence and Uniqueness* as well as part *(iii) Perturbation* assume approximation of $z_0(\varepsilon)$ of order $2k$. Part *(iv) Linearization*, part *(v) Partial Derivatives and Inverse*, as well as part *(vi) Determinant* are valid with respect to an arbitrary curve $z_0(\varepsilon)$ through $z = 0$ satisfying (3.9).

**2)** Theorem 2 (ii) concentrates on the description of all solutions of $G[z] = 0$ in the cone $C_k$ based on a certain parametrization of the cone given by (3.4). If one is interested to obtain all elements in arc space $X_\infty$ belonging to these solutions, i.e. having leading part $z = \varepsilon \cdot \bar{z}_1 + \cdots + \frac{1}{k!}\varepsilon^k \cdot \bar{z}_k$, then it is necessary to solve the system of undetermined coefficients in complete generality up to infinity for fixed $[\bar{z}_1, \dots \bar{z}_k] \in B^k$.

Now, as soon as a *k-transversal cone* is reached, this process can be performed, yielding for $k \geq 1,\ l \geq 1$, the following structure concerning the solutions of the system of undetermined coefficients

$$[\, z_1, \dots, z_k \,] = [\, \bar{z}_1, \dots, \bar{z}_k \,] \tag{3.10}$$

---

$$z_{k+1} = \underbrace{P_1}_{\in N^c} + q_1 \qquad q_1 \in N_{k+1}$$

$$z_{k+2} = \underbrace{P_2(q_1)}_{\in N^c} + q_2 \qquad q_2 \in N_{k+1}$$

$$\vdots \qquad\qquad \vdots$$

$$z_{k+l} = \underbrace{P_l(q_1, \dots, q_{l-1})}_{\in N^c} + q_l \qquad q_l \in N_{k+1}$$

---

$$z_{k+l+1} = \underbrace{P_{l+1}(q_1, \dots, q_l)}_{\in N_1^c \oplus \cdots \oplus N_k^c} + n_k \qquad n_k \in N_k$$

$$z_{k+l+2} = \underbrace{P_{l+2}(q_1, \dots, q_l, n_k)}_{\in N_1^c \oplus \cdots \oplus N_{k-1}^c} + n_{k-1} \qquad n_{k-1} \in N_{k-1}$$

$$\vdots \qquad\qquad \vdots$$

$$z_{k+l+k} = \underbrace{P_{l+k}(q_1, \dots, q_l, n_k, \dots, n_2)}_{\in N_1^c} + n_1 \qquad n_1 \in N_1$$

with operators $P_i(\cdot)$, $2 \le i \le l + k$ defined by the composition of multilinear mappings that can effectively be constructed by the recursion given in section 5. If $B$ and $\bar{B}$ are of finite dimensions, e.g. $G : \mathbb{C}^n \to \mathbb{C}^m$ with $n > m$, then the $P_i(\cdot)$ can be represented by polynomials of increasing degree with respect to chosen coordinates. For more details see [S3].

Note that moving $l$ towards infinity allows to calculate $\bar{X}_\infty$, i.e. the subset of arc space $X_\infty$ with leading coefficients $[\bar{z}_1, \dots, \bar{z}_k]$ according to

$$\bar{X}_\infty = \{\, (z_i)_{i \in \mathbb{N}} \;:\; [\, z_1, \dots, z_k \,] = [\, \bar{z}_1, \dots, \bar{z}_k \,],\ z_{k+l} = P_l(q_1, \dots, q_{l-1}) + q_l\,,\ q_l \in N_{k+1}\,,\ l \in \mathbb{N} \,\}\,.$$

**3)** In (iii), the stability of solutions $z(\varepsilon, p)$ and cone $C_k$ are treated with respect to higher order perturbations of $G[z]$. The results in (iii) are far from optimal. In some more detail, perturbations of $G[z]$ of order $\left\lfloor \frac{2k}{l} \right\rfloor + 1$, $1 \le l \le k$, $\bar{z}_l \ne 0$, cannot destroy the family $z(\varepsilon, p)$ of solution curves only varying $\varepsilon$-derivatives of $z(\varepsilon, p)$ above the order of $k$. The notion $\lfloor q \rfloor$ means next integer less or equal $q \in \mathbb{Q}$. This improvement is not shown in this paper, but follows from a close inspection of the recursion given in section 5. For example, if $\bar{z}_1 = 0$ and $\bar{z}_2 \ne 0$, i.e. $l = 2$, then $\left\lfloor \frac{2k}{2} \right\rfloor + 1 = k + 1$ and perturbations of order $k + 1$ will not destroy $z(\varepsilon, p)$, a good improvement.

Another improvement can be derived by considering the solvability condition $\mathrm{I}^k(\bar{z}_k, \dots, \bar{z}_1) \in R[\Delta^k(\bar{z}_{k-1}, \dots, \bar{z}_1)]$ from (2.4) more closely. Here, the matrix operator $\Delta^k$ depends on $G_0^1, \dots, G_0^k$, the inhomogeneity $\mathrm{I}^k$ depends on $G_0^1, \dots, G_0^{2k}$, whereas the nondegeneracy condition depends on $G_0^1, \dots, G_0^{k+1}$. Hence, any perturbation of order between $k + 2$ and $2k$ that does not move the inhomogeneity $\mathrm{I}^k$ out of the range $R[\Delta^k]$, will not destroy the solutions $z(\varepsilon, p)$. In specific examples and in the context of Newton polygons, this criterion turns out to be quite useful.

Finally, note that perturbations of order between $k + 2$ and $2k$ may destroy the solutions by destroying the approximation condition, but not the cone.

**4)** In (iv), the higher order terms $O(|\varepsilon|^{2k+2})$ of the cone mapping (3.8) are eliminated by an $RL$-transformation, where the coordinate transformation in $B$ is directly given by the linear bijection $A_\varepsilon \in GL[N^c, N^c]$ from (3.10), i.e. the linearization in (iv) is essentially obtained by the left transformation $\mathcal{I}_{\varepsilon,p}$ of the image space $\bar{B}$.

Alternatively, (3.8) may be linearized by pure right transformation along the following lines. Simply consider the equation

$$G[\, Z_k(\varepsilon, n^c, p)\,] \;=\; G[\, z_0(\varepsilon)\,] \;+\; \varepsilon^{2k+1} \cdot \hat{L}_{k+1} \cdot \varphi$$

with new variable $\varphi \in N^c$, or equivalently by (3.8) after cancelling out $\varepsilon^{2k+1} \neq 0$

$$\hat{L}_{k+1} \cdot n^c \;+\; \hat{L}_{k+1}^{k+1} \cdot p \;+\; O(|\varepsilon|) \;-\; \hat{L}_{k+1} \cdot \varphi \;=\; 0$$

where $\hat{L}_{k+1}^{k+1} \in L[B, \bar{B}]$ denotes the last component of $\hat{L}_{k+1}$. Now, by $\hat{L}_{k+1} \in GL[N^c, \bar{B}]$, the equation can uniquely be solved with respect to $n^c \in N^c$, implying linearization of the cone mapping (3.8) according to

$$G[\, Z_k(\, \varepsilon,\, n^c[\, \varepsilon, p, \varphi\,], p\,)\,] \;=\; G[\, z_0(\varepsilon)\,] \;+\; \varepsilon^{2k+1} \cdot \hat{L}_{k+1} \cdot \varphi$$

with new coordinates $\varphi \in N^c$. Depending from the question to be solved, one may decide, which version of linearization to be more appropriate.

**5)** Concerning part (v), we notice that the approximation $z_0(\varepsilon)$ can more generally be replaced by the complete cone $Z(\varepsilon, n_c, p)$, analogously to part (vi). This means that the cone is characterized by homogenous conditions with respect to derivatives. In particular, $G'[Z_k(\varepsilon, n^c, p)] \in GL[N^c, \bar{B}]$ in case of $\varepsilon \neq 0$, i.e. $z = 0$ is at an isolated singularity with respect to the complete cone and the operator norm of the inverse is uniformly limited in the cone by

$$\| \, G_{N^c}[\, Z(\varepsilon, n_c, p)\,]^{-1} \, \| = O\big(\, |\varepsilon|^{-k}\,\big) \quad as \quad \varepsilon \to 0\,. \tag{3.11}$$

This follows straightforward from (3.8) and the proof in section 6.

In the context of an existing trivial solution curve $y = 0$ of $G[\lambda, y] = L(\lambda)y + O(\|y\|^2)$, the notion of $\lambda = 0$ to be a *generalized algebraic eigenvalue of the operator family* $L(\lambda)$ *of order* $k$ is introduced in [L1], if the norm of the inverse of the Jacobian with respect to $y$ is limited according to

$$\| \, G_y[\lambda\,, 0\,]^{-1} \, \| = O\big(\, |\lambda|^{-k}\,\big) \quad as \quad \lambda \to 0\,.$$

Hence by (3.11), we merely extend this concept to curves (or cones) in general position in $B$ and without assuming further Fredholm properties.

Reversely, in Corollary 1 below, we give sufficient conditions ascertaining that an isolated singularity with operator norm of the inverse exploding by $O\big(|\varepsilon|^{-k}\big)$ gives rise to a $k$-transversal cone.

**6)** In (vi), also the case of $\chi = 0$, i.e. $\bar{B} = R_1 = R[G_0^1]$, may well arise. Then, $G_0^1$ is already surjective and we simply meet the constellation of the standard implicit function theorem. In this situation, Theorem 2 (vi) only says that the determinant with respect to the complement $N^c$ is different from zero in the cone, where in fact, this determinant is different from zero everywhere.

In case of $\chi$ even, the sign of the determinant is equal in both of the half-cones and zero at $z = 0$. If $\chi$ is odd, then the sign of the determinant is different in the two half-cones.

In this sense, the characteristic number $\chi$ is a property of a *k-transversal cone*. In the next section, the implications of (vi) are studied in some more detail.

**7)** Let us interpret one of the results from the viewpoint of an operator family. Given an arbitrary curve $z_0(\varepsilon) = \varepsilon \cdot \bar{z}_1 + \cdots$ of the form (3.9) and perform an expansion along this curve according to

$$G[\,z_0(\varepsilon) + z\,] = G[\,z_0(\varepsilon)\,] + G'[\,z_0(\varepsilon)\,] \cdot z + O(\,\|z\|^2\,)\,.$$

Further, assume the $\varepsilon$-dependent operator family $G'[z_0(\varepsilon)] \in L[B,\bar{B}]$ to be singular at $\varepsilon = 0$, i.e. $R[G_0^1] \neq \bar{B}$, but regular, when disturbing to $\varepsilon \neq 0$. Then, one may ask, how many of the summands within the $\varepsilon$-expansion of the linear part

$$G'[\,z_0(\varepsilon)\,] \cdot z = \underbrace{\{\,G_0^1 + \varepsilon \cdot G_0^2 \bar{z}_1 + \tfrac{1}{2}\varepsilon^2 \cdot (\,G_0^2 \bar{z}_2 + G_0^3 \bar{z}_1^2\,) + \cdots + \tfrac{1}{k!}\varepsilon^k \cdot (\,\cdots\,) + \cdots\}}_{=:\,P^k(\varepsilon)} \cdot z \qquad (3.12)$$

are in fact needed for ensuring regularity of $G'[z_0(\varepsilon)]$ for $\varepsilon \neq 0$? This sort of question remembers $k$-jet determinacy (of an operator family). Now, if a *$k$-transversal cone* $C_k$ can be build up according to (3.3), (3.4), $(N_2)$, then by application of the singular transformation $z = A_\varepsilon \cdot n^c$ to (3.12), we obtain

$$G'[\,z_0(\varepsilon)\,] \cdot A_\varepsilon \cdot n^c = \varepsilon^{2k+1} \cdot \hat{L}_{k+1} \cdot n^c + O\big(\,|\varepsilon|^{2k+2}\,\big)$$

with regularity for $\varepsilon \neq 0$ due to $\hat{L}_{k+1} \in GL[N^c,\bar{B}]$. Hence, the effect of the singular transformation, i.e. the effect of the cone, is to concentrate the regularity of the operator family $G'[z_0(\varepsilon)]$ for $\varepsilon \neq 0$ in one leading operator $\hat{L}_{k+1}$, where this process involves exactly first $k+1$ summands within the expansion (3.12), i.e. the Taylor polynomial $P^k(\varepsilon)$ of degree $k$ of the smooth operator family $G'[\,z_0(\cdot)\,] \in C^{q-1}\big(\mathcal{B}_{\delta_1}(0), L[B,\bar{B}]\big)$ has to be taken into consideration.

In addition, the recursion of building up a *$k$-transversal cone* $C_k$ represents a constructive method for calculation of the order $k$ given in (3.11). Maybe, this concept is applicable to further families of operators, supposed the singularity is isolated and of non-flat type.

In the last part of the results section, we state some sufficient conditions for a curve $z_0(\varepsilon)$ to be surrounded by a *submersive $k$-transversal cone.* But first, let us collect some necessary conditions following from the existence of a *$k$-transversal cone* $C_k$.

Concerning the level sets of the submersion $G[z]$ in the cone $C_k$ with respect to the curve $z_0(\varepsilon)$, we will obtain the equivalence

$$G[\,z\,] = G[\,z_0(\varepsilon)\,]\,, z \in C_k \qquad (3.13)$$

$$\Leftrightarrow \quad z = z_0(\varepsilon) + \underbrace{[\,A_\varepsilon^1 \;\cdots\; A_\varepsilon^{k+1}\,] \cdot \begin{pmatrix} \bar{n}_1^c(\varepsilon, n_{k+1}) \\ \vdots \\ \bar{n}_{k+1}^c(\varepsilon, n_{k+1}) + n_{k+1} \end{pmatrix}}_{=:\,\varepsilon^{k+1} \cdot \bar{A}(\varepsilon, n_{k+1})}$$

with smooth mappings $\bar{n}_i^c(\varepsilon, n_{k+1}) \in C^{q-2k-1}\big(\mathcal{B}_{\delta_1}(0) \times U_{k+1}, N_i^c\big), U_{k+1} \subset N_{k+1},\; i = 1, \dots k+1$.

Thus, for every $\varepsilon \neq 0$ fixed, the set in $B$ mapped to $G[z_0(\varepsilon)]$ can locally be parametrized by a subset $U_{k+1}$ of the smallest kernel $N_{k+1}$ of the filtration (blue box in figure 1). In addition, we have $\bar{A}(0, n_{k+1}) = n_{k+1}$, i.e. when passing from $\varepsilon = 0$ to $\varepsilon \neq 0$, the level sets are given by scaled perturbations of $N_{k+1}$. In some more detail, differentiation of the level set identity

$$G[\,z_0(\varepsilon) + \varepsilon^{k+1} \cdot \bar{A}(\varepsilon, n_{k+1})\,] = G[\,z_0(\varepsilon)\,]$$

in the direction $N_{k+1}$ of the level sets implies

$$G'[\, z_0(\varepsilon)\,] \cdot [\, I_B + O(\,|\varepsilon|\,)\,] \cdot n_{k+1} = 0 \tag{3.14}$$

and we see that a necessary condition for *k-transversality* of a cone $C_k$ requires the kernel $N_{k+1} \subset N_1 = N\big[G'[\, z_0(0)\,]\big] = N[G_0^1]$ to be smoothly continued from $\varepsilon = 0$ to $\varepsilon \neq 0$.

Now, collecting necessary conditions of a *k-transversal cone*, we obtain from (3.3) the split condition $B = N^c \oplus N_{k+1}$ with operator norm of the $N^c$-inverse limited by (3.11) and $N_{k+1}$ satisfying (3.14). Essentially, these properties are also sufficient for *k-transversality* of a given cone $C_k$ according to the following Corollary.

**Corollary 1 :** Given a curve $z_0(\varepsilon)$ of the form (3.9) with corresponding cone $C_k$ from (3.6) and decomposition (3.3) with $R_{k+1}^c = \{0\}$. Then (i), (ii) and (iii) are equivalent.

(i) The cone $C_k$ is *k-transversal*

(ii) There exist closed subspaces $U \neq \{0\},\ V \neq \{0\}$ with $B = U \oplus V$ satisfying $G_U[z_0(\varepsilon)] := G'[z_0(\varepsilon)]_{|U} \in GL[U, \bar{B}], \varepsilon \neq 0$ and

$$\|\, G_U[\, z_0(\varepsilon)\,]^{-1} \,\| = O\big(\,|\varepsilon|^{-k}\,\big) \tag{3.15}$$

$$G'[\, z_0(\varepsilon)\,] \cdot [\, I_B + \varepsilon \cdot Q(\varepsilon)\,] \cdot v \;= 0,\;\; v \in V \tag{3.16}$$

with $Q(\varepsilon) \in C\big(\mathcal{B}_{\delta_1}(0), L[V, B]\big)$.

(iii) There exist closed subspaces $U \neq \{0\},\ V \neq \{0\}$ with $B = U \oplus V$ satisfying $P_U^k(\varepsilon) := P^k(\varepsilon)_{|U} \in GL[U, \bar{B}], \varepsilon \neq 0$ and

$$\|\, P_U^k(\varepsilon)^{-1} \,\| = O\big(\,|\varepsilon|^{-k}\,\big) \tag{3.17}$$

$$\|\, P^k(\varepsilon) \cdot [\, I_B + \varepsilon \cdot Q(\varepsilon)\,]_{|V} \,\| \; = O\big(\,|\varepsilon|^{k+1}\,\big) \tag{3.18}$$

with $Q(\varepsilon) \in C\big(\mathcal{B}_{\delta_1}(0), L[V, B]\big)$.

Note that according to (iii), *k-transversality* of a cone $C_k$ merely depends on the Taylor polynomial $P^k(\varepsilon)$ of degree $k$ of the operator family $G'[z_0(\varepsilon)]$. The subspaces $U$, $V$ and the family of operators $Q(\varepsilon)$ in (ii), (iii) can be chosen identically.

## *4. Secondary Bifurcation by Degree Theory*

In this section, we impose further restrictions on $G[z] = 0$ for obtaining further properties besides existence, uniqueness and stability of solution curves as stated in Theorem 2.

First, within the fine resolution (3.3) of cone $C_k$ , we assume $N^c = N_1^c \oplus N_2^c \oplus \cdots \oplus N_{k+1}^c$ to be of finite dimension, implying by use of the characteristic number $\chi \geq 0$ from Theorem 2 (vi) for each half-cone constant conditions with respect to the sign of the determinant. Hence, the application of degree theory will be possible, eventually proving secondary bifurcation from $z = 0$.

For example, if we start with an equation $H[u] = 0,\ H \in C^q(D, \bar{D}),\ D, \bar{D}$ Banach spaces, with $H[0] = 0$ and linearization $H'[0] \in L[D, \bar{D}]$ to be a Fredholm operator (of arbitrary index), then after performing a Lyapunov-Schmidt reduction a reduced system of the form $G[z] = 0,\ G \in C^q(B, \bar{B})$ occurs with finite dimensions of the Banach spaces $B$, $\bar{B}$, implying finite dimension of $N^c$ as well.

Further, assume the curve $z_0(\varepsilon)$ from (3.9) to be an approximation of order $2k$, as well as the cone $C_k$ to be $k$-transversal according to a simplified decomposition (3.3) of the form

$$\begin{array}{l} \overbrace{\hphantom{N_1^c \oplus N_2^c \oplus \cdots \oplus N_{k+1}^c}}^{dim(N^c)<\infty} \\ B = N_1^c \oplus N_2^c \oplus \cdots \oplus N_{k+1}^c \oplus N_{k+1} \quad with \quad N_{k+1} = \{0\} \oplus \{\bar{z}_l\} \\ \quad\;\; \uparrow \qquad \uparrow \qquad\qquad \uparrow \\ \quad\; \boxed{S_1} \;\; \boxed{S_2} \qquad \boxed{S_{k+1}} \\ \quad\;\; \downarrow \qquad \downarrow \qquad\qquad \downarrow \\ \bar{B} = R_1 \oplus R_2 \oplus \cdots \oplus \underbrace{R_{k+1}}_{\neq 0} \end{array} \tag{4.1}$$

In particular, by $R_{k+1} \neq \{0\}$, $k$ is supposed to be minimal in the sense that for smaller values of $k$ the Banach space $\bar{B}$ is not yet build up by the direct sum of the ranges. Moreover, the characteristic number $\chi$ satisfies $\chi \geq k$ by Theorem 2 (vi) under consideration of $N_{k+1}^c \neq \{0\}$. Finally, due to $P_{k+1} = \{0\}$, no extra parameters occur in the system.

Thus, by Theorem 2 (ii), precisely one solution curve $z(\varepsilon) = z_0(\varepsilon) + A_\varepsilon \cdot n^c(\varepsilon)$, $A_\varepsilon \in GL[N^c, N^c]$ of $G[z] = 0$ exists in $C_k$, implying that the topological degree with respect to the map $\bar{G}_\varepsilon[n^c] \coloneqq G[z_0(\varepsilon) + n^c]$, $n^c \in U_\varepsilon \subset N^c$, $|\varepsilon| \ll 1$, is given by

$$deg(\, \bar{G}_\varepsilon,\, U_\varepsilon \,) = sign\; det\{\, G_{N^c}[\, z_0(\varepsilon)\,]\,\} \tag{4.2}$$

under the restriction $\mathbb{K} = \mathbb{R}$. Hence, by Theorem 2 (vi), we have constant degree $+1$ or $-1$ in each half-cone, ending up with the following result.

**Corollary 2 :** Suppose $\mathbb{K} = \mathbb{R}$ and a curve $z_0(\varepsilon)$ with $\bar{z}_l \neq 0$, $1 \leq l \leq k$ to be an approximation of order $2k$, as well as the cone $C_k$ to be *k-transversal* according to (4.1). Then, if both, $l \geq 1$ and $\chi \geq 1$ are *odd*, a continuum of solutions of $G[z] = 0$ emanates from $z = 0$, existing outside the cone $C_k$.

Note that by (4.2), assumption $\chi$ *odd* and Theorem 2 (vi), the degree changes sign when moving from one half-cone to the other, and by $l$ *odd*, we obtain a well-defined straight movement in $B$ in the direction of $\bar{z}_l$. Hence, a cylindrical construction, as depicted in figure 3, can be set up in a straightforward way to ensure the continuum of solutions by homotopy invariance of the degree. In detail, a construction of this kind is performed in [L1] in the context of a known solution curve $y = 0$ of $G[\lambda, y] = L(\lambda)y + O(\|y\|^2)$.

Compared to the results of the Madrid school, we do not think that much new is given by Corollary 2 from a theoretical point of view. On the other hand, our procedure is working with respect to solution curves in general position within Banach space $B$, i.e. there is no need to split the Banach space $B$ into a product space $B \simeq \Lambda \times Y$ with $dim(\Lambda) = 1$ and/or to transfer the solution curve to an axis.

In addition, using the external parameter $\varepsilon$ for parametrization of the solution curve, the characteristic number $\chi$ and the topological degree in each half-cone can be calculated with respect to solution curves that cannot be moved to an axis. Compare a typical cusp curve, as depicted in figure 3, arising quite frequently in singularity theory.

Let us now look at a simple example for application of Corollary 2.

**Example:** Given the real polynomial equation

$$G[x,y] = -xy^3 + x^5 + y^6 = 0 \tag{4.3}$$

composed of monomials of order 4, 5 and 6. Then, a corresponding solution curve of the form

$$y(x) = x^{4/3} \cdot [\, 1 + r(x) \,] \tag{4.4}$$

is directly obtained by plugging the ansatz $y = x^{4/3}[\, 1 + r \,]$ into (4.3), splitting off $x^5$ and solving the remainder with respect to the new variable $r$ by use of the implicit function theorem.

Let us now repeat the derivation of this solution curve, using the Banach space procedure of this paper, with the aim to investigate the solution curve with respect to secondary bifurcation. First, when neglecting the highest order term $y^6$ in (4.3), we obtain the equation $x(-y^3 + x^4) = 0$, which motivates to start with the approximation

$$z_0(\varepsilon) = \begin{pmatrix} x \\ y \end{pmatrix} = \begin{pmatrix} \varepsilon^3 \\ \varepsilon^4 \end{pmatrix} = \tfrac{1}{3!}\varepsilon^3 \cdot \begin{pmatrix} 3! \\ 0 \end{pmatrix} + \tfrac{1}{4!}\varepsilon^4 \cdot \begin{pmatrix} 0 \\ 4! \end{pmatrix} =: \tfrac{1}{3!}\varepsilon^3 \cdot \bar{z}_3 + \tfrac{1}{4!}\varepsilon^4 \cdot \bar{z}_4 \,.$$

Then, constructing fine resolution and cones by (3.3), (3.4), we reach $k$-transversality by straightforward calculation according to

$$\begin{array}{ccccccccc}
B = \mathbb{R}^2 = & \binom{0}{0} & \oplus \cdots \oplus & \binom{0}{0} & \oplus & \overbrace{\binom{0}{1}}^{= N_{12}^c} & \oplus & \overbrace{\binom{3!}{0}}^{= \bar{z}_3} \\
& \uparrow & & \uparrow & & \uparrow & & \\
& \boxed{[0\ 0]} & & \boxed{[0\ 0]} & & \boxed{S_{12} = [0\ \ c]} & & \\
& \downarrow & & \downarrow & & \downarrow & & \\
\bar{B} = \mathbb{R} \ = & \{0\} & \oplus \cdots \oplus & \{0\} & \oplus & \{1\} & &
\end{array}$$

with $c \neq 0$. Hence, we obtain transversality of order $k = 11$, characteristic number given by $\chi = 11$ according to Theorem 2 (vi), as well as index of leading coefficient given by $l = 3$. Further, the behaviour of the determinant along the approximation $z_0(\varepsilon)$ reads by Theorem 2 (vi)

$$det\{\, G_y[\, z_0(\varepsilon) \,] \,\} = \varepsilon^{11} \cdot s(\varepsilon) \ , \ \ s(\varepsilon) \neq 0,$$

and $z_0(\varepsilon)$ satisfies approximation of order $2k$ by

$$|\, G[\, z_0(\varepsilon) \,] \,| = |\, y^6{}_{\,|\, y=\varepsilon^4} \,| = |\, \varepsilon^{24} \,| = O\big(\, |\varepsilon|^{2k+1} \,\big)$$

under consideration of $2k + 1 = 23$. Thus, we can be sure by Theorem 2 (ii) that a unique solution curve $z(\varepsilon)$ exists within the cone $C_{11}$, in this way repeating (4.4) by external $\varepsilon$-parametrization. Finally, due to $\chi = 11$ *odd* and $l = 3$ *odd*, Corollary 2 ascertains the existence of a continuum of secondary solutions emanating from $z = 0$ outside the cone $C_{11}$.

Again by direct calculation, these secondary solutions are given by $x(y) = y^3 \cdot [1 + p(y)]$ and repeating the $\varepsilon$-procedure from above with respect to $x(y)$, we obtain $k = \chi = 3$ *odd* and $l = 1$ *odd*, yielding now by Corollary 2 the solution curve $y(x)$ from (4.4) as secondary bifurcation. In principle, each curve of a bifurcation diagram can be used as a starting point for investigation of secondary bifurcation by Corollary 2.

Note also, that the Milnor number $\mu$ of the singularity (4.3) is easily calculated by the $k$ or $\chi$ values $k_1 = \chi_1 = 11$ and $k_2 = \chi_2 = 3$ of the two solution curves according to the formula

$$\mu = \chi_1 + \chi_2 - ord(G) + 1 = 11 + 3 - 4 + 1 = 11\,.$$

In general, as derived in [S3], the formula

$$\mu = \chi_1 + \cdots + \chi_\tau - ord(G) + 1$$

is valid, if all segments of the Newton polygon factorize completely with multiplicities 1. Here, $\tau$ denotes the number of different solution curves through the singularity. See also [S3] for further applications of Corollary 2 in higher dimensions, as well as applications to simple $ADE$-singularities [AGV].

Finally, we note that we could also try to calculate the *generalized algebraic multiplicity* $\chi$ associated to $y(x)$ from (4.4) by first moving $y(x)$ to the $x$-axis according to

$$\bar{G}[\,x, \bar{y}\,] := G[\,x,\; x^{4/3} \cdot [\,1 + r(x)] + \bar{y}\,]\,,$$

yielding the trivial solution $\bar{y} = 0$ with respect to $\bar{G}[x, \bar{y}] = 0$ as desired. However, we lose smoothness according to $\bar{G} \notin C^2$ and the standard transformation of $y(x)$ to the $x$-axis does not work for calculation of the generalized algebraic multiplicity $\chi$.

We do not say that it is impossible to transform $y(x)$ to an axis by choosing appropriate coordinates within $B = \mathbb{R}^2$; we only say that the standard transformation does not work, supporting the idea that procedures concerning general position might be useful.

## *5. Construction of Fine Resolution and Cones*

In this section, we show in detail the construction of the fine resolution (3.3) and the cone (3.4). The construction is kept independent of any approximation properties of the curve $z_0(\varepsilon)$ from (3.9). In this sense, a cone-like fibration of a curve satisfying $\bar{z}_l \neq 0, 1 \leq l \leq k$, is constructed.

Using $S_1 = G_0^1$ and $S_2 = P_{R_1^c}\, 2G_0^2 \bar{z}_1|_{N_1}$ from (3.2), we start with $k = 1$ and corresponding decomposition

$$\begin{array}{ccccccc} B = & N_1^c & \oplus & N_2^c & \oplus & N_2 & \quad with \quad N_2 = P_2 \oplus \{\bar{z}_{l,2}\} \\ & \uparrow & & \uparrow & & & \\ & \boxed{S_1} & & \boxed{S_2} & & & \\ & \downarrow & & \downarrow & & & \\ \bar{B} = & R_1 & \oplus & R_2 & \oplus & R_2^c & \end{array} \tag{5.1}$$

Then, the first cone $C_1$ is defined by the map

$$Z_1(\varepsilon, n_1^c, n_2^c, p_2) := z_0(\varepsilon) + [\; \tfrac{1}{3!}\varepsilon^3 \quad \tfrac{1}{2!}\varepsilon^2\,] \cdot \overbrace{\begin{pmatrix} I_B & \\ & I_B \end{pmatrix}}^{=:\, \boxed{\hat{M}_2}} \cdot \begin{pmatrix} n_1^c \\ n_2^c + p_2 \end{pmatrix}$$

$$= z_0(\varepsilon) + \tfrac{1}{2}\varepsilon^2 \cdot (n_2^c + p_2) + \tfrac{1}{6}\varepsilon^3 \cdot n_1^c\,, \tag{5.2}$$

where the complete construction merely depends on the first coefficient $\bar{z}_1$ and the first two derivatives $G_0^1$ and $G_0^2$. In figure 5 two possible realizations of the cone $C_1$ are drawn.

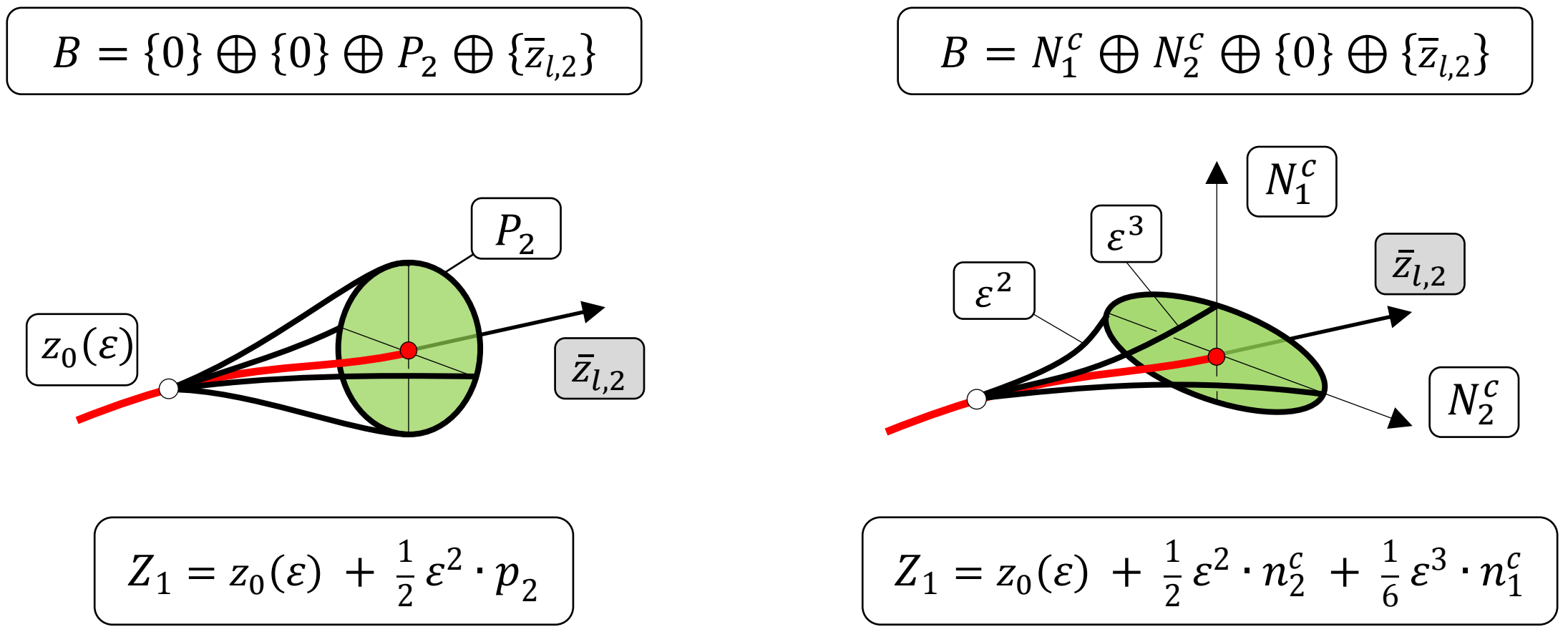

Figure 5 : The cone $C_1$ in two different realizations.

In the left diagram, both $N_1^c$ and $N_2^c$ are zero and the cone is exclusively build up by the subspace $P_2$, which is homogenously shrinking by order of $\varepsilon^2$ in each direction as $\varepsilon \to 0$. On contrary, in the right diagram, $P_2$ equals zero and the cone is composed of $N_1^c$ and $N_2^c$ shrinking by order of $\varepsilon^3$ and $\varepsilon^2$ respectively.

Now, for general $k \geq 2$, the decomposition of $B$, $\bar{B}$ and the cone $C_k$ are defined by the following recursion starting with $k = 2$. Remember, the basis for the construction is given by the system of undetermined coefficients (2.2), obtained from plugging the power series ansatz $z = \varepsilon \cdot \bar{z}_1 + \frac{1}{2!}\varepsilon^2 \cdot \bar{z}_2 + \cdots$ into the power series expansion of $G$ at $z = 0$. Then, the basic formula (2.3) has to be used intensively, implying complex, but constructive formulas.

First, two simple triangular schemes of the form

| | | | | | | | | | |
|---|---|---|---|---|---|---|---|---|---|
| | $\boxed{d_{0,1} = 1}$ | | | | \| | $\boxed{c_{0,1} = 1}$ | | | |
| $n = 1$ | $d_{1,1} = 1$ | | | | \| | $c_{1,1} = 1$ | | | |
| $n = 1$ | $d_{2,1} = 1$ | $\boxed{d_{2,2} = \frac{3}{2}}$ | | | \| | $c_{2,1} = 1$ | $\boxed{c_{2,2} = 1}$ | | |
| $n = 2$ | $d_{3,1} = 1$ | $d_{3,2} = \frac{4}{3}$ | | | \| | $c_{3,1} = 1$ | $c_{3,2} = \frac{3}{2}$ | | |
| $n = 2$ | $d_{4,1} = 1$ | $d_{4,2} = \frac{5}{4}$ | $\boxed{d_{4,3} = \frac{5}{3}}$ | | \| | $c_{4,1} = 1$ | $c_{4,2} = \frac{4}{2}$ | $\boxed{c_{4,3} = 1}$ | |
| $n = 3$ | $d_{5,1} = 1$ | $d_{5,2} = \frac{6}{5}$ | $d_{5,3} = \frac{6}{4}$ | | \| | $c_{5,1} = 1$ | $c_{5,2} = \frac{5}{2}$ | $c_{5,3} = \frac{5}{3}$ | |
| $n = 3$ | $d_{6,1} = 1$ | $d_{6,2} = \frac{7}{6}$ | $d_{6,3} = \frac{7}{5}$ | $\boxed{d_{6,4} = \frac{7}{4}}$ | \| | $c_{6,1} = 1$ | $c_{6,2} = \frac{6}{2}$ | $c_{6,3} = \frac{5}{2}$ | $\boxed{c_{6,4} = 1}$ |
| | ⋮ | ⋮ | ⋮ | ⋮ | \| | ⋮ | ⋮ | ⋮ | ⋮ |
| | $l = 1$ | $l = 2$ | $l = 3$ | $\cdots$ | \| | $l = 1$ | $l = 2$ | $l = 3$ | $\cdots$ |

are defined column by column according to

$$
\begin{array}{llll}
l \geq 1 & : & \boxed{d_{2l-2,l} := \frac{2l-1}{l}} & \boxed{c_{2l-2,l} := 1} \\
n \geq l & : & d_{2n-1,l} := \frac{2n}{2n+1-l} & c_{2n-1,l} := c_{2n-2,l} \cdot d_{2n-2,l} \\
n \geq l & : & d_{2n,l} := \frac{2n+1}{2n+2-l} & c_{2n,l} := c_{2n-1,l} \cdot d_{2n-1,l}
\end{array}
\tag{5.3}
$$

with the rows of the schemes partly comprised to diagonal matrices

$$
\begin{array}{llll}
D^{2n-1} & := & Diag[\, d_{2n-1,1}, \cdots, d_{2n-1,n} \,] & \in \mathbb{R}^{n,n} \\
D^{2n} & := & Diag[\, d_{2n,1}, \cdots, d_{2n,n} \,] & \in \mathbb{R}^{n,n} \\
C^{2n-1} & := & Diag[\, c_{2n-1,1}, \cdots, c_{2n-1,n} \,] & \in \mathbb{R}^{n,n} \\
C^{2n} & := & Diag[\, c_{2n,1}, \cdots, c_{2n,n} \,] & \in \mathbb{R}^{n,n} \,.
\end{array}
\tag{5.4}
$$

Every element of the $d$-scheme is explicitly defined, whereas the $c$-scheme is iteratively defined using the $d$-scheme. With respect to proofs, this kind of $c$-scheme definition is quite convenient. Nevertheless, the $c$-scheme can also be expressed explicitly using binomial coefficients according to

$$
c_{2n-1,l} \cdot \binom{2(l-1)}{l-1} = \binom{2n-1}{l-1} \qquad and \qquad c_{2n,l} \cdot \binom{2(l-1)}{l-1} = \binom{2n}{l-1}.
$$

Then, using the abbreviations

$$
E^2 := \begin{pmatrix} E_{1,1} & E_{1,2} \\ E_{2,1} & E_{2,2} \end{pmatrix} := \begin{pmatrix} \boxed{I_B} & -S_1^{-1} P_{R_1} \bar{S}_2 \\ 0 & \boxed{I_B} \end{pmatrix} \in L[\, B^2, B^2 \,]
$$

$$
\begin{pmatrix} a^3 & \bar{a}^3 \\ A^3 & \bar{A}^3 \end{pmatrix} := (C^3)^{-1} \cdot E^2 \cdot C^3 \in L[\, B^2, B^2 \,]
\tag{5.5}
$$

$$
[\, W_3^4, \ W_2^4 \,](\bar{z}_1) := [\, 4G_0^2 \bar{z}_1, \ 6G_0^3 \bar{z}_1^2 \,] \in L[\, B^2, \bar{B} \,] \qquad and \qquad M^3 := I_B \in L[\, B, B \,],
$$

the recursion, continued with $k = 2$, is well defined according to the following notation

$$
\bar{S}_{k+1} := [\, W_{2k-1}^{2k}, \dots, W_k^{2k} \,](\bar{z}_{k-1}, \dots, \bar{z}_1) \cdot \begin{pmatrix} & \bar{a}^{2k-1} \\ M^{2k-1} \cdot & \bar{A}^{2k-1} \end{pmatrix} + \frac{(2k)!}{(k!)^2} G_0^2 \bar{z}_k \in L[\, B, \bar{B} \,]
$$

$$
S_{k+1} := P_{R_k^c} \bar{S}_{k+1} \,|_{N_k} \in L[\, N_k, R_k^c \,]
\tag{5.6}
$$

with corresponding decomposition

$$
\begin{array}{lllll}
N_k = N_{k+1}^c & \oplus & N[\, S_{k+1}(\bar{z}_k, \dots, \bar{z}_1) \,] & =: & N_{k+1}^c \oplus N_{k+1} \\
R_k^c = R[\, S_{k+1}(\bar{z}_k, \dots, \bar{z}_1) \,] & \oplus & R_{k+1}^c & =: & R_{k+1} \oplus R_{k+1}^c
\end{array}
$$

and all subspaces assumed to be closed and all projections assumed to be continuous.

Hence, setting $S_0 := 0$, the linear mappings $S_i$, $i = 1, \dots, k+1$, are iteratively defined in such a way that the mapping $S_i$ acts between the kernel of the previous mapping $N[S_{i-1}] = N_{i-1}$ and a complement $R^c_{i-1}$ of the range $R[S_{i-1}] = R_{i-1}$ of this mapping.

By construction, we end up with the following decompositions of $B$ and $\bar{B}$ at iteration step $k$.

$$
\begin{array}{l}
B = N^c_1 \oplus \overbrace{N^c_2 \oplus \cdots \oplus \overbrace{N^c_{k+1} \oplus N_{k+1}}^{=N_k}}^{=N_1 \;\; (\vdots)} \quad (\text{whole } N^c_1 \oplus \cdots \oplus N_{k+1} = N_0) \qquad \text{with} \quad N_{k+1} = P_{k+1} \oplus \{\bar{z}_{l,k+1}\} \\
\quad\;\; \uparrow \qquad \uparrow \qquad\qquad\;\; \uparrow \\
\quad\; \boxed{S_1} \;\; \boxed{S_2} \qquad\quad \boxed{S_{k+1}} \\
\quad\;\; \downarrow \qquad \downarrow \qquad\qquad\;\; \downarrow \\
\bar{B} = R_1 \oplus \underbrace{R_2 \oplus \cdots \oplus \underbrace{R_{k+1} \oplus R^c_{k+1}}_{=R^c_k}}_{=R^c_1 \;\; (\vdots)} \quad (\text{whole } R_1 \oplus \cdots \oplus R^c_{k+1} = R^c_0)
\end{array}
\tag{5.7}
$$

Note also that with $\bar{S}_1 := G^1_0$ and $\bar{S}_2 := 2G^2_0 \bar{z}_1$, the relations

$$S_i := P_{R^c_{i-1}}\, \bar{S}_i \,|_{N_{i-1}} \in L[\, N_{i-1}, R^c_{i-1} \,] \qquad \text{and} \qquad S_i \in GL[\, N^c_i, R_i \,], \quad i = 1, \dots, k+1$$

are in general valid, where the bijectivity of $S_i$ between the subspaces $N^c_i$ and $R_i$ is depicted, as usual, by arrows within (5.7).

Now, for closing the recursion (5.6) of the linear operators, the following formulas concerning $\bar{a}^{2k+1}$, $\bar{A}^{2k+1}$ and $M^{2k+1}$, $\widehat{M}_{k+1}$ have to be established according to

$$E_{k+1,k+1} := I_B \in L[\, B, B \,]$$

$$E_{i,k+1} := -S_i^{-1} P_{R_i} \sum_{\nu=i+1}^{k+1} \bar{S}_\nu E_{\nu,k+1} \in L[\, B, B \,], \quad i = k, \dots, 1$$

$$E^{k+1} := \begin{pmatrix} E_{1,1} & \cdots & E_{1,k+1} \\ & \ddots & \vdots \\ & & E_{k+1,k+1} \end{pmatrix} \in L[\, B^{k+1}, B^{k+1} \,] \tag{5.8}$$

$$\begin{pmatrix} a^{2k} & \bar{a}^{2k} \\ A^{2k} & \bar{A}^{2k} \end{pmatrix} := \left(C^{2k}\right)^{-1} \cdot E^k \cdot C^{2k} \in L[\, B^k, B^k \,]$$

$$\begin{pmatrix} a^{2k+1} & \bar{a}^{2k+1} \\ A^{2k+1} & \bar{A}^{2k+1} \end{pmatrix} := \left(C^{2k+1}\right)^{-1} \cdot E^{k+1} \cdot C^{2k+1} \in L[\, B^{k+1}, B^{k+1} \,]$$

with

$$a^{2k} \in L[\, B^{k-1}, B \,], \quad \bar{a}^{2k} \in L[\, B, B \,], \quad A^{2k} \in L[\, B^{k-1}, B^{k-1} \,], \quad \bar{A}^{2k} \in L[\, B, B^{k-1} \,]$$

$$a^{2k+1} \in L[\, B^k, B \,], \quad \bar{a}^{2k+1} \in L[\, B, B \,], \quad A^{2k+1} \in L[\, B^k, B^k \,], \quad \bar{A}^{2k+1} \in L[\, B, B^k \,]$$

and

$$M^{2k+1} := \begin{pmatrix} & & (a^{2k} & \bar{a}^{2k}) \\ & a^{2k-1} \cdot & (A^{2k} & \bar{A}^{2k}) \\ M^{2k-1} \cdot & A^{2k-1} \cdot & (A^{2k} & \bar{A}^{2k}) \end{pmatrix}_{\overline{k+1}} \in \ L[\, B^k, B^k \,], \tag{5.9}$$

where the index $\overline{k+1}$ denotes cancelling of last row of the $(k+1) \times k$ matrix in brackets. Using (5.8), (5.9) implies continuation of the linear operators in (5.6) without the $W$-Operators that are defined below in a non-iterative way. Next, setting

$$\boxed{\widehat{M}_{k+1}} := \left( \begin{array}{c|c} I_B & \\ & M^{2k+1} \end{array} \right), \tag{5.10}$$

the cone $C_k$ is given by

$$Z_k(\varepsilon, n_1^c, \dots, n_{k+1}^c, p_{k+1}) := z_0(\varepsilon) + \left[ \frac{\varepsilon^{2k+1}}{(2k+1)!} \cdots \frac{\varepsilon^{k+1}}{(k+1)!} \right] \cdot \boxed{\widehat{M}_{k+1}} \cdot \begin{pmatrix} n_1^c \\ \vdots \\ n_{k+1}^c + p_{k+1} \end{pmatrix}. \tag{5.11}$$

Finally, a simple calculation shows the upper triangularity

$$M_{i,j}^{2k+1} = 0\,, \quad i = 2, \dots, k\,, \quad j = 1, \dots, i-1 \qquad and \qquad M_{i,i}^{2k+1} = I_B\,, \quad i = 1, \dots, k \tag{5.12}$$

implying

$$M^{2k+1} \in GL[\, B^k, B^k \,] \tag{5.13}$$

with dependency of $M^{2k+1}$ from $z_1, \dots, z_{k-1}$ and $G_0^1, \dots, G_0^k$. In addition, we note that the iteratively defined mappings $E_{i,k+1}$ in (5.8) can also be written in an explicit way according to

$$E_{k+1,k+1} = I_B\,, \qquad E_{k,k+1} = -S_k^{-1}\, P_{R_k}\, \bar{S}_{k+1}$$

$$E_{i,k+1} := -S_i^{-1}\, P_{R_i} \left( I_B + \sum_{\nu=1}^{k-i} (-1)^\nu \sum_{i+1 \le n_1 < \cdots < n_\nu \le k} \prod_{\tau=1}^{\nu} \bar{S}_{n_\tau}\, S_{n_\tau}^{-1}\, P_{R_{n_\tau}} \right) \bar{S}_{k+1}\,, \quad i = k-1, \dots, 1.$$

It remains to state formulas concerning the $W$-Operators in (5.6). In contrast to the cone $C_k$ in (5.11), the formulas of the $W$-Operators are not defined recursively, but in an explicit way, exploiting the system of undetermined coefficients in combination with (2.3).

Finally, we should not forget, to supply explicit formulas concerning section 2. In some more detail $\Delta^k(z_{k-1}, \dots, z_1) \in L[B^k, \bar{B}^k]$ and $\mathrm{I}^k(z_k, \dots, z_1) \in \bar{B}^k$ from (2.4) as well as $W^{2k+1}(z_k, \dots, z_1)$ and $R^{2k+1}(z_k, \dots, z_1)$ from (2.12) are needed. Note also, that the formulas concerning $T^1[z_1]$, $T^2[z_2, z_1]$ and $T^3[z_3, z_2, z_1]$ are already stated in (2.2) and it is sufficient to supply formulas with respect to $T^{2k}[z_{2k}, \dots, z_1]$ and $T^{2k+1}[z_{2k+1}, \dots, z_1]$ for $k \ge 2$.

First, even components of system (2.2) are given for $k \ge 2$ by

$$T^{2k}[z_{2k}, \dots, z_1] = [\, W_{2k}^{2k}, \dots, W_k^{2k} \,](z_{k-1}, \dots, z_1) \cdot \begin{pmatrix} z_{2k} \\ \vdots \\ z_k \end{pmatrix} + R^{2k}(z_{k-1}, \dots, z_1) + \frac{(2k)!}{2(k!)^2}\, G_0^2 z_k^2 \ \in \ \bar{B} \tag{5.14}$$

with

$$W_{\mu}^{2k}(z_{k-1},\dots,z_1) := \frac{1}{\mu!}\sum_{\beta=1}^{2k} G_0^{\beta} \sum_{\substack{n_1+\cdots+n_k+1=\beta \\ 1\cdot n_1+\cdots+(k-1)\cdot n_{k-1}+\mu=2k}} \frac{(2k)!}{n_1!\cdots n_{k-1}!}\prod_{\tau=1}^{k-1}\left(\frac{1}{\tau!}\,z_\tau\right)^{n_\tau} \in L[B,\bar{B}] \tag{5.15}$$

for $\mu = 2k,\dots,k$, as well as

$$R^{2k}(z_{k-1},\dots,z_1) := \sum_{\beta=1}^{2k} G_0^{\beta} \sum_{\substack{n_1+\cdots+n_{k-1}=\beta \\ 1\cdot n_1+\cdots+(k-1)\cdot n_{k-1}=2k}} \frac{(2k)!}{n_1!\cdots n_{k-1}!}\prod_{\tau=1}^{k-1}\left(\frac{1}{\tau!}\,z_\tau\right)^{n_\tau} \in \bar{B} \tag{5.16}$$

closing (5.6) also with respect to the $W$-Operators.

Next, odd components of system (2.2) read for $k \geq 2$

$$T^{2k+1}[z_{2k+1},\dots,z_1] = \underbrace{\left[W_{2k+1}^{2k+1},\dots,W_{k+1}^{2k+1}\right](z_k,\dots,z_1)}_{=\,W^{2k+1}(z_k,\dots,z_1)}\cdot\begin{pmatrix} z_{2k+1} \\ \vdots \\ z_{k+1}\end{pmatrix} + R^{2k+1}(z_k,\dots,z_1) \in \bar{B} \tag{5.17}$$

with

$$W_{\mu}^{2k+1}(z_k,\dots,z_1) := \frac{1}{\mu!}\sum_{\beta=1}^{2k+1} G_0^{\beta} \sum_{\substack{n_1+\cdots+n_k+1=\beta \\ 1\cdot n_1+\cdots+k\cdot n_k+\mu=2k+1}} \frac{(2k+1)!}{n_1!\cdots n_k!}\prod_{\tau=1}^{k}\left(\frac{1}{\tau!}\,z_\tau\right)^{n_\tau} \in L[B,\bar{B}] \tag{5.18}$$

for $\mu = 2k+1,\dots,k+1$, whereas the corresponding inhomogeneity is given by

$$R^{2k+1}(z_k,\dots,z_1) := \sum_{\beta=1}^{2k+1} G_0^{\beta} \sum_{\substack{n_1+\cdots+n_k=\beta \\ 1\cdot n_1+\cdots+k\cdot n_k=2k+1}} \frac{(2k+1)!}{n_1!\cdots n_k!}\prod_{\tau=1}^{k}\left(\frac{1}{\tau!}\,z_\tau\right)^{n_\tau} \in \bar{B}\,. \tag{5.19}$$

Note that by (5.17)-(5.19), the remainder equation (2.12) from section 2 is well defined and it remains to state explicit formulas concerning the system of equations $T^{2k} = \cdots = T^{k+1} = 0$ from (2.4). Now, using $\Delta^1(z_1) = G_0^1 \in L[B,\bar{B}]$ and $\mathrm{I}^1(z_1) = G_0^2 z_1^2 \in \bar{B}$ in case of $k=1$ (cf. section 2), as well as combining even and odd components from (5.14) and (5.17), the upper triangular operator $\Delta^k(z_{k-1},\dots,z_1) \in L[B^k,\bar{B}^k]$ with associated inhomogeneity $\mathrm{I}^k(z_k,\dots,z_1) \in \bar{B}^k$ read for $k \geq 2$

$$\Delta_{i,i}^k(z_k,\dots,z_1) := G_0^1 \in L[B,\bar{B}]\,,\quad i=1,\dots,k \tag{5.20}$$

$$\Delta_{i,j}^k(z_k,\dots,z_1) := 0 \in L[B,\bar{B}]\,,\quad i=2,\dots,k,\quad j=1,\dots,i-1$$

$$\Delta_{k-i+1,k-j+1}^k(z_k,\dots,z_1) := \frac{1}{(k+j)!}\sum_{\beta=1}^{k+i} G_0^{\beta} \sum_{\substack{n_1+\cdots+n_k+1=\beta \\ 1\cdot n_1+\cdots+k\cdot n_k+(k+j)\cdot 1=k+i}} \frac{k!}{n_1!\cdots n_k!}\prod_{\tau=1}^{k}\left(\frac{1}{\tau!}\,z_\tau\right)^{n_\tau}$$

$$\in L[B,\bar{B}]\,,\quad i=2,\dots,k,\quad j=1,\dots,i-1,$$

and

$$I_{k-i+1}^k(z_k,\dots,z_1) := \sum_{\beta=1}^{k+i} G_0^{\beta} \sum_{\substack{n_1+\cdots+n_k=\beta \\ 1\cdot n_1+\cdots+k\cdot n_k=k+i}} \frac{k!}{n_1!\cdots n_k!}\prod_{\tau=1}^{k}\left(\frac{1}{\tau!}\,z_\tau\right)^{n_\tau} \in \bar{B}\,,\quad i=1,\dots,k\,.$$

For later use, we also note that by (5.15) and (5.18) we obtain for $k \geq 1$

$$W_{2k+1}^{2k+1}(z_k, \dots, z_0) = W_{2k}^{2k}(z_{k-1}, \dots, z_0) = G_0^1 \in L[B, \bar{B}] . \tag{5.21}$$

## *6. Proof of Main Results*

In the following Lemma, the main ingredients for proving Theorem 2 are comprised.

**Lemma 1 :** For $k \geq 1$ and arbitrary $[z_k, \dots, z_1] \in B^k$ we obtain

(i) $\quad N[\, \Delta^k(z_{k-1}, \dots, z_1) \,] = R[\, M^{2k+1}{}_{|\, N_1 \times \cdots \times N_k} \,]$

(ii) $\quad W^{2k+1}(z_k, \dots, z_1) \cdot \boxed{\widehat{M}_{k+1}} = [\, \bar{S}_1 \mid \bar{S}_2(z_1) \mid \cdots \mid \bar{S}_{k+1}(z_k, \dots, z_1) \,] \cdot C^{2k+1}$

(iii) $\quad N[\, [\, \bar{S}_1 \mid \bar{S}_2(z_1) \mid \cdots \mid \bar{S}_{k+1}(z_k, \dots, z_1) \,]_{|\, N_0 \times \cdots \times N_k} \,] = R[\, E^{k+1}{}_{|\, N_1 \times \cdots \times N_{k+1}} \,]$

(iv) $$M^{2k+1} \cdot \begin{pmatrix} d_{2k+1,2} & & \\ & \ddots & \\ & & d_{2k+1,k+1} \end{pmatrix} = \begin{pmatrix} d_{2k+1,2} & & \\ & \ddots & \\ & & d_{2k+1,k+1} \end{pmatrix} \cdot \begin{pmatrix} & a^{2k+1} \\ M^{2k+1} \cdot & A^{2k+1} \end{pmatrix}_{\overline{k+1}}$$

(v) $\quad$ Assume a curve $z_0(\varepsilon)$ as in (3.9) with $\bar{z}_l \neq 0,\ 1 \leq l \leq k$. Then,

$$\|\, G[\, z_0(\varepsilon) \,] \,\| = O(\, |\varepsilon|^{2k+1} \,) \implies \bar{z}_l \in N_{k+1}$$

**Proof of Theorem 2 (i) and (ii) :** We start with an arbitrary curve of the form

$$z_0(\varepsilon) = \varepsilon \cdot \bar{z}_1 + \cdots + \tfrac{1}{k!}\, \varepsilon^k \cdot \bar{z}_k + \tfrac{1}{(k+1)!}\, \varepsilon^{k+1} \cdot \bar{z}_{k+1} + \cdots + \tfrac{1}{(2k+1)!}\, \varepsilon^{2k+1} \cdot \bar{z}_{2k+1} + \varepsilon^{2k+2} \cdot r(\varepsilon)$$

with $\bar{z}_l \neq 0,\ 1 \leq l \leq k$ and $z_0(\varepsilon)$ of class $C^q,\ q \geq 2k+2$ and $k \geq 1$.

As already mentioned, the main idea of the ansatz is given by supplementing the coefficients $[\bar{z}_{k+1}, \dots, \bar{z}_{2k+1}]$ of $z_0(\varepsilon)$ with free coefficients $[b_{k+1}, \dots,\ b_{2k+1}]$ according to

$$G[\, z_0(\varepsilon) + \left[ \tfrac{\varepsilon^{2k+1}}{(2k+1)!} \cdots \tfrac{\varepsilon^{k+1}}{(k+1)!} \right] \cdot \begin{pmatrix} b_{2k+1} \\ \vdots \\ b_{k+1} \end{pmatrix} ]$$

$$= G[\, \varepsilon \cdot \bar{z}_1 + \cdots + \frac{1}{k!}\, \varepsilon^k \cdot \bar{z}_k$$

$$+ \frac{1}{(k+1)!}\, \varepsilon^{k+1} \cdot (\bar{z}_{k+1} + b_{k+1}) + \cdots + \frac{1}{(2k+1)!}\, \varepsilon^{2k+1} \cdot (\bar{z}_{2k+1} + b_{2k+1}) + \varepsilon^{2k+2} \cdot r(\varepsilon) \,]$$

$$\overset{(2.1)}{=} \sum_{i=1}^{k} \frac{1}{i!}\, \varepsilon^i \cdot T^i[\bar{z}_i, \dots, \bar{z}_1] + \sum_{i=k+1}^{2k} \frac{1}{i!}\, \varepsilon^i \cdot T^i[\bar{z}_i + b_i, \dots, \bar{z}_1]$$

$$+ \frac{1}{(2k+1)!}\, \varepsilon^{2k+1} \cdot T^{2k+1}[\bar{z}_{2k+1} + b_{2k+1}, \dots, \bar{z}_1] + \varepsilon^{2k+2} \cdot r_1(\varepsilon, b_{2k+1}, \dots, b_{k+1})$$

$$= \sum_{i=1}^{k} \frac{1}{i!}\, \varepsilon^i \cdot T^i[\bar{z}_i, \dots, \bar{z}_1] + \left[ \frac{\varepsilon^{2k}}{(2k)!} \cdots \frac{\varepsilon^{k+1}}{(k+1)!} \right] \cdot \begin{pmatrix} T^{2k}[\bar{z}_{2k} + b_{2k}, \dots, \bar{z}_1] \\ \vdots \\ T^{k+1}[\bar{z}_{k+1} + b_{k+1}, \dots, \bar{z}_1] \end{pmatrix}$$

$$+\ \frac{1}{(2k+1)!}\,\varepsilon^{2k+1}\cdot T^{2k+1}[\bar z_{2k+1}+b_{2k+1},\dots,\bar z_1]\ +\ \varepsilon^{2k+2}\cdot r_1(\varepsilon,b_{2k+1},\dots,b_{k+1})$$

$$\overset{(2.4)}{\underset{}{\overset{(2.12)}{=}}}\ \sum_{i=1}^{k}\frac{1}{i!}\,\varepsilon^i\cdot T^i[\bar z_i,\dots,\bar z_1]\ +\left[\frac{\varepsilon^{2k}}{(2k)!}\cdots\frac{\varepsilon^{k+1}}{(k+1)!}\right]\cdot\left\{\Delta^k\cdot\begin{pmatrix}\bar z_{2k}+b_{2k}\\ \vdots\\ \bar z_{k+1}+b_{k+1}\end{pmatrix}+\ \mathrm{I}^k(\bar z_k,\dots,\bar z_1)\right\}$$

$$+\ \frac{1}{(2k+1)!}\,\varepsilon^{2k+1}\cdot\left\{W^{2k+1}(\bar z_k,\dots,\bar z_1)\cdot\begin{pmatrix}\bar z_{2k+1}+b_{2k+1}\\ \vdots\\ \bar z_{k+1}+b_{k+1}\end{pmatrix}+\ R^{2k+1}(\bar z_k,\dots,\bar z_1)\right\}$$

$$+\,\varepsilon^{2k+2}\cdot r_1(\varepsilon,b_{2k+1},\dots,b_{k+1})=G[\cdot]\,. \tag{6.1}$$

Next, restrict $[b_{2k+1},\dots,b_{k+1}]\in B^{k+1}$ to $B\times N[\Delta^k]$ and perform an expansion of $r_1(\varepsilon,\cdot)$ with respect to $[b_{2k+1},\dots,b_{k+1}]$ to obtain

$$G[\cdot]=\sum_{i=1}^{k}\frac{1}{i!}\,\varepsilon^i\cdot T^i[\bar z_i,\dots,\bar z_1]\ +\left[\frac{\varepsilon^{2k}}{(2k+1)!}\cdots\frac{\varepsilon^{k+1}}{(k+1)!}\right]\cdot\begin{pmatrix}T^{2k}[\bar z_{2k},\dots,\bar z_1]\\ \vdots\\ T^{k+1}[\bar z_{k+1},\dots,\bar z_1]\end{pmatrix}$$

$$+\ \frac{1}{(2k+1)!}\,\varepsilon^{2k+1}\cdot\left\{T^{2k+1}[\bar z_{2k+1},\dots,\bar z_1]\ +\ W^{2k+1}(\bar z_k,\dots,\bar z_1)\cdot\begin{pmatrix}b_{2k+1}\\ \vdots\\ b_{k+1}\end{pmatrix}\right\} \tag{6.2}$$

$$+\,\varepsilon^{2k+2}\cdot\{\,r_1(\varepsilon,0,\dots,0)\ +\ r_2(\varepsilon,b_{2k+1},\dots,b_{k+1})\,\}$$

with smooth remainder operators $r_1(\cdot)$ and $r_2(\cdot)$ with $r_2(\varepsilon,0,\dots,0)=0$. Now, with $b_{2k+1}=\cdots=b_{k+1}=0$, formulas (6.1), (6.2) obviously simplify to $G[\cdot]=G[z_0(\varepsilon)]$ and we end up with

$$G[\cdot]=G[z_0(\varepsilon)]\ +\ \frac{1}{(2k+1)!}\,\varepsilon^{2k+1}\cdot W^{2k+1}(\bar z_k,\dots,\bar z_1)\cdot\begin{pmatrix}b_{2k+1}\\ \vdots\\ b_{k+1}\end{pmatrix}\ +\ \varepsilon^{2k+2}\cdot r_2(\varepsilon,b_{2k+1},\dots,b_{k+1})$$

$$\underset{(i)}{\overset{(5.10)}{=}}\ G[z_0(\varepsilon)]\ +\ \frac{1}{(2k+1)!}\,\varepsilon^{2k+1}\cdot W^{2k+1}(\bar z_k,\dots,\bar z_1)\cdot\boxed{\widehat M_{k+1}}\cdot\begin{pmatrix}n_0\\ \vdots\\ n_k\end{pmatrix}\ +\ \varepsilon^{2k+2}\cdot r_3(\varepsilon,n_0,\dots,n_k)$$

$$\underset{(ii)}{=}\ G[z_0(\varepsilon)]\ +\ \varepsilon^{2k+1}\cdot\overbrace{\frac{1}{(2k+1)!}\cdot[\bar S_1\,|\cdots|\,\bar S_{k+1}(\bar z_k,\dots,\bar z_l)]\cdot C^{2k+1}}^{=:\,\hat L_{k+1}}\cdot\begin{pmatrix}n_0\\ \vdots\\ n_k\end{pmatrix}+\varepsilon^{2k+2}\cdot r_3(\varepsilon,n_0,\dots,n_k)$$

$$=G[\,z_0(\varepsilon)\,]\ +\ \varepsilon^{2k+1}\cdot\hat L_{k+1}\cdot\begin{pmatrix}n_0\\ \vdots\\ n_k\end{pmatrix}\ +\ \varepsilon^{2k+2}\cdot r_3(\varepsilon,n_0,\dots,n_k) \tag{6.3}$$

by invoking Lemma 1 (i), (ii), the definition (5.10) of $\widehat M_{k+1}$ and $[n_0,\dots,n_k]\in N_0\times\cdots\times N_k\subset B^{k+1}$.

Then, using fine resolution (5.7), the free variables $[n_0,\dots,n_k]$ in (6.3) can further be restricted to $U_1^c\times\cdots\times U_k^c\times\left(U_{k+1}^c\oplus U_{k+1}^p\right)\subset N_1^c\times\cdots\times N_k^c\times(N_{k+1}^c\oplus P_{k+1})$, hence mapping the cone $C_k$ from $B$ to $\bar B$ according to

$$G[\,Z_k(\varepsilon, n_1^c, \dots, n_{k+1}^c, p_{k+1})\,] \overset{(5.11)}{=} G[\,z_0(\varepsilon) + \underbrace{\left[\tfrac{\varepsilon^{2k+1}}{(2k+1)!} \cdots \tfrac{\varepsilon^{k+1}}{(k+1)!}\right] \cdot \boxed{\widehat{M}_{k+1}}}_{= A_\varepsilon = [\,A_\varepsilon^1 \cdots A_\varepsilon^{k+1}\,]} \cdot \begin{pmatrix} n_1^c \\ \vdots \\ n_{k+1}^c + p_{k+1} \end{pmatrix}]$$

$$\overset{(6.3)}{=} G[\,z_0(\varepsilon)\,] + \varepsilon^{2k+1} \cdot \hat{L}_{k+1} \cdot \begin{pmatrix} n_1^c \\ \vdots \\ n_{k+1}^c + p_{k+1} \end{pmatrix} + \varepsilon^{2k+2} \cdot r_3(\varepsilon, n_1^c, \dots, n_{k+1}^c + p_{k+1}) \tag{6.4}$$

$$= \sum_{i=l}^{2k} \frac{1}{i!}\, \varepsilon^i \cdot T^i[\bar{z}_i, \dots, \bar{z}_l]$$

$$+ \varepsilon^{2k+1} \cdot \{ \frac{1}{(2k+1)!} T^{2k+1}[\bar{z}_{2k+1}, \dots, \bar{z}_l] + \hat{L}_{k+1} \cdot \begin{pmatrix} n_1^c \\ \vdots \\ n_{k+1}^c + p_{k+1} \end{pmatrix} + \varepsilon \cdot r_4(\varepsilon, n_1^c, \dots, n_{k+1}^c, p_{k+1}) \}$$

$$=: \sum_{i=l}^{2k} \frac{1}{i!}\, \varepsilon^i \cdot T^i[\bar{z}_i, \dots, \bar{z}_l] + \varepsilon^{2k+1} \cdot H(\varepsilon, n_1^c, \dots, n_{k+1}^c, p_{k+1})$$

and by direct inspection of (5.6), (5.7) we obtain

$$\hat{L}_{k+1} \in GL[\,N^c,\, R_1 \oplus R_2 \oplus \cdots \oplus R_{k+1}\,]. \tag{6.5}$$

Now, (6.4) delivers (3.8) as well as Theorem 2 (i).

Concerning Theorem 2 (ii), assume $z_0(\varepsilon)$ now to be an approximation of order $2k$. Then, we obtain from (6.4) the blown-up remainder equation

$$H(\varepsilon, n_1^c, \dots, n_{k+1}^c, p_{k+1}) = \tag{6.6}$$

$$\frac{1}{(2k+1)!} T^{2k+1}[\bar{z}_{2k+1}, \dots, \bar{z}_l] + \hat{L}_{k+1} \cdot \begin{pmatrix} n_1^c \\ \vdots \\ n_{k+1}^c + p_{k+1} \end{pmatrix} + \varepsilon \cdot r_4(\varepsilon, n_1^c, \dots, n_{k+1}^c, p_{k+1}) = 0$$

with $H(\varepsilon, n_1^c, \dots, n_{k+1}^c, p_{k+1}) \in C^{q-2k-1}(\,\mathcal{B}_{\delta_1}(0) \times U_1^c \times \cdots \times U_{k+1}^c \times U_{k+1}^p, \bar{B}\,)$ and $P_{k+1}$ chosen as a direct complement of $\bar{z}_l$, due to Lemma 1 (v) and $N_{k+1} = P_{k+1} \oplus \{\bar{z}_l\}$.

Finally, suppose $Z_k$ to define a *k-transversal cone*, then $\hat{L}_{k+1} \in GL[N^c, \bar{B}]$ by (6.5), yielding for $\varepsilon = 0$ the base solutions

$$H(0, n_1^c, \dots, n_{k+1}^c, p_{k+1}) = 0 \quad \Leftrightarrow$$

$$\begin{pmatrix} n_1^c \\ \vdots \\ n_{k+1}^c \end{pmatrix}(0, p_{k+1}) \overset{(6.3)}{:=} -\hat{L}_{k+1}^{-1} \cdot \frac{1}{(2k+1)!}\left[\,T^{2k+1}[\bar{z}_{2k+1}, \dots, \bar{z}_l] + c_{2k+1,k+1} \cdot \bar{S}_{k+1}(\bar{z}_k, \dots, \bar{z}_l) \cdot p_{k+1}\,\right]$$

and the implicit function theorem delivers locally unique and smooth continuation of these solutions to $\varepsilon \neq 0$, thus finishing the proof of Theorem 2 (ii), possibly after reduction of $\delta_1 > 0$.

**Proof of Theorem 2 (iii) :** The approximation condition of order $2k$ concerns the equations $T^1[\bar{z}_1] = \cdots = T^{2k}[\bar{z}_{2k}, \dots, \bar{z}_1] = 0$, merely depending from $G_0^1, \dots, G_0^{2k}$ according to (5.14) and (5.17). Hence, perturbations of $G[z]$ of order $O(\|z\|^{2k+1})$ cannot destroy the approximation property of $z_0(\varepsilon)$.

Further, by tedious, but simple inspection of the iteration from section 5, we see that the construction of a $k$-transversal cone $C_k$ only is affected by $G_0^1, \dots, G_0^{k+1}$, completing the proof of (iii).

**Proof of Theorem 2 (iv) :** Setting $p := p_{k+1}$, define in view of (6.4)

$$\bar{G}_{\varepsilon,p}[\, n^c \,] := G[\, Z_k(\varepsilon, n^c, p) \,] = G[\, z_0(\varepsilon) + A_\varepsilon \cdot n^c + A_\varepsilon^{k+1} \cdot p \,]\,,$$

implying $\bar{G}_{\varepsilon,p} : N^c \longrightarrow \bar{B}$ to be a local diffeomorphism for $\varepsilon \neq 0$, due to supposed transversality of the cone $C_k$, i.e. $\hat{L}_{k+1} \in GL[N^c, \bar{B}]$. Now, with $\bar{r}_3(\varepsilon, n_1^c, \dots, n_{k+1}^c, p) := r_3(\varepsilon, n_1^c, \dots, n_{k+1}^c + p)$, as well as defining subsequently

$$G_{\varepsilon,p}[\, z \,] := G[\, z_0(\varepsilon) + z + A_\varepsilon^{k+1} \cdot p \,]$$

$$\mathcal{J}_{\varepsilon,p}[\, \bar{b} \,] := \bar{b} - G[\, z_0(\varepsilon) \,] - \varepsilon^{2k+1} \cdot \tfrac{c_{2k+1,k+1}}{(2k+1)!} \cdot \bar{S}_{k+1} \cdot p - \varepsilon^{2k+2} \cdot \bar{r}_3\big(\, \varepsilon,\ \bar{G}_{\varepsilon,p}^{-1}[\, \bar{b} \,],\ p \,\big) \tag{6.7}$$

$$=:\ \bar{b} + \varepsilon \cdot r_5\big(\varepsilon, p, \bar{b}\big)\,,$$

we obtain by direct calculation from (6.4)

$$\mathcal{J}_{\varepsilon,p} \circ G_{\varepsilon,p} \circ A_\varepsilon \cdot n^c = \varepsilon^{2k+1} \cdot \hat{L}_{k+1} \cdot n^c\,,$$

ending the proof of Theorem 2 (iv).

**Proof of Theorem 2 (v) :** First note that the operator $A_\varepsilon$ from (6.4) is defined by the matrix operator $\widehat{M}_{k+1}$ from (5.8), (5.10) and (5.12), showing that $A_\varepsilon$ can be applied to every element from $B^{k+1}$. Hence, using the upper triangularity of $\widehat{M}_{k+1}$ according to (5.10) and (5.12), we obtain

$$A_\varepsilon \cdot \begin{pmatrix} n_0 \\ \vdots \\ n_k \end{pmatrix} = \left[ \tfrac{\varepsilon^{2k+1}}{(2k+1)!} \cdot I_B \,\middle|\, \tfrac{\varepsilon^{2k}}{(2k)!} \cdot I_B \,\middle|\, \tfrac{\varepsilon^{2k-1}}{(2k-1)!} \cdot [\, I_B + O(\,|\varepsilon|\,)\,] \,\middle|\, \cdots \,\middle|\, \tfrac{\varepsilon^{k+1}}{(k+1)!} \cdot [\, I_B + O(\,|\varepsilon|\,)\,] \right] \cdot \begin{pmatrix} n_0 \\ \vdots \\ n_k \end{pmatrix}$$

$$= A_\varepsilon^1 \cdot n_0 + \cdots + A_\varepsilon^{k+1} \cdot n_k \tag{6.8}$$

with $A_\varepsilon^i \in L[N_{i-1}, B]$. In addition, the existence of constants $0 < c_i < d_i, i = 1, \dots, k+1$ with

$$\begin{array}{ccccc} c_1 \cdot |\varepsilon|^{2k+1} \cdot \|\, n_0 \,\| & \leq & \|\, A_\varepsilon^1 \cdot n_0 \,\| & \leq & d_1 \cdot |\varepsilon|^{2k+1} \cdot \|\, n_0 \,\| \\ \cdots & \cdots & \cdots & \cdots & \cdots \\ c_{k+1} \cdot |\varepsilon|^{k+1} \cdot \|\, n_k \,\| & \leq & \left\|\, A_\varepsilon^{k+1} \cdot n_k \,\right\| & \leq & d_{k+1} \cdot |\varepsilon|^{k+1} \cdot \|\, n_k \,\| \end{array} \tag{6.9}$$

can obviously be assured. Then, by (6.1) and (6.3)

$$G[\, z_0(\varepsilon) \ + \underbrace{\left[ \tfrac{\varepsilon^{2k+1}}{(2k+1)!} \ \cdots \ \tfrac{\varepsilon^{k+1}}{(k+1)!} \right] \cdot \boxed{\widehat{M}_{k+1}}}_{=\, A_\varepsilon \,=\, [\, A_\varepsilon^1 \ \cdots \ A_\varepsilon^{k+1} \,]} \cdot \begin{pmatrix} n_0 \\ \vdots \\ n_k \end{pmatrix} ]$$

$$= G[\, z_0(\varepsilon) \,] \ + \ \varepsilon^{2k+1} \cdot \hat{L}_{k+1} \cdot \begin{pmatrix} n_0 \\ \vdots \\ n_k \end{pmatrix} \ + \ \varepsilon^{2k+2} \cdot r_3(\varepsilon, n_0, \dots, n_k) \tag{6.10}$$

and differentiation with respect to $N_0 \times \cdots \times N_k$ at $n_0 = \cdots = n_k = 0$ implies by chain rule

$$G'[\,z_0(\varepsilon)\,]\cdot A_\varepsilon \cdot \begin{pmatrix} n_0 \\ \vdots \\ n_k \end{pmatrix} = \varepsilon^{2k+1}\cdot \hat{L}_{k+1}\cdot \begin{pmatrix} n_0 \\ \vdots \\ n_k \end{pmatrix} + \varepsilon^{2k+2}\cdot \underbrace{r_{3,n_0,\dots,n_k}(\varepsilon,0,\dots,0)}_{=:\,H_\varepsilon=[\,H_\varepsilon^1\ \cdots\ H_\varepsilon^{k+1}\,]}\cdot \begin{pmatrix} n_0 \\ \vdots \\ n_k \end{pmatrix} \qquad (6.11)$$

with a smooth family of bounded operators $H_\varepsilon$. Then, using the definition of $\hat{L}_{k+1}$ in (6.3), the equation can further be splitted with norms of partial derivatives satisfying for $i = 1, \dots, k+1$

$$\frac{\|\, G'[\,z_0(\varepsilon)\,]\cdot A_\varepsilon^i \cdot n_{i-1} \,\|}{\|\, A_\varepsilon^i \cdot n_{i-1} \,\|} = \frac{\|\, \varepsilon^{2k+1}\cdot \frac{c_{2k+1,i}}{(2k+1)!}\cdot \bar{S}_i \cdot n_{i-1} + \varepsilon^{2k+2}\cdot H_\varepsilon^i \cdot n_{i-1} \,\|}{\|\, A_\varepsilon^i \cdot n_{i-1} \,\|}$$

$$\overset{(6.9)}{\leq} \frac{|\varepsilon|^{2k+1}\cdot \|\, \frac{c_{2k+1,i}}{(2k+1)!}\cdot \bar{S}_i + \varepsilon\cdot H_\varepsilon^i \,\| \cdot \cancel{\|\, n_{i-1} \,\|}}{c_i \cdot |\varepsilon|^{2k+2-i}\cdot \cancel{\|\, n_{i-1} \,\|}} \leq \underbrace{e_i}_{>0}\cdot |\varepsilon|^{i-1} = O(\,|\varepsilon|^{i-1}\,)$$

ending up with

$$\|\, G'[\,z_0(\varepsilon)\,]_{|\,R[\,A_\varepsilon^1|_{N_0}\,]} \,\| = O(\,|\varepsilon|^0\,)\,,\ \ \cdots\ \ ,\|\, G'[\,z_0(\varepsilon)\,]_{|\,R\left[\,A_\varepsilon^{k+1}|_{N_k}\,\right]} \,\| = O\big(\,|\varepsilon|^k\big)\,.$$

Now, due to $N_1^c \subset N_0\,, \dots, N_{k+1}^c \subset N_k$, the first part of Theorem 2 (v) is shown.

In the next step, let us estimate the operator norm of the inverse of the Jacobian with respect to the complement $N^c$ along the approximation $z_0(\varepsilon)$.

First, from the settings (5.7), (5.9), it is not too difficult to see that for $\varepsilon \neq 0$ we obtain $A_\varepsilon \in GL[N^c, N^c]$ and even in some more detail $A_\varepsilon^1 \in L[N_1^c, N_1^c]$, $A_\varepsilon^2 \in L[N_2^c, N_2^c]$ as well as $A_\varepsilon^i \in L[N_i^c, N_1^c \oplus \cdots \oplus N_{i-2}^c \oplus N_i^c]$ in case of $3 \leq i \leq k+1$.

Now, differentiation of (6.4) with respect to $N^c$ implies at every point of the cone

$$G'[\,Z_k(\varepsilon, n^c, p_{k+1})\,]\cdot A_\varepsilon \cdot \bar{n}^c = \varepsilon^{2k+1}\cdot [\,\hat{L}_{k+1} + \varepsilon\cdot H_\varepsilon\,]\cdot \bar{n}^c, \quad \bar{n}^c \in N^c \qquad (6.12)$$

and by $\hat{L}_{k+1} \in GL[N^c, \bar{B}]$ one obtains $G'[\,Z_k(\varepsilon, n^c, p_{k+1})\,] \in GL[N^c, \bar{B}]$ in case of $\varepsilon \neq 0$. Further, using (6.12), we obtain for arbitrary $\bar{n}^c \in N^c$ the lower bound

$$\frac{\|\, G'[\,Z_k(\varepsilon, n^c, p_{k+1})\,]\cdot A_\varepsilon \cdot \bar{n}^c \|}{\|\, A_\varepsilon \cdot \bar{n}^c \|} = |\varepsilon|^{2k+1}\cdot \frac{\|\, \hat{L}_{k+1}\cdot n^c + \varepsilon\cdot H_\varepsilon \cdot \bar{n}^c \|}{\|\, A_\varepsilon \cdot \bar{n}^c \|}$$

$$\geq |\varepsilon|^{2k+1}\cdot \frac{\left|\, \overbrace{\|\, \hat{L}_{k+1}\cdot \bar{n}^c \|}^{\geq d\cdot \|\, \bar{n}^c \,\|} - |\varepsilon| \cdot \overbrace{\|\, H_\varepsilon \cdot \bar{n}^c \|}^{\leq e\cdot \|\, \bar{n}^c \,\|} \,\right|}{\underbrace{\|\, A_\varepsilon \cdot \bar{n}^c \|}_{\leq f\cdot |\varepsilon|^{k+1}\cdot \|\, \bar{n}^c \,\|}}$$

$$\geq |\varepsilon|^{2k+1}\cdot \frac{(\,d - |\varepsilon|\cdot e\,)\cdot \cancel{\|\, \bar{n}^c \|}}{f\cdot |\varepsilon|^{k+1}\cdot \cancel{\|\, \bar{n}^c \|}} \geq \underbrace{g}_{>0}\cdot |\varepsilon|^k$$

with constants $d > 0$ by $\hat{L}_{k+1} \in GL[N^c, \bar{B}]$ and $f > 0$ by (6.8), (6.9) and $N_1^c \subset N_0\,, \dots, N_{k+1}^c \subset N_k$. Note that the constants are uniformly valid with respect to the complete cone.

Now, $R[A_\varepsilon|_{N^c}] = N^c$ in case of $\varepsilon \neq 0$ and we obtain

$$\|\, G'[\,Z_k(\varepsilon, n^c, p_{k+1})\,]\cdot \bar{n}^c \| \geq g\cdot |\varepsilon|^k \cdot \|\, \bar{n}^c \|\,.$$

Thus, the linearization with respect to the complement $N^c$, i.e.

$$G_{N^c}[\,Z_k(\varepsilon, n^c, p_{k+1})\,] := G'[\,Z_k(\varepsilon, n^c, p_{k+1})\,]_{|N^c}$$

satisfies

$$\|\,G_{N^c}[\,Z_k(\varepsilon, n^c, p_{k+1})\,]^{-1}\,\| \;=\; \sup_{\|\,\bar{b}\,\|=1} \|\,G_{N^c}[\,Z_k(\varepsilon, n^c, p_{k+1})\,]^{-1}\cdot\bar{b}\,\|$$

$$= \sup_{\|\,\bar{n}^c\,\|=1} \|\,G_{N^c}[\,\cdot\,]^{-1}\cdot\frac{G_{N^c}[\,\cdot\,]\cdot\bar{n}^c}{\|\,G_{N^c}[\,\cdot\,]\cdot\bar{n}^c\|}\,\| \;=\; \sup_{\|\,n^c\,\|=1}\frac{\|\,\bar{n}^c\|}{\|\,G_{N^c}[\,\cdot\,]\cdot\bar{n}^c\|} \;\le\; \frac{\cancel{\|\,\bar{n}^c\|}}{g\cdot|\varepsilon|^k\cdot\cancel{\|\,\bar{n}^c\|}}\,,$$

yielding $\|G_{N^c}[\,Z_k(\varepsilon, n^c, p_{k+1})\,]^{-1}\| = O\big(|\varepsilon|^{-k}\big)$ at every point of a *$k$-transversal cone*. In particular, we obtain along the center line $\|G_{N^c}[\,z_0(\varepsilon)\,]^{-1}\| = O\big(|\varepsilon|^{-k}\big)$, as desired.

**Proof of Theorem 2 (vi) :** Finally, assume $N^c$ to be of finite dimension. Then the linear mappings $A_\varepsilon$ from (6.4) and $\hat{L}_{k+1}$ from (6.3) can be represented by matrices according to

$$\boxed{A_\varepsilon} = \begin{pmatrix} \boxed{\frac{\varepsilon^{2k+1}}{(2k+1)!}\cdot I_{d_1}} & 0 & O\big(|\varepsilon|^{2k}\big) & \cdots & O\big(|\varepsilon|^{k+2}\big) \\ & \boxed{\frac{\varepsilon^{2k}}{(2k)!}\cdot I_{d_2}} & 0 & \ddots & \vdots \\ & & \boxed{\frac{\varepsilon^{2k-1}}{(2k-1)!}\cdot I_{d_3}} & \ddots & O\big(|\varepsilon|^{k+2}\big) \\ & & & \ddots & 0 \\ & & & & \boxed{\frac{\varepsilon^{k+1}}{(k+1)!}\cdot I_{d_{k+1}}} \end{pmatrix} \in \mathbb{K}^{d,d}$$

and

$$\boxed{\hat{L}_{k+1}} = \frac{1}{(2k+1)!}\cdot\begin{pmatrix} \boxed{I_{d_1}} & O(1) & \cdots & O(1) \\ & \boxed{I_{d_2}} & \ddots & \vdots \\ & & \ddots & O(1) \\ & & & \boxed{I_{d_{k+1}}} \end{pmatrix}\cdot C^{2k+1} \in \mathbb{K}^{d,d}.$$

Here, $d_i := dim(N_i^c) \ge 0$, $d := d_1 + \cdots + d_{k+1} \ge 1$ and canonical bases chosen in $R_i$ according to $R_i = R[S_i]$, $S_i \in GL[N_i^c, R_i]$, $i = 1, \ldots, k+1$. Obviously,

$$det\{\boxed{A_\varepsilon}\} = c_1\cdot\varepsilon^{(2k+1)\cdot d_1 \,+\, (2k)\cdot d_2 + \cdots + \,(k+1)\cdot d_{k+1}},\;\; c_1 \neq 0$$

$$det\{\,\varepsilon^{2k+1}\cdot[\,\boxed{\hat{L}_{k+1}}\; + \; O(|\varepsilon|)\,]\,\} = c_2\cdot\varepsilon^{(2k+1)\cdot d}\cdot[\,1 + O(|\varepsilon|)\,]\,,\; c_2 \neq 0$$

and by determinant multiplication theorem, equation (6.12) implies

$$det\{\,G'[\,Z_k(\varepsilon, n^c, p_{k+1})\,]_{\,|\,N^c}\,\}$$

$$= c_3\cdot\varepsilon^{(2k+1)\cdot(d_1 \,+\, \cdots \,+\, d_{k+1}) \;-\; (2k+1)\cdot d_1 \,-\, (2k)\cdot d_2 \,-\, \cdots \,-\, (k+1)\cdot d_{k+1}}\cdot[\,1 + O(|\varepsilon|)\,]$$

$$= c_3\cdot\varepsilon^{1\cdot d_2 \,+\, \cdots \,+\, k\cdot d_{k+1}}\cdot[\,1 + O(|\varepsilon|)\,] \;=\; c_3\cdot\varepsilon^{\chi}\cdot[\,1 + O(|\varepsilon|)\,]\,,\;\; c_3 \neq 0$$

with $\chi = 1 \cdot dim(N_2^c) + \cdots + k \cdot dim(N_{k+1}^c)$, finishing the proof of Theorem 2.

**Proof of Corollary 1 (i)→(ii) :** First, we choose $U := N^c$, $V := N_{k+1}$ implying (3.15) by Theorem 2 (v) and it remains to show (3.16). Now, for calculation of the level sets in the cone $C_k$ with respect to $z_0(\varepsilon)$, we use identity (6.10) with $(n_0, \dots, n_{k-1}, n_k)$ replaced by $(n_1^c, \dots, n_k^c, n_{k+1}^c + n_{k+1})$, yielding the level set requirement

$$G[\, z_0(\varepsilon) + A_\varepsilon \cdot \begin{pmatrix} n_1^c \\ \vdots \\ n_{k+1}^c + n_{k+1} \end{pmatrix} ]$$

$$\overset{(6.10)}{=} G[\, z_0(\varepsilon)\, ] + \varepsilon^{2k+1} \cdot \hat{L}_{k+1} \cdot \begin{pmatrix} n_1^c \\ \vdots \\ n_{k+1}^c + n_{k+1} \end{pmatrix} + \varepsilon^{2k+2} \cdot r_3(\varepsilon, n_1^c, \dots, n_{k+1}^c + n_{k+1}) \overset{!}{=} G[\, z_0(\varepsilon)\, ]$$

$$\underset{(6.3)}{\overset{\varepsilon \neq 0}{\Leftrightarrow}} \quad \hat{L}_{k+1} \cdot \begin{pmatrix} n_1^c \\ \vdots \\ n_{k+1}^c \end{pmatrix} + \frac{c_{2k+1,k+1}}{(2k+1)!} \cdot \bar{S}_{k+1} \cdot n_{k+1} + \varepsilon \cdot r_3(\varepsilon, n_1^c, \dots, n_{k+1}^c + n_{k+1})$$

$$=: H(\varepsilon, n_1^c, \dots, n_{k+1}^c, n_{k+1}) = 0 \,.$$

Obviously, we have solutions $H(\varepsilon, 0, \dots, 0, 0) = 0$ and we can solve at $\varepsilon = 0$ according to

$$H(0, n_1^c, \dots, n_{k+1}^c, n_{k+1}) = 0 \quad \Leftrightarrow \quad \begin{pmatrix} \bar{n}_1^c \\ \vdots \\ \bar{n}_{k+1}^c \end{pmatrix} (0, n_{k+1}) := -\hat{L}_{k+1}^{-1} \cdot \frac{c_{2k+1,k+1}}{(2k+1)!} \cdot \bar{S}_{k+1} \cdot n_{k+1}$$

with $\bar{n}_{k+1}^c(0, n_{k+1}) = 0$. Hence, by $H_{n^c}(0, n_1^c, \dots, n_{k+1}^c, n_{k+1}) = \hat{L}_{k+1} \in GL[N^c, \bar{B}]$, we get smooth level sets satisfying

$$G[\, z_0(\varepsilon) + A_\varepsilon \cdot \begin{pmatrix} \bar{n}_1^c(\varepsilon, n_{k+1}) \\ \vdots \\ \bar{n}_{k+1}^c(\varepsilon, n_{k+1}) + n_{k+1} \end{pmatrix} ] = G[\, z_0(\varepsilon)\, ]$$

with

$$(\bar{n}_1^c, \dots, \bar{n}_{k+1}^c)(\varepsilon, 0) = 0 \quad and \quad \bar{n}_{k+1}^c(0, n_{k+1}) = 0 \,.$$

Now, differentiation in the direction of the level sets at $n_{k+1} = 0$ implies

$$G'[\, z_0(\varepsilon)\, ] \cdot A_\varepsilon \cdot \begin{pmatrix} \bar{n}_{1,n_{k+1}}^c(\varepsilon, 0) \\ \vdots \\ \bar{n}_{k+1,n_{k+1}}^c(\varepsilon, 0) + I_B \end{pmatrix}$$

$$\overset{(6.8)}{=} G'[\, z_0(\varepsilon)\, ] \cdot \frac{\varepsilon^{k+1}}{(k+1)!} \cdot [\, \bar{n}_{k+1,n_{k+1}}^c(\varepsilon, 0) + I_B + \varepsilon \cdot Q_1(\varepsilon)\, ]$$

$$= G'[\, z_0(\varepsilon)\, ] \cdot \frac{\varepsilon^{k+1}}{(k+1)!} \cdot [\, I_B + \varepsilon \cdot Q_2(\varepsilon)\, ] = 0$$

with smooth family of bounded linear operators $Q_2(\varepsilon) \in L[V, B]$ yielding (3.16) as desired.

**Proof of Corollary 1 (ii)→(i) :** The proof from (ii) back to (i) uses a contradiction argument. For this purpose, assume the cone $C_k$ from (5.11) *not* to be *k-transversal*, i.e. $R_{k+1}^c \neq \{0\}$. Then, considering decomposition (5.7) and the assumption $G_U[z_0(\varepsilon)] \in GL[U, \bar{B}], \varepsilon \neq 0$, we obtain the decompositions

$$\begin{array}{c} B \;=\; N^c \oplus N_{k+1} \\ \uparrow \\ \boxed{\varepsilon^{2k+1} \cdot \hat{L}_{k+1}} \\ \downarrow \\ \bar{B} \;=\; R \;\oplus\; \underbrace{R_{k+1}^c}_{\neq \{0\}} \end{array} \qquad \text{and} \qquad \begin{array}{c} B \;=\; U \;\oplus\; V \\ \uparrow \\ \boxed{G_U[\,z_0(\varepsilon)\,]} \\ \downarrow \\ \bar{B} \end{array} \tag{6.13}$$

with $N^c = N_1^c \oplus \cdots \oplus N_{k+1}^c$, $R = R_1 \oplus \cdots \oplus R_{k+1}$ and $\hat{L}_{k+1} \in GL[N^c, R]$ by construction. Now, $U$ is in bjection with $\bar{B}$ and $N^c$ is in bijection with $R \subsetneq \bar{B}$. Hence, the subspace $N_{k+1}$ cannot be contained in $V$, i.e. we have $N_{k+1} \not\subset V$.

Next, choose $\bar{n}_{k+1} \notin V$ and consider the partial derivative along $z_0(\varepsilon)$ in direction of

$$G'[\,z_0(\varepsilon)\,] \cdot A_\varepsilon \cdot \underbrace{\begin{pmatrix} \bar{a}^{2k+1} \\ \\ \bar{A}^{2k+1} \end{pmatrix}}_{=:\,\bar{m}} \cdot \bar{n}_{k+1} = G'[\,z_0(\varepsilon)\,] \cdot A_\varepsilon \cdot \begin{pmatrix} \bar{m}_1 \\ \vdots \\ \bar{m}_{k+1} \end{pmatrix} \cdot \bar{n}_{k+1} \tag{6.14}$$

with $\bar{m}_1 \in L[B, N_1^c], \dots, \bar{m}_k \in L[B, N_k^c]$ as well as $\bar{m}_{k+1} = I_B$ by the definitions in (5.8). Then, we will see that from (3.15), (3.17) and supposed nontransversality of the cone, the contradiction

$$\underbrace{\alpha_1}_{>0} \cdot |\varepsilon|^{2k+1} \leq \|\, G'[\,z_0(\varepsilon)\,] \cdot A_\varepsilon \cdot \bar{m} \cdot \bar{n}_{k+1} \,\| \leq \underbrace{\alpha_2}_{\geq 0} \cdot |\varepsilon|^{2k+2}, \tag{6.15}$$

will result, thus completing the proof.

First, for all $n_{k+1} \in N_{k+1}$, we obtain by (5.8) and (6.8)

$$\begin{aligned} A_\varepsilon \cdot \bar{m} \cdot n_{k+1} &= \left[ \frac{\varepsilon^{2k+1}}{(2k+1)!} \cdot I_B \middle| \cdots \middle| \frac{\varepsilon^{k+1}}{(k+1)!} \cdot [\, I_B + O(\,|\varepsilon|\,)\,] \right] \cdot \bar{m} \cdot n_{k+1} \\ &= \frac{\varepsilon^{k+1}}{(k+1)!} \cdot n_{k+1} + \varepsilon^{k+2} \cdot L_\varepsilon \cdot n_{k+1} \end{aligned} \tag{6.16}$$

with a smooth family of bounded linear operators $L_\varepsilon$ yielding

$$\begin{aligned} \|\, A_\varepsilon \cdot \bar{m} \cdot n_{k+1} \,\| &\geq |\varepsilon|^{k+1} \cdot \left| \frac{1}{(k+1)!} \cdot \|\, n_{k+1} \,\| - \|\, \varepsilon \cdot L_\varepsilon \cdot n_{k+1} \,\| \right| \\ &\geq |\varepsilon|^{k+1} \cdot \left[ \frac{1}{(k+1)!} - |\varepsilon| \cdot \|L_\varepsilon\| \right] \cdot \|\, n_{k+1} \,\| \;\geq\; \underbrace{d_1}_{>0} \cdot |\varepsilon|^{k+1} \cdot \|\, n_{k+1} \,\| \end{aligned} \tag{6.17}$$

for $d_1 > 0$ and $\varepsilon$ chosen sufficiently small. Then, considering (6.11) and (6.17), we see

$$\frac{\|\, G'[\,z_0(\varepsilon)\,] \cdot A_\varepsilon \cdot \overbrace{\begin{pmatrix} \bar{m}_1 \\ \vdots \\ \bar{m}_{k+1} \end{pmatrix}}^{\in N_1^c \times \cdots \times N_k^c \times N_{k+1}} \cdot n_{k+1} \,\|}{\|\, A_\varepsilon \cdot \bar{m} \cdot n_{k+1} \,\|} \overset{(6.11)}{=} \frac{\|\, \varepsilon^{2k+1} \cdot [\, \hat{L}_{k+1} \cdot \bar{m} \cdot n_{k+1} + \varepsilon \cdot H_\varepsilon \cdot \bar{m} \cdot n_{k+1} \,] \,\|}{\|\, A_\varepsilon \cdot \bar{m} \cdot n_{k+1} \,\|}$$

$$\overset{(6.17)}{\leq} \frac{\| \varepsilon^{2k+1} \cdot [\, \overbrace{\hat{L}_{k+1} \cdot \bar{m} \cdot n_{k+1}}^{=0} + \varepsilon \cdot H_\varepsilon \cdot \bar{m} \cdot n_{k+1} \,] \|}{d_1 \cdot |\varepsilon|^{k+1} \cdot \| n_{k+1} \|} \tag{6.18}$$

$$= \frac{|\varepsilon|^{2k+2} \cdot \| H_\varepsilon \cdot \bar{m} \cdot n_{k+1} \|}{d_1 \cdot |\varepsilon|^{k+1} \cdot \| n_{k+1} \|} \leq \frac{|\varepsilon|^{k+1} \cdot \| H_\varepsilon \| \cdot \| \bar{m} \| \cdot \cancel{\| n_{k+1} \|}}{d_1 \cdot \cancel{\| n_{k+1} \|}} \leq \underbrace{d_2}_{\geq 0} \cdot |\varepsilon|^{k+1} .$$

Note that $\hat{L}_{k+1} \cdot \bar{m} \cdot n_{k+1} = 0$ is a consequence of Lemma 1 (iii), implying by the definition of $\hat{L}_{k+1}$ in (6.3)

$$\hat{L}_{k+1} \cdot \begin{pmatrix} n_0 \\ \vdots \\ n_k \end{pmatrix} = 0$$

$$\Leftrightarrow \begin{pmatrix} n_0 \\ \vdots \\ n_k \end{pmatrix} = \left(C^{2k+1}\right)^{-1} \cdot E^{k+1} \cdot C^{2k+1} \cdot \begin{pmatrix} n_1 \\ \vdots \\ n_{k+1} \end{pmatrix} \overset{(5.8)}{=} \begin{pmatrix} a^{2k+1} & \bar{a}^{2k+1} \\ A^{2k+1} & \bar{A}^{2k+1} \end{pmatrix} \cdot \begin{pmatrix} n_1 \\ \vdots \\ n_{k+1} \end{pmatrix},$$

i.e. the kernel of $\hat{L}_{k+1} \in L[N_0 \times \cdots \times N_k, \bar{B}]$ is given by the range of the matrix operator in brackets applied to $N_1 \times \cdots \times N_{k+1}$, where the blue marked column of this operator was abbreviated in (6.14) by $\bar{m}$.

Now, in the last step, fix $\bar{n}_{k+1} \notin V$ and use the split condition $B = U \oplus V$ to obtain $\bar{n}_{k+1} = u_0 + v_0$, $u_0 \neq 0$. Further, considering (6.16), we set

$$A_\varepsilon \cdot \bar{m} \cdot \bar{n}_{k+1} = \varepsilon^{k+1} \cdot [\, \tfrac{1}{(k+1)!} \cdot I_B + \varepsilon \cdot L_\varepsilon \,] \cdot \bar{n}_{k+1} =: \tfrac{\varepsilon^{k+1}}{(k+1)!} \cdot b_\varepsilon$$

with $b_\varepsilon \in B$ and $b_0 = \bar{n}_{k+1} = u_0 + v_0$. Next, it is easy to see by implicit function theorem that $b_\varepsilon$ can be splitted according to $b_\varepsilon = u_\varepsilon + [\, I_B + \varepsilon \cdot Q(\varepsilon) \,] \cdot v_\varepsilon$ with $u_\varepsilon \in U$ and $v_\varepsilon \in V$ implying

$$A_\varepsilon \cdot \bar{m} \cdot \bar{n}_{k+1} = \tfrac{\varepsilon^{k+1}}{(k+1)!} \cdot \{\, u_\varepsilon + [\, I_B + \varepsilon \cdot Q(\varepsilon) \,] \cdot v_\varepsilon \,\} .$$

Then, using assumptions (3.15) and (3.16), we end up with

$$\| G'[\, z_0(\varepsilon) \,] \cdot A_\varepsilon \cdot \bar{m} \cdot \bar{n}_{k+1} \| = \tfrac{|\varepsilon|^{k+1}}{(k+1)!} \cdot \| G'[\, z_0(\varepsilon) \,] \cdot \{\, u_\varepsilon + [\, I_B + \varepsilon \cdot Q(\varepsilon) \,] \cdot v_\varepsilon \,\} \|$$

$$\geq \tfrac{|\varepsilon|^{k+1}}{(k+1)!} \cdot \left| \underbrace{\| G'[\, z_0(\varepsilon) \,] \cdot u_\varepsilon \|}_{\geq c_1 \cdot |\varepsilon|^k \, \|u_\varepsilon\| \; by \; (3.15)} - \underbrace{\| G'[\, z_0(\varepsilon) \,] \cdot [\, I_B + \varepsilon \cdot Q(\varepsilon) \,] \cdot v_\varepsilon \|}_{=0 \; by \; (3.16)} \right| \tag{6.19}$$

$$\geq \underbrace{c_1}_{>0} \cdot \underbrace{\|u_\varepsilon\|}_{>0} \cdot \tfrac{|\varepsilon|^{2k+1}}{(k+1)!} \geq \alpha_1 \cdot |\varepsilon|^{2k+1}$$

with $\alpha_1 > 0$ chosen sufficiently small, yielding the left hand side of contradiction (6.15).

Finally, from (6.18) evaluated with $n_{k+1} = \bar{n}_{k+1}$, we obtain the right hand side of (6.15) according to

$$\| G'[\, z_0(\varepsilon) \,] \cdot A_\varepsilon \cdot \bar{m} \cdot \bar{n}_{k+1} \| \leq \underbrace{d_2}_{\geq 0} \cdot |\varepsilon|^{k+1} \cdot \| A_\varepsilon \cdot \bar{m} \cdot \bar{n}_{k+1} \|$$

$$\leq\ d_2 \cdot |\varepsilon|^{k+1} \cdot \underbrace{\| A_\varepsilon \|}_{\leq d_3 \cdot |\varepsilon|^{k+1}} \cdot \underbrace{\| \bar{m} \cdot \bar{n}_{k+1} \|}_{>0} \ \leq\ d_2 \cdot |\varepsilon|^{k+1} \cdot \underbrace{d_3}_{\geq 0} \cdot |\varepsilon|^{k+1} \cdot \| \bar{m} \cdot \bar{n}_{k+1} \| \tag{6.20}$$

$$\leq\ \underbrace{\alpha_2}_{\geq 0} \cdot |\varepsilon|^{2k+2} ,$$

finishing the proof back from (ii) to (i).

**Proof of Corollary 1 (ii)→(iii) :** First by definition (3.12) of $P^k(\varepsilon)$, we have

$$G'[\, z_0(\varepsilon)\, ] = P^k(\varepsilon) + \varepsilon^{k+1} \cdot R(\varepsilon) \tag{6.21}$$

with a smooth family of bounded linear operators $R(\varepsilon)$. Now, from (3.15) and the triangle inequality one obtains

$$\| G'[\, z_0(\varepsilon)\, )\, ] \cdot u \,\| \geq \underbrace{c_1}_{>0} \cdot |\varepsilon|^k \cdot \| u \|$$

$$\left\| P^k(\varepsilon) \cdot u \right\| + |\varepsilon|^{k+1} \cdot \| R(\varepsilon) \cdot u \| \ \geq \| P^k(\varepsilon) \cdot u \ + \varepsilon^{k+1} \cdot R(\varepsilon) \cdot u \| \geq\ c_1 \cdot |\varepsilon|^k \cdot \| u \|$$

implying

$$\left\| P^k(\varepsilon) \cdot u \right\| \geq\ c_1 \cdot |\varepsilon|^k \cdot \| u \| -\ |\varepsilon|^{k+1} \cdot \| R(\varepsilon) \cdot u \|$$

$$\geq\ c_1 \cdot |\varepsilon|^k \cdot \| u \| -\ |\varepsilon|^{k+1} \cdot \| R(\varepsilon) \| \cdot \| u \| = |\varepsilon|^k \cdot (\, c_1 - |\varepsilon| \cdot \| R(\varepsilon) \|\, ) \cdot \| u \|$$

$$\geq\ \underbrace{c_2}_{>0} \cdot |\varepsilon|^k \cdot \| u \| ,$$

as well as the first claim (3.17) of (iii). Next, from (3.16) and (6.21) we see

$$0 = G'[\, z_0(\varepsilon)\, ] \cdot [\, I_B + \varepsilon \cdot Q(\varepsilon)\, ]|_V \ =\ [\, P^k(\varepsilon) + \varepsilon^{k+1} \cdot R(\varepsilon)\, ] \cdot [\, I_B + \varepsilon \cdot Q(\varepsilon)\, ]|_V$$

$$\Longrightarrow \quad \| P^k(\varepsilon)[\, I_B + \varepsilon \cdot Q(\varepsilon)\, ]|_V \,\| = |\varepsilon|^{k+1} \cdot\ \| R(\varepsilon) \cdot [\, I_B + \varepsilon \cdot Q(\varepsilon)\, ]|_V \,\| \ \leq\ \ c_3 \cdot |\varepsilon|^{k+1},$$

yielding the second claim (3.18) and (iii) is shown.

**Proof of Corollary 1 (iii)→(i) :** For closing the loop, we prefer to derive (i) from (iii). Essentially, this is merely a repetition of the proof from (ii) to (i). First, consider (3.17) to obtain

$$\left\| P^k(\varepsilon) \cdot u \right\| \geq \underbrace{c_4}_{>0} \cdot |\varepsilon|^k \cdot \| u \|$$

$$\Longrightarrow \quad \| G'[\, z_0(\varepsilon)\, )\, ] \cdot u \,\| = \| P^k(\varepsilon) \cdot u \ + \varepsilon^{k+1} \cdot R(\varepsilon) \cdot u \|$$

$$\geq \left| \underbrace{\| P^k(\varepsilon) \cdot u \,\|}_{\geq c_4 \cdot |\varepsilon|^k \| u \|} - |\varepsilon|^{k+1} \cdot \underbrace{\| R(\varepsilon) \cdot u \,\|}_{\leq \| R(\varepsilon) \| \cdot \| u \|} \right|$$

$$\geq |\varepsilon|^k \cdot (\, c_4 - |\varepsilon| \cdot \| R(\varepsilon) \|\, ) \cdot \| u \| \geq \underbrace{c_5}_{>0} \cdot |\varepsilon|^k \cdot \| u \|$$

as well as (3.15) from (ii). Now, we enter the proof from (ii) to (i) at (6.13) and follow it up until (6.19), where the assumption (3.16) of (ii) is needed for the first time, which has now to be replaced by the assumption (3.18) from (iii). Hence, we have to repeat (6.19) using a slight variation. In some more detail, we start with the estimation

$$\| G'[\,z_0(\varepsilon)\,] \cdot A_\varepsilon \cdot \bar{m} \cdot \bar{n}_{k+1} \| = \frac{|\varepsilon|^{k+1}}{(k+1)!} \cdot \| G'[\,z_0(\varepsilon)\,] \cdot \{\, u_\varepsilon + [\, I_B + \varepsilon \cdot Q(\varepsilon)\,] \cdot v_\varepsilon \,\} \|$$

$$\geq \frac{|\varepsilon|^{k+1}}{(k+1)!} \cdot \big|\, \| G'[\,z_0(\varepsilon)\,] \cdot u_\varepsilon \| - \| G'[\,z_0(\varepsilon)\,] \cdot [\, I_B + \varepsilon \cdot Q(\varepsilon)\,] \cdot v_\varepsilon \| \,\big| \tag{6.22}$$

$$\overset{(6.21)}{=} \frac{|\varepsilon|^{k+1}}{(k+1)!} \cdot \left| \underbrace{\| G'[\,z_0(\varepsilon)\,] \cdot u_\varepsilon \|}_{\geq c_6 \cdot |\varepsilon|^k \| u_\varepsilon \| \text{ by } (3.15)} - \| [\, P^k(\varepsilon) + \varepsilon^{k+1} \cdot R(\varepsilon)] \cdot [\, I_B + \varepsilon \cdot Q(\varepsilon)\,] \cdot v_\varepsilon \| \right| .$$

Now, we use (3.18) from (iii) according to

$$\| [\, P^k(\varepsilon) + \varepsilon^{k+1} \cdot R(\varepsilon)] \cdot [\, I_B + \varepsilon \cdot Q(\varepsilon)\,] \cdot v_\varepsilon \|$$

$$\leq \| P^k(\varepsilon) \cdot [\, I_B + \varepsilon \cdot Q(\varepsilon)\,] \cdot v_\varepsilon \| + |\varepsilon|^{k+1} \cdot \| R(\varepsilon) \cdot [\, I_B + \varepsilon \cdot Q(\varepsilon)\,] \cdot v_\varepsilon \| \tag{6.23}$$

$$\overset{(3.18)}{\leq} c_7 \cdot |\varepsilon|^{k+1} \cdot \| v_\varepsilon \| + |\varepsilon|^{k+1} \cdot \| R(\varepsilon) \cdot [\, I_B + \varepsilon \cdot Q(\varepsilon)\,] \| \cdot \| v_\varepsilon \|$$

$$\leq c_8 \cdot |\varepsilon|^{k+1} \cdot \| v_\varepsilon \|$$

Now, combining (6.22) and(6.23), we end up with

$$\| G'[\,z_0(\varepsilon)\,] \cdot A_\varepsilon \cdot \bar{m} \cdot \bar{n}_{k+1} \| \geq \frac{|\varepsilon|^{k+1}}{(k+1)!} \cdot \left[ \underbrace{c_6}_{>0} \cdot |\varepsilon|^k \cdot \underbrace{\| u_\varepsilon \|}_{>0} - \underbrace{c_8}_{\geq 0} \cdot |\varepsilon|^{k+1} \cdot \| v_\varepsilon \| \right]$$

$$\geq \alpha_1 \cdot |\varepsilon|^{2k+1}$$

with $\alpha_1 > 0$ chosen sufficiently small, thus repeating (6.19) as desired. The right hand side of contradiction (6.15) follows from (6.20), again without change, yielding (i) from (iii).

## *7. Proof of Lemma 1*

Let us now turn to the proof of Lemma 1.

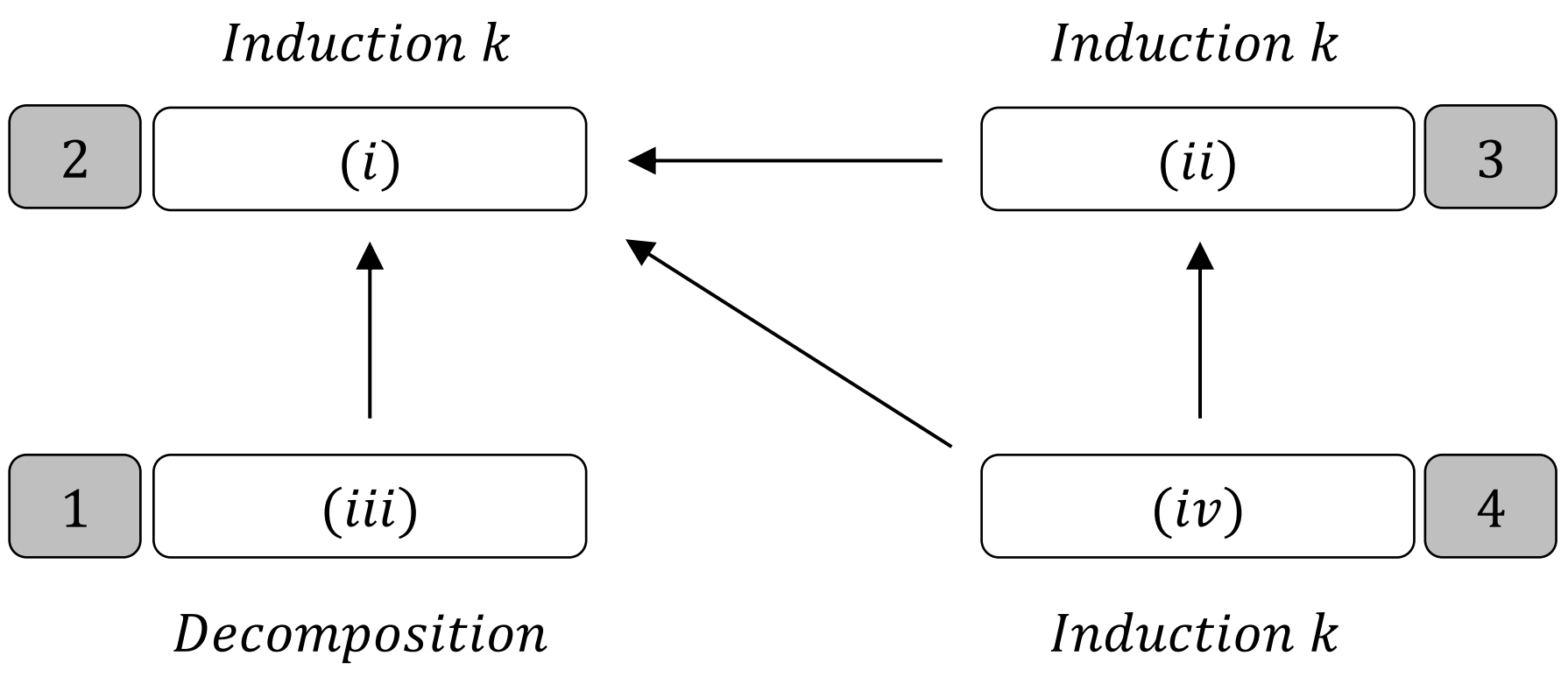

Figure 6 : The way to prove Lemma 1.

First, we obtain (iii) by direct calculation using decomposition (5.7) and some essential definitions of the iteration process within section 5. Next, (i) is obtained by an induction argument with respect to $k$, relying on the validity of (ii), (iii) and (iv). Thirdly, (ii) is proved, again by induction as well as using (iv) as principal argument. The most technical part (iv) is shown by induction on $k$, whereas (v) is obtained by a simple calculation using (2.5) and (i).

**Proof of Lemma 1 (iii) :** Dropping arguments, we obtain by decomposition (5.7)

$$[\,\bar{S}_1 \,|\cdots|\, \bar{S}_{k+1}] \cdot \begin{pmatrix} n_0 \\ \vdots \\ n_k \end{pmatrix} = 0 \in \bar{B}$$

$$\Leftrightarrow \quad P_{R_i}[\,\bar{S}_1 \,|\cdots|\, \bar{S}_{k+1}] \cdot \begin{pmatrix} n_0 \\ \vdots \\ n_k \end{pmatrix} = 0\,, i = 1, \dots, k+1 \quad \wedge \quad P_{R_{k+1}^c}[\,\bar{S}_1\,, \cdots, \bar{S}_{k+1}\,] \cdot \begin{pmatrix} n_0 \\ \vdots \\ n_k \end{pmatrix} = 0$$

$$\Leftrightarrow \quad \begin{pmatrix} \boxed{S_1} & P_{R_1}\bar{S}_2 & \cdots & P_{R_1}\bar{S}_{k+1} \\ & \ddots & \ddots & \vdots \\ & & \boxed{S_k} & P_{R_k}\bar{S}_{k+1} \\ & & & \boxed{S_{k+1}} \end{pmatrix} \cdot \begin{pmatrix} n_0 \\ \vdots \\ n_{k-1} \\ n_k \end{pmatrix} = 0$$

$$\wedge \quad \bar{S}_1 n_0 + \cdots + \bar{S}_{k+1} n_k \in R_1 \oplus \cdots \oplus R_{k+1}\ ,$$

where the definitions of linear $S$-mappings in (5.6) imply $\bar{S}_1 n_0 \in R_1, \dots, \bar{S}_{k+1} n_k \in R_1 \oplus \cdots \oplus R_{k+1}$ and we can restrict to the matrix operator equation. Then, by the definitions of the operators $E_{i,j}$ in (5.8), the following equivalences result from bottom up solution of the triangular system

$$\Leftrightarrow \quad \begin{cases} n_k &= \bar{n}_{k+1} &= E_{k+1,k+1} \cdot \bar{n}_{k+1}\,, & \bar{n}_{k+1} \in N_{k+1} \\ n_{k-1} &= \bar{n}_k - S_k^{-1} P_{R_k} \bar{S}_{k+1} \cdot \bar{n}_{k+1} &= \begin{bmatrix} E_{k,k} & E_{k,k+1} \end{bmatrix} \cdot \begin{pmatrix} \bar{n}_k \\ \bar{n}_{k+1} \end{pmatrix}, & \bar{n}_k \in N_k \\ \vdots & \vdots & \vdots & \vdots \end{cases}$$

$$\Leftrightarrow \quad \begin{pmatrix} n_0 \\ \vdots \\ n_k \end{pmatrix} = \begin{pmatrix} E_{1,1} & \cdots & E_{1,k+1} \\ & \ddots & \vdots \\ & & E_{k+1,k+1} \end{pmatrix} \cdot \begin{pmatrix} \bar{n}_1 \\ \vdots \\ \bar{n}_{k+1} \end{pmatrix}, \quad \begin{pmatrix} \bar{n}_1 \\ \vdots \\ \bar{n}_{k+1} \end{pmatrix} \in N_1 \times \cdots \times N_{k+1}$$

with the last equivalence inferred from a simple induction argument, thus finishing the proof of Lemma 1 (iii).

For later use, some easy to show preliminaries are added, based on (5.18), (5.15) and (5.3)

$$d_{m,1} = c_{m,1} = 1\,, \quad m \geq 0 \tag{7.1}$$

$$\frac{d_{2m,2+j}}{d_{2m-1,1+j}} = d_{2m,2}\,, \quad \frac{d_{2m+1,2+j}}{d_{2m,1+j}} = d_{2m+1,2}\,, \quad m \geq 1\,, \quad j = 0, \dots, m-1 \tag{7.2}$$

$$\frac{d_{2m+1,3+j}}{d_{2m-1,1+j}} = d_{2m,2} \cdot d_{2m+1,2}\,, \quad m \geq 2\,, \quad j = 0, \dots, m-2 \tag{7.3}$$

$$[\,W_{2m+1}^{2m+1}, \dots, W_{m+2}^{2m+1}\,](z_m, \dots, z_0) = [\,W_{2m}^{2m}, \dots, W_{m+1}^{2m}\,](z_{m-1}, \dots, z_0) \cdot D^{2m}\,, \quad m \geq 2 \tag{7.4}$$

$$W_{m+1}^{2m+1}(z_m,\dots,z_0) = [\, W_m^{2m}(z_{m-1},\dots,z_0) \;+\; \frac{(2m)!}{(m!)^2}\, G_0^2 \cdot z_m \,] \cdot d_{2m,m+1}\,, \qquad m \geq 2 \qquad (7.5)$$

$$[\, W_{2m}^{2m},\dots,W_{m+1}^{2m} \,](z_{m-1},\dots,z_0) = [\, W_{2m-1}^{2m-1},\dots,W_m^{2m-1} \,](z_{m-1},\dots,z_0)\cdot D^{2m-1}, \qquad m \geq 2 \qquad (7.6)$$

Note that the basic formulas (2.6) and (2.7), coupling high and low order derivatives within the system of undetermined coefficients, are a direct consequence of (7.1)-(7.6). This follows by straightforward calculation.

**Proof of Lemma 1 (i) :** Next, according to figure 6, (i) is derived from (ii), (iii) and (iv) by induction on $k \geq 1$. In case of $k = 1$, (i) is valid according to

$$T^2[z_2, z_1] \overset{(2.2)}{=} G_0^1 \cdot z_2 \;+\; G_0^2 \cdot z_1^2 = \Delta^1 \cdot z_2 + I^1(z_1)$$

$$M^3 \overset{(5.5)}{=} I_B \quad \textit{and} \quad N_1 = N[G_0^1]\,.$$

Now, suppose (i) with $k \geq 1$. We show the validity for $k+1$. According to (5.14), (5.15), (5.17) and (5.18), the coefficients $T^{2k+1},\dots,T^{k+1}$ can be written in the following form

$$(\, T^{2(k+1)} \quad T^{2k+1} \quad T^{2k} \quad \cdots \quad T^{k+2} \quad T^{k+1} \,)^T =$$

$$\underbrace{\left(\begin{array}{ccccc|c}
W_{2(k+1)}^{2(k+1)} & W_{2k+1}^{2(k+1)} & \cdots & \cdots & W_{k+2}^{2(k+1)} & W_{k+1}^{2(k+1)} \\
 & W_{2k+1}^{2k+1} & \cdots & \cdots & W_{k+2}^{2k+1} & W_{k+1}^{2k+1} \\
 & & W_{2k}^{2k} & \cdots & W_{k+2}^{2k} & W_{k+1}^{2k} \\
 & & & \ddots & \vdots & \vdots \\
 & & & & W_{k+2}^{k+2} & W_{k+1}^{k+2} \\
- & - & - & - & - & - \\
 & & & & & W_{k+1}^{k+1}
\end{array}\right)}_{=:\,\boxed{W}} \cdot \begin{pmatrix} z_{2(k+1)} \\ z_{2k+1} \\ z_{2k} \\ \vdots \\ z_{k+2} \\ - \\ z_{k+1} \end{pmatrix} \;+\; \bar{R}(z_k,\dots z_1)\,.$$

Note that the red marked elements of the matrix $\boxed{W}$ denote $\Delta^k(z_{k-1},\dots,z_1)$, whereas the matrix $\boxed{W}$ without last row and last column delivers $\Delta^{k+1}(z_k,\dots,z_1)$.

By induction hypothesis, the kernel of $\Delta^k$ is characterized by the range of $M^{2k+1}$ and our aim is to characterize the kernel of $\Delta^{k+1}$ by $M^{2k+3}$.

First, we look for elements in $N[\Delta^k] = R\big[M^{2k+1}|_{N_1\times\cdots\times N_k}\big]$, satisfying additionally the equation

$$(\, W_{2k+1}^{2k+1} \quad W_{2k}^{2k+1} \quad \cdots \quad W_{k+1}^{2k+1} \,)\cdot \left( \begin{array}{c} z_{2k+1} \\ M^{2k+1}\cdot \begin{pmatrix} n_1 \\ \vdots \\ n_k \end{pmatrix} \end{array} \right) = 0\,, \qquad (7.7)$$

or equivalently by setting $z_{2k+1} = n_0 \in B$ and use of Lemma 1 (ii), (iii)

$$\underbrace{\left( W_{2k+1}^{2k+1} \quad W_{2k}^{2k+1} \quad \cdots \quad W_{k+1}^{2k+1} \right)}_{=W^{2k+1}} \cdot \underbrace{\begin{pmatrix} I_B & \\ & M^{2k+1} \end{pmatrix}}_{=\widehat{M}_{k+1}} \cdot \begin{pmatrix} n_0 \\ n_1 \\ \vdots \\ n_k \end{pmatrix} = 0$$

$$\overset{(ii)}{\Leftrightarrow} \quad [\, \bar{S}_1 \mid \bar{S}_2(z_1) \mid \cdots \mid \bar{S}_{k+1}(z_k, \dots, z_1) \,] \cdot C^{2k+1} \cdot \begin{pmatrix} n_0 \\ n_1 \\ \vdots \\ n_k \end{pmatrix} = 0$$

$$\overset{(iii)}{\Leftrightarrow} \begin{pmatrix} n_0 \\ n_1 \\ \vdots \\ n_k \end{pmatrix} = \left(C^{2k+1}\right)^{-1} \cdot E^{k+1} \cdot \begin{pmatrix} \bar{n}_1 \\ \bar{n}_2 \\ \vdots \\ \bar{n}_{k+1} \end{pmatrix}, \quad \bar{n}_i \in N_i\,, \quad 1 \le i \le k+1$$

$$\Leftrightarrow \begin{pmatrix} n_0 \\ n_1 \\ \vdots \\ n_k \end{pmatrix} = \left(C^{2k+1}\right)^{-1} \cdot E^{k+1} \cdot C^{2k+1} \cdot \begin{pmatrix} \bar{n}_1 \\ \bar{n}_2 \\ \vdots \\ \bar{n}_{k+1} \end{pmatrix} \overset{(5.8)}{=} \begin{pmatrix} a^{2k+1} & \bar{a}^{2k+1} \\ A^{2k+1} & \bar{A}^{2k+1} \end{pmatrix} \cdot \begin{pmatrix} \bar{n}_1 \\ \bar{n}_2 \\ \vdots \\ \bar{n}_{k+1} \end{pmatrix}.$$

Hence, the elements from $N[\Delta^k]$ that can be extended to solutions of (7.7) are exactly given by

$$\begin{pmatrix} z_{2k+1} \\ z_{2k} \\ \vdots \\ z_{k+1} \end{pmatrix} = \begin{pmatrix} I_B & \\ & M^{2k+1} \end{pmatrix} \cdot \begin{pmatrix} a^{2k+1} & \bar{a}^{2k+1} \\ A^{2k+1} & \bar{A}^{2k+1} \end{pmatrix} \cdot \begin{pmatrix} n_1 \\ n_2 \\ \vdots \\ n_{k+1} \end{pmatrix} = \begin{pmatrix} I_B & * & * \\ & \ddots & * \\ & & I_B \end{pmatrix} \cdot \begin{pmatrix} n_1 \\ n_2 \\ \vdots \\ n_{k+1} \end{pmatrix},$$

where all matrices are of upper tridiagonal type with identity operator in the diagonal. Now, setting $n_{k+1} = 0$ within these solutions, we obtain the kernel of the matrix operator

$$\overline{\Delta}^k := \begin{pmatrix} W_{2k+1}^{2k+1} & \cdots & \cdots & W_{k+2}^{2k+1} \\ & W_{2k}^{2k} & \cdots & W_{k+2}^{2k} \\ & & \ddots & \vdots \\ & & & W_{k+2}^{k+2} \end{pmatrix}$$

according to

$$\begin{pmatrix} z_{2k+1} \\ \vdots \\ z_{k+2} \end{pmatrix} = \begin{pmatrix} I_B & \\ & M^{2k+1} \end{pmatrix}_{\overline{k+1}} \begin{pmatrix} a^{2k+1} \\ A^{2k+1} \end{pmatrix} \cdot \begin{pmatrix} n_1 \\ \vdots \\ n_k \end{pmatrix}. \qquad (7.8)$$

For further reasoning, the diagonal matrices $D^i$ and $C^i$ from (5.4) without the first element $D_{1,1}^i = C_{1,1}^i = 1$ are denoted by $\bar{D}^i$ and $\bar{C}^i$ respectively. Then, we look for elements in $N[\overline{\Delta}^k]$, satisfying additionally the equation

$$\left( W_{2(k+1)}^{2(k+1)} \quad W_{2k+1}^{2(k+1)} \quad \cdots \quad W_{k+2}^{2(k+1)} \right) \cdot \begin{pmatrix} z_{2(k+1)} \\ \begin{pmatrix} I_B & \\ & M^{2k+1} \end{pmatrix}_{\overline{k+1}} \begin{pmatrix} a^{2k+1} \\ A^{2k+1} \end{pmatrix} \cdot \begin{pmatrix} n_1 \\ \vdots \\ n_k \end{pmatrix} \end{pmatrix} = 0\,, \qquad (7.9)$$

or equivalently by setting $z_{2(k+1)} = n_0 \in B$ and use of (7.6), (5.3) with $m = k+1$ as well as Lemma 1 (ii), (iv)

$$\left( W^{2(k+1)}_{2(k+1)} \quad W^{2(k+1)}_{2k+1} \quad \cdots \quad W^{2(k+1)}_{k+2} \right) \cdot \begin{pmatrix} I_B & \\ & a^{2k+1} \\ & M^{2k+1} \cdot A^{2k+1} \end{pmatrix}_{\overline{k+2}} \begin{pmatrix} n_0 \\ n_1 \\ \vdots \\ n_k \end{pmatrix}$$

$$\overset{(7.6)}{=} \underbrace{\left( W^{2k+1}_{2k+1} \quad W^{2k+1}_{2k} \quad \cdots \quad W^{2k+1}_{k+1} \right)}_{=W^{2k+1}} \cdot D^{2k+1} \cdot \begin{pmatrix} I_B & \\ & a^{2k+1} \\ & M^{2k+1} \cdot A^{2k+1} \end{pmatrix}_{\overline{k+2}} \begin{pmatrix} n_0 \\ n_1 \\ \vdots \\ n_k \end{pmatrix}$$

$$\overset{(ii)}{=} \overbrace{[\,\bar{S}_1 \mid \bar{S}_2(z_1) \mid \cdots \mid \bar{S}_{k+1}(z_k, \dots, z_1)\,] \cdot C^{2k+1} \cdot \begin{pmatrix} I_B & \\ & M^{2k+1,-1} \end{pmatrix}}^{=W^{2k+1}}$$

$$\cdot D^{2k+1} \cdot \begin{pmatrix} I_B & \\ & a^{2k+1} \\ & M^{2k+1} \cdot A^{2k+1} \end{pmatrix}_{\overline{k+2}} \begin{pmatrix} n_0 \\ n_1 \\ \vdots \\ n_k \end{pmatrix}$$

$$\overset{(iv)}{=} [\,\bar{S}_1 \mid \bar{S}_2(z_1) \mid \cdots \mid \bar{S}_{k+1}(z_k, \dots, z_1)\,] \cdot \begin{pmatrix} I_B & \\ & \bar{C}^{2k+1} \cdot M^{2k+1,-1} \end{pmatrix}$$

$$\cdot \begin{pmatrix} I_B & \\ & \underbrace{\bar{D}^{2k+1} \cdot \bar{D}^{2k+1,-1}}_{=I} \cdot M^{2k+1} \cdot \bar{D}^{2k+1} \end{pmatrix} \cdot \begin{pmatrix} n_0 \\ \vdots \\ n_k \end{pmatrix}$$

$$= [\,\bar{S}_1 \mid \bar{S}_2(z_1) \mid \cdots \mid \bar{S}_{k+1}(z_k, \dots, z_1)\,] \cdot \begin{pmatrix} I_B & \\ & \bar{C}^{2k+1} \cdot \bar{D}^{2k+1} \end{pmatrix} \cdot \begin{pmatrix} n_0 \\ \vdots \\ n_k \end{pmatrix}$$

$$\overset{(5.3)}{=} [\,\bar{S}_1 \mid \bar{S}_2(z_1) \mid \cdots \mid \bar{S}_{k+1}(z_k, \dots, z_1)\,] \cdot C^{2k+2} \cdot \begin{pmatrix} n_0 \\ \vdots \\ n_k \end{pmatrix} = 0\,.$$

Now, again by Lemma 1 (iii), we obtain

$$\begin{pmatrix} n_0 \\ \vdots \\ n_k \end{pmatrix} = \left(C^{2k+2}\right)^{-1} \cdot E^{k+1} \cdot \begin{pmatrix} \bar{n}_1 \\ \vdots \\ \bar{n}_{k+1} \end{pmatrix}, \quad \bar{n}_i \in N_i\,, \quad 1 \le i \le k+1$$

$$\Leftrightarrow \begin{pmatrix} n_0 \\ \vdots \\ n_k \end{pmatrix} = \left(C^{2k+2}\right)^{-1} \cdot E^{k+1} \cdot C^{2k+2} \cdot \begin{pmatrix} \bar{n}_1 \\ \vdots \\ \bar{n}_{k+1} \end{pmatrix} \overset{(5.8)}{=} \begin{pmatrix} a^{2k+2} & \bar{a}^{2k+2} \\ A^{2k+2} & \bar{A}^{2k+2} \end{pmatrix} \cdot \begin{pmatrix} \bar{n}_1 \\ \vdots \\ \bar{n}_{k+1} \end{pmatrix}$$

and the elements from $N[\bar{\Delta}^k]$ that can be extended to $N[\Delta^{k+1}]$ are exactly given by

$$\begin{pmatrix} z_{2k+2} \\ z_{2k+1} \\ \vdots \\ z_{k+2} \end{pmatrix} \overset{(7.8)}{=} \begin{pmatrix} I_B & & \\ & a^{2k+1} & \\ & & M^{2k+1} \cdot A^{2k+1} \end{pmatrix}_{\overline{k+2}} \begin{pmatrix} a^{2k+2} & \bar{a}^{2k+2} \\ A^{2k+2} & \bar{A}^{2k+2} \end{pmatrix} \cdot \begin{pmatrix} \bar{n}_1 \\ \vdots \\ \bar{n}_{k+1} \end{pmatrix}$$

$$\overset{(5.9)}{=} \begin{pmatrix} (a^{2k+2} \quad \bar{a}^{2k+2}) \\ a^{2k+1} \cdot (A^{2k+2} \quad \bar{A}^{2k+2}) \\ M^{2k+1} \cdot A^{2k+1} \cdot (A^{2k+2} \quad \bar{A}^{2k+2}) \end{pmatrix}_{\overline{k+2}} \cdot \begin{pmatrix} \bar{n}_1 \\ \vdots \\ \bar{n}_{k+1} \end{pmatrix} = M^{2k+3} \cdot \begin{pmatrix} \bar{n}_1 \\ \vdots \\ \bar{n}_{k+1} \end{pmatrix},$$

thus finishing the induction.

**Proof of Lemma 1 (ii) :** In the next step, (ii) is accomplished by induction for $k \geq 1$ using (iv). In case of $k = 1$, the identity (ii) reads

$$W^3(z_1) \cdot \begin{pmatrix} I_B & \\ & M^3 \end{pmatrix} = [\, \bar{S}_1 \mid \bar{S}_2(z_1) \,] \cdot C^3,$$

yielding by (2.2), (5.17), (5.5), (5.7) and (5.3)

$$[\, G_0^1 \mid 3G_0^2 z_1 \,] \cdot \begin{pmatrix} I_B & \\ & I_B \end{pmatrix} = [\, G_0^1 \mid 2G_0^2 z_1 \,] \cdot \begin{pmatrix} 1 & \\ & \frac{3}{2} \end{pmatrix},$$

which is true by inspection. Now, suppose identity (ii) with $k$ replaced by $k-1 \geq 1$

$$W^{2k-1}(z_{k-1}, \dots, z_1) \cdot \begin{pmatrix} I_B & \\ & M^{2k-1} \end{pmatrix} = [\, \bar{S}_1 \mid \cdots \mid \bar{S}_k(z_{k-1}, \dots, z_1) \,] \cdot C^{2k-1}\,. \tag{7.10}$$

Our aim is to show (ii) with $k$, where the equality with respect to the first component is obviously valid under consideration of (5.21), (7.1) and it remains to show

$$\left[\, W_{2k}^{2k+1}, \dots, W_{k+1}^{2k+1} \,\right](z_k, \dots, z_1) \cdot M^{2k+1} \tag{7.11}$$

$$= [\, \bar{S}_2 \mid \cdots \mid \bar{S}_{k+1} \,](z_k, \dots, z_1) \cdot \begin{pmatrix} c_{2k+1,2} & & \\ & \ddots & \\ & & c_{2k+1,k+1} \end{pmatrix}.$$

Further, with index $k_s$ and $k_z$ denoting column $k$ and row $k$ of a matrix respectively, as well as index $k|$ denoting a matrix without column $k$, the left hand side of (7.11) reads

$$\{ \left[\, W_{2k}^{2k+1}, \dots, W_{k+1}^{2k+1} \,\right](z_k, \dots, z_1) \cdot M_{k|}^{2k+1}\,,\ \left[\, W_{2k}^{2k+1}, \dots, W_{k+1}^{2k+1} \,\right](z_k, \dots, z_1) \cdot M_{k_s}^{2k+1} \}$$

$$\overset{(5.12)}{=} \{ \left[\, W_{2k}^{2k+1}, \dots, W_{k+2}^{2k+1} \,\right](z_k, \dots, z_1) \cdot M_{k|,\bar{k}}^{2k+1}\,,$$

$$\left[\, W_{2k}^{2k+1}, \dots, W_{k+2}^{2k+1} \,\right](z_k, \dots, z_1) \cdot M_{k_s,\bar{k}}^{2k+1} + W_{k+1}^{2k+1}(z_k, \dots, z_1) \cdot I_B \} \tag{7.12}$$

$$\overset{(7.4)}{\overset{(7.5)}{=}} \{ [\, W_{2k-1}^{2k}, \dots, W_{k+1}^{2k} \,](z_{k-1}, \dots, z_1) \cdot \begin{pmatrix} d_{2k,2} & & \\ & \ddots & \\ & & d_{2k,k} \end{pmatrix} \cdot M_{k|,\bar{k}}^{2k+1} ,$$

$$[\, W_{2k-1}^{2k}, \dots, W_{k+1}^{2k} \,](z_{k-1}, \dots, z_1) \cdot \begin{pmatrix} d_{2k,2} & & \\ & \ddots & \\ & & d_{2k,k} \end{pmatrix} \cdot M_{k_s,\bar{k}}^{2k+1}$$

$$+ \; [\, W_k^{2k}(z_{k-1}, \dots, z_1) + \frac{(2k)!}{(k!)^2} \, G_0^2 z_k \,] \cdot d_{2k,k+1} \, \} \, .$$

The right hand side of (7.11) can be represented by induction hypothesis (7.10) according to

$$\{ [\, \bar{S}_2 \mid \cdots \mid \bar{S}_k \,](z_{k-1}, \dots, z_1) \cdot \begin{pmatrix} c_{2k+1,2} & & \\ & \ddots & \\ & & c_{2k+1,k} \end{pmatrix} , \; \bar{S}_{k+1}(z_k, \dots, z_1) \cdot c_{2k+1,k+1} \, \} \qquad (7.13)$$

$$\overset{(7.10)}{\overset{(5.6)}{\overset{(5.17)}{=}}} \{ [\, W_{2k-2}^{2k-1}, \dots, W_k^{2k-1} \,](z_{k-1}, \dots, z_1) \cdot M^{2k-1} \cdot \begin{pmatrix} c_{2k-1,2} & & \\ & \ddots & \\ & & c_{2k-1,k} \end{pmatrix}^{-1} \cdot \begin{pmatrix} c_{2k+1,2} & & \\ & \ddots & \\ & & c_{2k+1,k} \end{pmatrix} ,$$

$$[\, W_{2k-1}^{2k}, \dots, W_k^{2k} \,](z_{k-1}, \dots, z_1) \cdot \begin{pmatrix} & \bar{a}^{2k-1} \\ M^{2k-1} \cdot & \bar{A}^{2k-1} \end{pmatrix} \cdot c_{2k+1,k+1} \; + \; \frac{(2k)!}{(k!)^2} \, G_0^2 z_k \cdot c_{2k+1,k+1} \, \}$$

$$\overset{(7.6)}{\overset{(5.3)}{=}} \{ [\, W_{2k-1}^{2k}, \dots, W_{k+1}^{2k} \,](z_{k-1}, \dots, z_1) \cdot \begin{pmatrix} d_{2k-1,2} & & \\ & \ddots & \\ & & d_{2k-1,k} \end{pmatrix}^{-1} \cdot M^{2k-1} \cdot \begin{pmatrix} c_{2k-1,2} & & \\ & \ddots & \\ & & c_{2k-1,k} \end{pmatrix}^{-1}$$

$$\cdot \begin{pmatrix} c_{2k-1,2} & & \\ & \ddots & \\ & & c_{2k-1,k} \end{pmatrix} \cdot \begin{pmatrix} d_{2k-1,2} & & \\ & \ddots & \\ & & d_{2k-1,k} \end{pmatrix} \cdot \begin{pmatrix} d_{2k,2} & & \\ & \ddots & \\ & & d_{2k,k} \end{pmatrix} ,$$

$$[\, W_{2k-1}^{2k}, \dots, W_{k+1}^{2k} \,](z_{k-1}, \dots, z_1) \cdot \begin{pmatrix} & \bar{a}^{2k-1} \\ M^{2k-1} \cdot & \bar{A}^{2k-1} \end{pmatrix}_{\bar{k}} \cdot \overbrace{c_{2k,k+1}}^{=1} \cdot d_{2k,k+1}$$

$$+ \; [\, W_k^{2k}(z_{k-1}, \dots, z_1) \; + \; \frac{(2k)!}{(k!)^2} \, G_0^2 z_k \,] \cdot 1 \cdot d_{2k,k+1} \, \} \, .$$

Comparing (7.12) and (7.13), sufficient conditions for equality within (7.11) are given by

$$\begin{pmatrix} d_{2k,2} & & \\ & \ddots & \\ & & d_{2k,k} \end{pmatrix} \cdot M_{k|,\bar{k}}^{2k+1}$$

$$= \begin{pmatrix} d_{2k-1,2} & & \\ & \ddots & \\ & & d_{2k-1,k} \end{pmatrix}^{-1} \cdot M^{2k-1} \cdot \begin{pmatrix} d_{2k-1,2} & & \\ & \ddots & \\ & & d_{2k-1,k} \end{pmatrix} \cdot \begin{pmatrix} d_{2k,2} & & \\ & \ddots & \\ & & d_{2k,k} \end{pmatrix}$$

as well as

$$\begin{pmatrix} d_{2k,2} & & \\ & \ddots & \\ & & d_{2k,k} \end{pmatrix} \cdot M^{2k+1}_{k_s,\bar{k}} = \begin{pmatrix} & \bar{a}^{2k-1} \\ M^{2k-1} \cdot & \bar{A}^{2k-1} \end{pmatrix}_{\bar{k}} \cdot d_{2k,k+1}$$

or combined under consideration of (iv) with $k$ replaced by $k-1$ and (5.8), (5.9)

$$\begin{pmatrix} d_{2k,2} & & \\ & \ddots & \\ & & d_{2k,k} \end{pmatrix} \cdot M^{2k+1}_{\bar{k}}$$

$$= [\underbrace{\begin{pmatrix} & a^{2k-1} \\ M^{2k-1} \cdot & A^{2k-1} \end{pmatrix}_{\bar{k}} \cdot}_{by\ (iv)} \begin{pmatrix} d_{2k,2} & & \\ & \ddots & \\ & & d_{2k,k} \end{pmatrix}, \begin{pmatrix} & \bar{a}^{2k-1} \\ M^{2k-1} \cdot & \bar{A}^{2k-1} \end{pmatrix}_{\bar{k}} \cdot d_{2k,k+1}]$$

$$\overset{(5.9)}{\Leftrightarrow} \begin{pmatrix} d_{2k,2} & & \\ & \ddots & \\ & & d_{2k,k} \end{pmatrix} \cdot \left( \overbrace{\begin{pmatrix} & & (a^{2k} & \bar{a}^{2k}) \\ & a^{2k-1} \cdot & (A^{2k} & \bar{A}^{2k}) \\ M^{2k-1} \cdot & A^{2k-1} \cdot & (A^{2k} & \bar{A}^{2k}) \end{pmatrix}_{\overline{k+1}}}^{= M^{2k+1}} \right)_{\bar{k}}$$

$$= \begin{pmatrix} & (a^{2k-1} & \bar{a}^{2k-1}) \\ M^{2k-1} \cdot & (A^{2k-1} & \bar{A}^{2k-1}) \end{pmatrix}_{\bar{k}} \cdot \begin{pmatrix} d_{2k,2} & & \\ & \ddots & \\ & & d_{2k,k+1} \end{pmatrix}$$

$$\Leftrightarrow \begin{pmatrix} d_{2k,2} & & \\ & \ddots & \\ & & d_{2k,k} \end{pmatrix} \cdot \begin{pmatrix} & (a^{2k} & \bar{a}^{2k}) \\ \begin{pmatrix} & a^{2k-1} \\ M^{2k-1} \cdot & A^{2k-1} \end{pmatrix}_{\bar{k}} \cdot & (A^{2k} & \bar{A}^{2k}) \end{pmatrix}_{\bar{k}}$$

$$= \begin{pmatrix} & (a^{2k-1} & \bar{a}^{2k-1}) \\ M^{2k-1} \cdot & (A^{2k-1} & \bar{A}^{2k-1}) \end{pmatrix}_{\bar{k}} \cdot \begin{pmatrix} d_{2k,2} & & \\ & \ddots & \\ & & d_{2k,k+1} \end{pmatrix}.$$

Now, using again (iv), we obtain

$$\begin{pmatrix} d_{2k,2} & & \\ & \ddots & \\ & & d_{2k,k} \end{pmatrix} \cdot \begin{pmatrix} & (a^{2k} & \bar{a}^{2k}) \\ \overbrace{\begin{pmatrix} d_{2k-1,2} & & \\ & \ddots & \\ & & d_{2k-1,k} \end{pmatrix}^{-1} \cdot M^{2k-1} \cdot \begin{pmatrix} d_{2k-1,2} & & \\ & \ddots & \\ & & d_{2k-1,k} \end{pmatrix}}^{by\ (iv)} \cdot & (A^{2k} & \bar{A}^{2k}) \end{pmatrix}_{\bar{k}}$$

$$= \begin{pmatrix} \left(a^{2k-1} \quad \bar{a}^{2k-1}\right) \\ M^{2k-1} \cdot \left(A^{2k-1} \quad \bar{A}^{2k-1}\right) \end{pmatrix}_{\bar{k}} \cdot \begin{pmatrix} d_{2k,2} & & \\ & \ddots & \\ & & d_{2k,k+1} \end{pmatrix}. \tag{7.14}$$

Further, (5.3) and (5.8) yield

$$\begin{pmatrix} a^{2k} & \bar{a}^{2k} \\ A^{2k} & \bar{A}^{2k} \end{pmatrix} = \left(C^{2k}\right)^{-1} \cdot E^k \cdot C^{2k} \tag{7.15}$$

$$= \left(D^{2k-1}\right)^{-1} \cdot \left(C^{2k-1}\right)^{-1} \cdot E^k \cdot C^{2k-1} \cdot D^{2k-1} = \left(D^{2k-1}\right)^{-1} \cdot \begin{pmatrix} a^{2k-1} & \bar{a}^{2k-1} \\ A^{2k-1} & \bar{A}^{2k-1} \end{pmatrix} \cdot D^{2k-1}$$

and for later use, we obtain by a similar calculation

$$\begin{pmatrix} a^{2k} & \bar{a}^{2k} \\ A^{2k} & \bar{A}^{2k} \end{pmatrix} = D^{2k} \cdot \begin{pmatrix} a^{2k+1} \\ A^{2k+1} \end{pmatrix}_{\overline{k+1}} \cdot \left(D^{2k}\right)^{-1}. \tag{7.16}$$

Then by (7.15), the first row in (7.14) is equivalent to

$$d_{2k,2} \cdot \overbrace{\left(d_{2k-1,1}\right)^{-1}}^{=1} \cdot \left(a^{2k-1} \quad \bar{a}^{2k-1}\right) \cdot \overbrace{\begin{pmatrix} d_{2k-1,1} & & \\ & \ddots & \\ & & d_{2k-1,k} \end{pmatrix}}^{=D^{2k-1}}$$

$$= \left(a^{2k-1} \quad \bar{a}^{2k-1}\right) \cdot \begin{pmatrix} d_{2k,2} & & \\ & \ddots & \\ & & d_{2k,k+1} \end{pmatrix}$$

$$\Leftrightarrow \quad 0 = \left(a^{2k-1} \quad \bar{a}^{2k-1}\right) \cdot \underbrace{\left[ d_{2k,2} \cdot \begin{pmatrix} d_{2k-1,1} & & \\ & \ddots & \\ & & d_{2k-1,k} \end{pmatrix} - \begin{pmatrix} d_{2k,2} & & \\ & \ddots & \\ & & d_{2k,k+1} \end{pmatrix} \right]}_{=0 \ by \ (7.2)},$$

which is obviously true using (7.2) with $j = 0, \dots, k-1$. Along the same lines of reasoning, the remaining rows $2, \dots k-1$ in (7.14) can be treated with $j = 1, \dots, k-2$ according to

$$d_{2k,2+j} \cdot \left(d_{2k-1,1+j}\right)^{-1} \cdot M_{j_z}^{2k-1} \cdot \overbrace{\left(A^{2k-1} \quad \bar{A}^{2k-1}\right) \cdot \begin{pmatrix} d_{2k-1,1} & & \\ & \ddots & \\ & & d_{2k-1,k} \end{pmatrix}}^{(7.15)}$$

$$= M_{j_z}^{2k-1} \cdot \left(A^{2k-1} \quad \bar{A}^{2k-1}\right) \cdot \begin{pmatrix} d_{2k,2} & & \\ & \ddots & \\ & & d_{2k,k+1} \end{pmatrix}$$

$$\Leftrightarrow \quad 0 = M_{j_z}^{2k-1} \cdot \left(A^{2k-1} \quad \bar{A}^{2k-1}\right) \cdot \underbrace{\Bigg[ \overbrace{\left(\frac{d_{2k,2+j}}{d_{2k-1,1+j}}\right)}^{=d_{2k,2}\ by\ (7.2)} \cdot \begin{pmatrix} d_{2k-1,1} & & \\ & \ddots & \\ & & d_{2k-1,k} \end{pmatrix} - \begin{pmatrix} d_{2k,2} & & \\ & \ddots & \\ & & d_{2k,k+1} \end{pmatrix} \Bigg]}_{=0\ by\ (7.2)}$$

in this way finishing the proof of Lemma 1 (ii).

**Proof of Lemma 1 (iv) :** The proof of (iv) is accomplished independently of (i)-(iii), using again an inductive argument for $k \geq 1$. For $k = 1$, identity (iv) is obviously true by

$$\underbrace{M^3}_{=I_B\ by\ (5.5)} \cdot d_{3,2} = d_{3,2} \cdot \underbrace{a^3}_{=I_B\ by\ (5.5)}$$

Then, suppose (iv) with $k$ replaced by $k-1 \geq 1$ yielding

$$M^{2k-1} \cdot \begin{pmatrix} d_{2k-1,2} & & \\ & \ddots & \\ & & d_{2k-1,k} \end{pmatrix} = \begin{pmatrix} d_{2k-1,2} & & \\ & \ddots & \\ & & d_{2k-1,k} \end{pmatrix} \cdot \begin{pmatrix} a^{2k-1} \\ M^{2k-1} \cdot A^{2k-1} \end{pmatrix}_{\overline{k}}, \tag{7.17}$$

whereas the identity to prove is given by

$$\overbrace{\begin{pmatrix} \left(a^{2k} \quad \bar{a}^{2k}\right) \\ \begin{pmatrix} a^{2k-1} \\ M^{2k-1} \cdot A^{2k-1} \end{pmatrix} \cdot \left(A^{2k} \quad \bar{A}^{2k}\right) \end{pmatrix}_{\overline{k+1}}}^{=M^{2k+1}} \cdot \begin{pmatrix} d_{2k+1,2} & & \\ & \ddots & \\ & & d_{2k+1,k+1} \end{pmatrix} \tag{7.18}$$

$$= \begin{pmatrix} d_{2k+1,2} & & \\ & \ddots & \\ & & d_{2k+1,k+1} \end{pmatrix} \cdot \begin{pmatrix} a^{2k+1} \\ \underbrace{\begin{pmatrix} \left(a^{2k} \quad \bar{a}^{2k}\right) \\ \begin{pmatrix} a^{2k-1} \\ M^{2k-1} \cdot A^{2k-1} \end{pmatrix} \cdot \left(A^{2k} \quad \bar{A}^{2k}\right) \end{pmatrix}_{\overline{k+1}}}_{=M^{2k+1}} \cdot A^{2k+1} \end{pmatrix}_{\overline{k+1}}$$

under consideration of (5.9). Further, using the induction hypothesis (7.17), the right hand side of (7.18) transforms according to

$$\begin{pmatrix} d_{2k+1,2} & & \\ & \ddots & \\ & & d_{2k+1,k+1} \end{pmatrix} \tag{7.19}$$

$$\cdot \begin{pmatrix} & a^{2k+1} \\ & \begin{pmatrix} a^{2k} & \bar{a}^{2k} \end{pmatrix} \cdot A^{2k+1} \\ \begin{pmatrix} d_{2k-1,2} & & \\ & \ddots & \\ & & d_{2k-1,k} \end{pmatrix}^{-1} \cdot M^{2k-1} \cdot \begin{pmatrix} d_{2k-1,2} & & \\ & \ddots & \\ & & d_{2k-1,k} \end{pmatrix} & \cdot \begin{pmatrix} A^{2k} & \bar{A}^{2k} \end{pmatrix} \cdot A^{2k+1} \end{pmatrix}_{\overline{k+1}} .$$

Next, the left hand side in (7.18) is transformed into (7.19), row by row. Concerning the first row, we obtain from the left hand side in (7.18) and the second relation in (7.2)

$$\begin{aligned} &\begin{pmatrix} a^{2k} & \bar{a}^{2k} \end{pmatrix} \cdot \begin{pmatrix} d_{2k+1,2} & & \\ & \ddots & \\ & & d_{2k+1,k+1} \end{pmatrix} \\ &\overset{(7.16)}{=} d_{2k,1} \cdot a^{2k+1} \cdot \begin{pmatrix} d_{2k,1} & & \\ & \ddots & \\ & & d_{2k,k} \end{pmatrix}^{-1} \cdot \begin{pmatrix} d_{2k+1,2} & & \\ & \ddots & \\ & & d_{2k+1,k+1} \end{pmatrix} \\ &\overset{(7.2)}{=} \underbrace{d_{2k,1}}_{=1} \cdot a^{2k+1} \cdot \begin{pmatrix} d_{2k,1} & & \\ & \ddots & \\ & & d_{2k,k} \end{pmatrix}^{-1} \cdot \begin{pmatrix} d_{2k,1} & & \\ & \ddots & \\ & & d_{2k,k} \end{pmatrix} \cdot d_{2k+1,2} \\ &= d_{2k+1,2} \cdot a^{2k+1}, \end{aligned}$$

obviously agreeing with the first row in (7.19). The second row from the left hand side in (7.18) implies

$$\begin{aligned} &a^{2k-1} \cdot \begin{pmatrix} A^{2k} & \bar{A}^{2k} \end{pmatrix} \cdot \begin{pmatrix} d_{2k+1,2} & & \\ & \ddots & \\ & & d_{2k+1,k+1} \end{pmatrix} \\ &\overset{(7.16)}{\underset{(7.15)}{=}} \underbrace{d_{2k-1,1}}_{=1} \cdot a^{2k} \cdot \begin{pmatrix} d_{2k-1,1} & & \\ & \ddots & \\ & & d_{2k-1,k-1} \end{pmatrix}^{-1} \cdot \begin{pmatrix} d_{2k,2} & & \\ & \ddots & \\ & & d_{2k,k} \end{pmatrix} \cdot A_{\bar{k}}^{2k+1} \cdot \begin{pmatrix} d_{2k,1} & & \\ & \ddots & \\ & & d_{2k,k} \end{pmatrix}^{-1} \\ &\qquad \cdot \begin{pmatrix} d_{2k+1,2} & & \\ & \ddots & \\ & & d_{2k+1,k+1} \end{pmatrix} \\ &\overset{(7.2)}{=} a^{2k} \cdot d_{2k,2} \cdot A_{\bar{k}}^{2k+1} \cdot d_{2k+1,2} \\ &\overset{(7.3),\, j=0}{=} d_{2k+1,3} \cdot a^{2k} \cdot A_{\bar{k}}^{2k+1} \\ &= d_{2k+1,3} \cdot \begin{pmatrix} a^{2k} & \bar{a}^{2k} \end{pmatrix} \cdot A^{2k+1}, \end{aligned}$$

yielding the second row in (7.19). Note that the last identity follows from the fact that the $k$-th row of $A^{2k+1}$ equals zero according to the definition in (5.8).

Finally, in case of $k \geq 3$, the remaining rows $3, \ldots k$ in (7.18), (7.19) can be treated along the same lines of reasoning with $j = 1, \ldots, k-2$ according to

$$M_{j_z}^{2k-1} \cdot A^{2k-1} \cdot \begin{pmatrix} A^{2k} & \bar{A}^{2k} \end{pmatrix} \cdot \begin{pmatrix} d_{2k+1,2} & & \\ & \ddots & \\ & & d_{2k+1,k+1} \end{pmatrix}$$

$$\overset{(7.16)}{\underset{(7.15)}{=}} M_{j_z}^{2k-1} \cdot \begin{pmatrix} d_{2k-1,2} & & \\ & \ddots & \\ & & d_{2k-1,k} \end{pmatrix} \cdot A^{2k} \cdot \begin{pmatrix} d_{2k-1,1} & & \\ & \ddots & \\ & & d_{2k-1,k-1} \end{pmatrix}^{-1}$$

$$\cdot \begin{pmatrix} d_{2k,2} & & \\ & \ddots & \\ & & d_{2k,k} \end{pmatrix} \cdot A_{\bar{k}}^{2k+1} \cdot \begin{pmatrix} d_{2k,1} & & \\ & \ddots & \\ & & d_{2k,k} \end{pmatrix}^{-1} \cdot \begin{pmatrix} d_{2k+1,2} & & \\ & \ddots & \\ & & d_{2k+1,k+1} \end{pmatrix}$$

$$\overset{(7.2)}{=} M_{j_z}^{2k-1} \cdot \begin{pmatrix} d_{2k-1,2} & & \\ & \ddots & \\ & & d_{2k-1,k} \end{pmatrix} \cdot A^{2k} \cdot d_{2k,2} \cdot A_{\bar{k}}^{2k+1} \cdot d_{2k+1,2}$$

$$\overset{(7.3)}{=} \frac{d_{2k+1,3+j}}{d_{2k-1,1+j}} \cdot M_{j_z}^{2k-1} \cdot \begin{pmatrix} d_{2k-1,2} & & \\ & \ddots & \\ & & d_{2k-1,k} \end{pmatrix} \cdot A^{2k} \cdot A_{\bar{k}}^{2k+1}$$

$$= \frac{d_{2k+1,3+j}}{d_{2k-1,1+j}} \cdot M_{j_z}^{2k-1} \cdot \begin{pmatrix} d_{2k-1,2} & & \\ & \ddots & \\ & & d_{2k-1,k} \end{pmatrix} \cdot \begin{pmatrix} A^{2k} & \bar{A}^{2k} \end{pmatrix} \cdot A^{2k+1},$$

finishing the proof of Lemma 1 (iv).

**Proof of Lemma 1 (v) :** We assume a curve $z_0(\varepsilon)$ as in (3.9) with $\|G[z_0(\varepsilon)]\| = O(|\varepsilon|^{2k+1})$, implying by (2.5) with $k$ replaced by $2k$

$$\begin{pmatrix} T^{2k}[\bar{z}_{2k}, \ldots, \bar{z}_1] \\ \vdots \\ T^1[\bar{z}_1] \end{pmatrix} = \Gamma^{-2k} \cdot \Delta^{2k}(\bar{z}_{2k-1}, \ldots, \bar{z}_1) \cdot \Gamma^{2k} \cdot \begin{pmatrix} \bar{z}_{2k} \\ \vdots \\ \bar{z}_1 \end{pmatrix} = \begin{pmatrix} 0 \\ \vdots \\ 0 \end{pmatrix}.$$

Then, by Lemma 1 (i), there exists $\bar{n}_i \in N_i, 1 \leq i \leq 2k$, satisfying

$$\Gamma^{2k} \cdot \begin{pmatrix} \bar{z}_{2k} \\ \vdots \\ \bar{z}_1 \end{pmatrix} = M^{4k+1} \cdot \begin{pmatrix} \bar{n}_1 \\ \vdots \\ \bar{n}_{2k} \end{pmatrix}$$

with upper tridiagonal matrix operator $M^{4k+1}$ and $diag[M^{4k+1}] = [I_B \cdots I_B]$. Now, using $\bar{z}_1 = \cdots = \bar{z}_{l-1} = 0$, we obtain $\bar{n}_{2k} = \cdots = \bar{n}_{2k-l+2} = 0$ and $\Gamma_l^{2k} \cdot \bar{z}_l = I_B \cdot \bar{n}_{2k-l+1} \in N_{2k-l+1} \subset N_{k+1}$ due to $l \leq k$.

*Matthias Stiefenhofer*

*University of Applied Sciences*

*87435 Kempten (Germany)*

*matthias.stiefenhofer@hs-kempten.de*